%% file: 0_main_arXiv_May2026.tex
\documentclass[11pt]{article}
\usepackage[letterpaper,hmargin=1.0in,vmargin=1.0in]{geometry}
\usepackage{amsmath,amsthm,amssymb,setspace,sectsty,bbm,graphicx,mathtools}
\usepackage[round]{natbib}
\usepackage[usenames,dvipsnames]{xcolor}
\usepackage[linktocpage=true,pagebackref=true,colorlinks,allcolors=blue,bookmarks=true,bookmarksopen,bookmarksnumbered]{hyperref}

\usepackage{booktabs} 
\usepackage{algorithm,algorithmicx}
\usepackage[noend]{algpseudocode}
\usepackage{subcaption}
\usepackage{bbm}
\usepackage{tabularx}
\usepackage{multirow}
\usepackage{makecell}
\usepackage{xcolor}
\usepackage{bm}

\newtheoremstyle{DStheorem}
{\topsep}
{\topsep}
{\itshape}
{0pt}
{\scshape}
{.}
{ }
{\thmname{#1}\thmnumber{ #2}\thmnote{ (#3)}}
\theoremstyle{DStheorem}
\newtheorem{theorem}{Theorem}[section]
\newtheorem{lemma}[theorem]{Lemma}
\newtheorem{claim}[theorem]{Claim}
\newtheorem{corollary}[theorem]{Corollary}
\newtheorem{observation}[theorem]{Observation}

\newtheorem{proposition}[theorem]{Proposition}
\newtheorem{definition}[theorem]{Definition}

\newtheorem{assumption}{Assumption}[section]
\let\oldproofname=\proofname
\renewcommand{\proofname}{\rm\sc{\oldproofname}}

\usepackage{scalerel,stackengine}
\usepackage{tikz}
\usetikzlibrary{arrows.meta,
            decorations.pathreplacing,
                    calligraphy,
                positioning}

\makeatletter
\def\BState{\State\hskip-\ALG@thistlm}
\makeatother
\def\one{{\mathbbm{1}}}
\def\RR{{\mathbb R}}

\def\NN{{\mathbb N}}

\def\ZZ{{\mathbb Z}}

\def\EE{{\mathbb E}}

\def\PP{{\mathbb P}}

\def\C{{\mathcal C}}

\def\S{{\mathcal S}}
\def\E{{\mathcal E}}

\def\A{{\mathcal A}}

\def\argmax{\text{argmax\,}}
\def\C{{\mathcal C}}

\newcommand{\Ma}{\mathsf{M}}
\def\S{{\mathcal S}}

\newcommand{\OPT}{\mathrm{OPT}}
\newcommand{\alg}{\text{ALG}}
\newcommand{\UB}{\text{UB}}

\newcommand{\revised}[1]{\textcolor{black}{#1}}
\definecolor{mysilver}{RGB}{220,220,220}

\newcommand{\fullystatic}{{\sf(Fully Static)}}
\newcommand{\fullyS}{{\sf(FS)}}
\newcommand{\fullyadaptive}{{\sf(Fully Adaptive)}}
\newcommand{\fullyA}{{\sf(FA)}}
\newcommand{\onesidedstatic}{{\sf(One-sided Static)}}
\newcommand{\Conesidedstatic}{({$\mathcal{C}$}{\sf-One-sided Static})}
\newcommand{\Sonesidedstatic}{({$\mathcal{S}$}{\sf-One-sided Static})}
\newcommand{\onesidedS}{{\sf(OS)}}
\newcommand{\ConesidedS}{({\mathcal{C}}\text{\sf-OS)}}
\newcommand{\SonesidedS}{({\mathcal{S}}\text{\sf-OS)}}
\newcommand{\onesidedadaptive}{{\sf(One-sided Adaptive)}}
\newcommand{\Conesidedadaptive}{({$\mathcal{C}$}\text{\sf-One-sided Adaptive})}
\newcommand{\Sonesidedadaptive}{({$\mathcal{S}$}\text{\sf-One-sided Adaptive})}
\newcommand{\onesidedA}{{\sf(OA)}}
\newcommand{\ConesidedA}{({\mathcal{C}}\text{\sf-OA)}}
\newcommand{\SonesidedA}{({\mathcal{S}}\text{\sf-OA)}}

\usepackage[inline]{enumitem}
\usepackage{hyperref} 
\usepackage{mathtools}
\usepackage{soul}

\allowdisplaybreaks

\sectionfont{\large} \subsectionfont{\normalsize}
\allowdisplaybreaks
\makeindex
\newcommand{\changelocaltocdepth}[1]{%
  \addtocontents{toc}{\protect\setcounter{tocdepth}{#1}}%
  \setcounter{tocdepth}{#1}%
}

\begin{document}
\begin{titlepage}

\title{Two-sided Assortment Optimization: \\ Adaptivity Gaps and Approximation Algorithms}

\author{%
Omar El Housni\thanks{Cornell Tech, Cornell University. Email: {\tt oe46@cornell.edu}.}
\and Ulysse Hennebelle\thanks{Cornell Tech, Cornell University. Email: {\tt uh34@cornell.edu}.}
\and Alfredo Torrico\thanks{ Daniel J. Epstein Department of Industrial and Systems Engineering, University of Southern California. Email: {\tt atorrico@usc.edu}.}}

\date{}
\maketitle

\setcounter{page}{200}
\thispagestyle{empty}

\begin{abstract}
Choice-based matching platforms have recently proliferated thanks to their application in labor markets, dating, accommodation and carpooling.
Due to correlated preferences, these platforms face the challenge of reducing choice congestion among popular options who receive more requests than they can handle, which may lead to suboptimal market outcomes. 
To address this challenge we introduce a two-sided assortment optimization framework under general choice preferences. The goal in this problem is to maximize the expected number of matches by deciding which assortments are displayed to the agents and the order in which they are shown. In this context, we identify several classes of policies that platforms can use in their design. Our goals are: (1) to measure the value that one class of policies has over another one, and (2) to approximately solve the optimization problem itself for a given class. For (1), we define the adaptivity gap as the worst-case ratio between the optimal values of two different policy classes. First, we show that the  gap between the class of policies that statically show assortments to one-side first and the class of policies that adaptively show assortments to one-side first is exactly $e/(e-1)$. Second,  we show that the gap between the latter class of policies and the fully adaptive class of policies that show assortments to agents one by one is exactly $2$. We also note that the worst policies are those who simultaneously show assortments to all the agents, in fact, we show that their adaptivity gap even with respect to one-sided static policies can be arbitrarily large. 
For (2), we first design a polynomial time algorithm that achieves a $1/4$ approximation factor within the class of policies that adaptively show assortments to agents one by one. Furthermore, when agents' preferences are governed by multinomial-logit models, we show that a 0.067 approximation factor can be obtained within the class of policies that show assortments to all agents at once.  We further generalize our results to constrained assortment settings, where we impose an upper bound on the size of the displayed assortments. Finally, we present a computational study to evaluate the empirical performance of our theoretical guarantees.
\end{abstract}

    
\end{titlepage}


\newpage
\pagestyle{plain}
\setcounter{page}{1}

\input{1_introduction}

\input{1-1_contributions}

\input{1-2_related_literature}
%
\input{2_problem_formulation}


\input{3_adaptivity_gaps}

\input{4_approx_fully_adaptive}

\input{5_approx_fully_static}

\input{6_extensions_cardinality}

\input{7_numerical_experiments}

\input{8_conclusions}



\addcontentsline{toc}{section}{Bibliography}
\bibliographystyle{plainnat}
\bibliography{bibliography_assortment} 

\changelocaltocdepth{1} 
\appendix
\input{999_appendix}
\end{document}

%% file: 1_introduction.tex
\section{Introduction}\label{sec:introduction}
Choice-based matching platforms have transformed the way we move around cities, how we connect and interact with others, and how we outsource tasks and search for jobs. Classic examples include freelancing platforms like Taskrabbit and Upwork, dating platforms such as Bumble and Tinder, lodging apps like Airbnb and Vrbo, and {carpooling} companies such as Blablacar.  These markets consists of two different sides (e.g., travelers and hosts) who are looking to match, i.e., to initiate an interaction or transaction. More importantly, the platform does not centrally dictate the pairs, instead, the matching process is choice-based in that users select profiles from a subset of displayed alternatives and a match is realized when two users, on different sides, mutually select each other. The platform, then, faces \emph{design} and \emph{efficiency} challenges since it serves as a mediator that controls how users interact and what they see, as we explain next:
\begin{itemize}
\item[-] {\bf Design.} Some crucial aspects in the design of a choice-based platform are: the number of selections users make and the direction of interactions. First, certain markets are designed such that users  (roughly speaking) select at most one option among those being displayed, e.g., \revised{users of the Airbnb platform} generally submit a single request. In other markets, such as dating apps, users select multiple profiles in one session. Second, the direction of interaction refers to which side initiates the matching process. Some markets allow for \emph{two-way interactions} where both sides can search the market and initiate an interaction, for example, in Upwork, both someone looking for a job as a freelancer and someone outsourcing a project can reach out first. Other platforms, due to privacy concerns or specific application needs, are designed with \emph{one-way interactions} where only one side (established in advance) initiates the request, e.g., travelers reach out to hosts in Airbnb, while the hosts respond to requests.
\item[-] {\bf Efficiency.} These matching markets account not only for the way users interact, but also for their preferences, which in turn affect what they see. In general, each user on the platform has preferences over the opposite side, for instance, \revised{freelancers with vast experience are likely to be preferred over those who appear less trained}. These preferences may lead to  \emph{choice congestion} which occurs when the most popular users receive more requests than they can handle, resulting in market inefficiencies. For each platform design, a natural tension \revised{arises} between (i) displaying relevant options to keep users engage, and (ii) reducing choice congestion among popular alternatives. Carefully balancing these two may translate {into} improved revenue and experience.
\end{itemize}
\vspace{0.5em}
In this work, we aim to provide a toolkit to address design and efficiency challenges faced by choice-based matching platforms where users make a single selection\footnote{We do not focus in this paper on markets where users make multiple selections simultaneously, such as dating apps.}, e.g., Airbnb and Upwork. 
To achieve this, we introduce a broad framework for two-sided assortment optimization: A topic that has recently gained attention as an algorithmic tool to reduce choice congestion \citep{ashlagi_etal19,torrico21,aouad2020online,rios2022improving,rios2023platform,ahmed2022parameterized}. 
Our framework allows us to make a comprehensive comparison between platform designs and develop algorithmic solutions for them. Informally, we model the platform as a bipartite graph where customers and suppliers (both collectively called agents) are looking to match with each other. We model the agents' preferences using \emph{discrete choice models} to capture substitution effects. This means that each agent probabilistically picks either one option from the displayed \emph{assortment} of alternatives or leaves the market. The platform aims to maximize the expected total number of matches between both sides. The matching process is controlled by a policy which, in each step, \revised{\emph{processes}} a group of agents, designs their assortments and observes their choices; the policy continues with the remaining individuals until no more are left. To put it simply, the policy is an abstraction of the platform design that controls how users interact and what they see. In this broad context, we focus on four main policy classes which capture not only one-way and two-way interactions, but also different levels of \emph{adaptivity}. This last feature is particularly relevant to distinguish platforms with limited resources which cannot constantly re-compute assortments versus those that need to adapt to agents' previous choices (e.g., Airbnb updates the listings once a host accepts a request).

Our main contribution is twofold: First, to address the design challenge, we consider the concept of \emph{adaptivity gap} as, \revised{roughly speaking, the largest possible ratio between the optimal expected number of matches achieved by two different policy classes, across all instances}. This allows us to effectively measure the benefit of adopting a platform design with one-way interactions relative to two-way interactions, as well as the gains of adaptivity. 
Second, to address efficiency, we design algorithms that, for a fixed class of policies determined by the platform’s design choices, solve the corresponding optimization problem of deciding which assortments to display and in what order.

%% file: 1-1_contributions.tex

\subsection{Our Results and Contributions}\label{sec:contributions}
In this work, we introduce a two-sided assortment optimization framework on a two-sided matching platform {comprised of} customers and suppliers with general choice preferences. \revised{In this context, a customer and a supplier get matched when, given the assortments displayed by the platform, they mutually select each other}. 
Within this framework, we identify four main policy classes \revised{that are characterized by how they process the set of agents in the platform, i.e., (i) which subset of agents will be shown an assortment and (ii) which assortment is presented to each of those agents}. First, we define \emph{fully static} policies, which simultaneously present assortments to all agents on the platform.
At the other extreme, we define \emph{fully adaptive} policies, where agents are processed one at a time, and the platform observes each agent’s choice before processing the next.
In between these two extremes, we define \emph{one-sided} policies, where all agents on one side must be processed before any on the other side. We distinguish between two types of one-sided policies: \emph{one-sided static} policies, where the assortments are simultaneously shown to the first side,
and \emph{one-sided adaptive} policies, where the platform can adapt the assortments for the first side based on the choices of previously processed agents on that same side.
We study the adaptivity gaps between these policy classes and design algorithms to solve the corresponding optimization problems.
We refer to Figure~\ref{fig:adaptivity_gap} for a summary of our main contributions, and in what follows, we present the details of each of our results.

\begin{figure}[htpb]
\centering
\scalebox{1}{
 \begin{tikzpicture}[
BC/.style args = {#1/#2/#3}{
        decorate,
        decoration={calligraphic brace, amplitude=6pt,
        pre =moveto, pre  length=1pt,
        post=moveto, post length=1pt,
        raise=#1,
              #2},
         thick,
        pen colour={#3}, text=#3
        },              ]
\draw[thick,-] (0,0) -- (3,0)  {};
\draw[thick,dotted] (3,0) -- (4,0)  {};
\draw[thick,-] (4,0) -- (15,0)  {};
\draw[-] (1,-0.2) -- node[below=2ex] {\footnotesize\fullystatic} (1,0.2);
\draw[-] (6,-0.2) -- node[above=2ex] {\footnotesize\onesidedstatic} (6,0.2);
\draw[-] (10,-0.2) -- node[below=2ex] {\footnotesize\onesidedadaptive} (10,0.2);
\draw[-] (14,-0.2) -- node[above=2ex] {\footnotesize\fullyadaptive} (14,0.2);
 \draw[arrows = {Stealth[scale=1.1]-|}] (1,0.5) -- node[above=2ex] {\footnotesize Theorem~\ref{thm:fullystatic_mnl_approx}: $0.067$ (for MNL)} (1,1.3);
 \draw[arrows = {Stealth[scale=1.1]-|}] (14,1) -- node[above=2ex] {\footnotesize Theorem~\ref{coro:approx_fully_adaptive}: $1/4$} (14,1.8);
 \draw[arrows = {Stealth[scale=1.1]-|}] (10,0.5) -- node[above=2ex] {\footnotesize $\begin{matrix}
 \text{Theorem~\ref{thm:approxfactor_greedy}:} & 1/2\\
\text{Corollary~\ref{coro:approx_onesided_adaptive}:} & 1/2-\epsilon
 \end{matrix}$
 } (10,1.3);
 \draw[arrows = {Stealth[scale=1]-|}] (6,0.9) -- node[above=2ex] {\footnotesize $1-1/e$, \citep{torrico21}} (6,1.6);
\draw [BC=1mm/mirror/black] 
    (1,-1.1) -- node[below=2ex] {\footnotesize Proposition~\ref{prop:gap_static_zero}: $\mathcal{O}(n)$} (6,-1.1);
\draw [BC=1mm/mirror/black] 
    (6,-1.1) -- node[below=2ex] {\footnotesize Theorem~\ref{thm:gap_onesided_adaptive}: $e/(e-1)$} (10,-1.1);
\draw[BC=1mm/mirror/black]    
(10,-1.1) -- node[below=2ex] {\footnotesize Theorem~\ref{thm:gap_fully_adaptive}: $2$} (14,-1.1);
\draw[BC=1mm/mirror/black]    
(6,-2.1) -- node[below=2ex] {\footnotesize Theorem~\ref{thm:gap_os_fa}: $2e/(e-1)$} (14,-2.1);
 \end{tikzpicture}}
 \caption{Main results for the unconstrained setting. Each adaptivity gap (presented below a bracket) corresponds to the ratio between the value on the right over the value on the left, for example, the adaptivity gap between \onesidedstatic\ and \onesidedadaptive\ is $e/(e-1)$. On  the top, we present our  guarantees for each class of problems (all guarantees are obtained in polynomial-time except  the pseudo-polynomial-time sampling method in Corollary~\ref{coro:approx_onesided_adaptive}).}

\label{fig:adaptivity_gap}
\end{figure}

\vspace{2mm}
\noindent
{\bf Adaptivity Gaps.} 
 Our main results on adaptivity gaps are given in Section~\ref{sec:adaptivity_gaps}, and  summarized  as  follows:

\begin{enumerate}
    \item[(a)] {\em Adaptivity gap between fully-static~and one-sided static policies:} In Proposition~\ref{prop:gap_static_zero}, we provide an example {showing} 
that the gap between these two policy classes can be arbitrarily large even under a very simple choice model. 
In practical terms, a matching platform would significantly benefit by allowing one side of the market to choose first, as this enables the platform to collect more information about agents' preferences.

\item[(b)] {\em Adaptivity gap between one-sided static and one-sided adaptive policies:} In Theorem~\ref{thm:gap_onesided_adaptive}, we show that the gap between these classes is exactly $e/(e-1)$, and this result is tight. The proof is structured in two parts. Lemma~\ref{lemma:LB_onesided_adaptive} establishes the $e/(e-1)$ upper bound using a novel \emph{linear programming relaxation} that captures the relationship between the choices made by the initiating side and the responses of the other side. This relaxation serves as an upper bound on the performance of any one-sided adaptive policy, as any such policy can be mapped to a feasible solution of the linear program with the same objective value.  
We note that the relaxation encodes two extreme probability distributions: the independent distribution (also known as product distribution), which corresponds to one-sided static policies, and the worst-case distribution, which bounds the optimal one-sided adaptive policy. Under mild assumptions on agents' preference structures, our proof leverages the concept of \emph{correlation gap}~\citep{agrawal2010correlation} to evaluate the objective values under these distinct distributions.  In Lemma~\ref{lemma:UB_onesided_adaptive}, we provide a non-trivial instance for which the adaptivity gap is asymptotically close to $e/(e-1)$. 

\item[(c)] {\em Adaptivity gap between one-sided adaptive and fully adaptive policies:} In Theorem~\ref{thm:gap_fully_adaptive}, we show that the gap between these policy classes is exactly $2$, and this result is tight. The proof again consists of two parts. Lemma~\ref{lemma:LB_fully_adaptive} establishes the upper bound by coupling an optimal fully adaptive policy with two feasible one-sided adaptive policies:  
(1) a policy that mimics the actions of the optimal policy when processing suppliers, while skipping the steps involving customers (customers are processed at the end, receiving assortments as determined by the optimal policy), and  
(2) the opposite policy that mimics the actions for customers and skips suppliers.  
We show that the total number of matches produced by these two one-sided adaptive policies is at least the number of matches achieved by the optimal fully adaptive policy. Lemma~\ref{lemma:UB_fully_adaptive} provides an instance where the adaptivity gap approaches 2.

\item[(d)] {\em Adaptivity gap between one-sided static and fully adaptive policies:} From results (b) and (c), it follows that the gap between these two policy classes is at most $2e/(e-1)$. In Theorem~\ref{thm:gap_os_fa}, we construct an instance where the adaptivity gap is asymptotically close to $2e/(e-1)$.
\end{enumerate}

\vspace{2mm}

\noindent
{\bf Approximation Algorithms.}  Our main algorithmic results for the unconstrained setting are given in Sections~\ref{sec:approx_fully_adaptive} and~\ref{sec:approx_fully_static}, and  summarized  as  follows:

\begin{itemize}
    \item [(a)] {\em Approximation algorithm for fully adaptive policies.}  
    We present a randomized polynomial-time algorithm that achieves a $1/4$-approximation for the fully adaptive policy class; see Theorem~\ref{coro:approx_fully_adaptive}. This result is obtained by leveraging our approximation guarantees for one-sided adaptive policies when the first side being processed is defined in advance.
    Specifically, we consider separately the class of one-sided adaptive policies when customers are processed first and, analogously, for suppliers.
    In Theorem~\ref{thm:approxfactor_greedy}, we design a greedy algorithm (\emph{Arbitrary Order Greedy}) that achieves a $1/2$-approximation with respect to the optimal policy within each of these classes.  The algorithm works as follows: assuming customers are processed first, the method sets an arbitrary order over the customer set. For each “arriving” customer, it displays an assortment of suppliers that maximizes the marginal probability that the customer is selected by a supplier, weighted by her probabilistic choices. After all customers have made their selections, each supplier is shown the subset of customers who initially chose them (i.e., their backlog). This algorithm builds on a greedy routine introduced in \citep{aouad2020online}, but our analysis differs in that it does not rely on the arrival order. Instead, we provide a new analysis using a primal-dual argument and the order-independent LP relaxation from Section~\ref{sec:onesided_adaptive}.  To obtain Theorem~\ref{coro:approx_fully_adaptive}
    we run the Arbitrary Order Greedy algorithm twice: once assuming customers initiate, and once assuming suppliers initiate. We then show that randomly selecting one of the two runs with equal probability achieves a $1/4$-approximation in expectation relative to the optimal fully adaptive policy.
Additionally to the results above, if the platform must implement a one-sided adaptive policy but has not predetermined which side initiates, we provide a pseudo-polynomial time sampling method (Corollary~\ref{coro:approx_onesided_adaptive}) that, with high probability, identifies the better initiating side and runs the greedy algorithm accordingly. This method achieves a $(1/2 - \epsilon)$-approximation against the best one-sided adaptive policy across both possible starting sides.

   \item [(b)] {\em Approximation algorithm for fully-static policies:} 
In Section~\ref{sec:approx_fully_static}, we focus on the class of fully static policies under the assumption that customer and supplier choice preferences follow multinomial-logit (MNL) models, with each agent having a distinct MNL model. We first show that the corresponding optimization problem is strongly NP-hard, even when all customers share the same MNL model and all suppliers share the same MNL model. Our main contribution in this section is the design of an algorithm that achieves a $0.067$-approximation factor (Theorem~\ref{thm:fullystatic_mnl_approx}). To establish this result, we categorize agents into groups based on whether their preference values are high or low. 
For instances where both customers and suppliers have low-valued preferences, we propose a linear programming relaxation and obtain a fractional solution, which we convert into a feasible one using independent randomized rounding with only a constant-factor loss. In scenarios where one group of agents (either customers or suppliers) has high preference values, we bound the problem with a submodular optimization relaxation, from which we can also extract a feasible solution with a constant-factor loss. The final $0.067$-approximation guarantee is derived by combining the guarantees across the three regimes and leveraging the subadditivity of the objective function.

\end{itemize}

\vspace{0.5em}
\noindent
{\bf Extensions to Constrained Assortment Models.} 
In Section~\ref{sec:cardinality_constraints}, we extend our framework to cardinality-constrained settings, where each agent can be shown at most a fixed number of alternatives. We distinguish between the \emph{two-way constrained} setting, where both customers and suppliers are subject to constraints, and the \emph{one-way constrained} setting, where constraints apply only to the side being processed first in one-sided policies. The two-way constrained model introduces new algorithmic and analytical challenges not present in the unconstrained case: constrained demand functions may no longer be submodular or efficiently computable, and previous LP relaxations must be adapted to handle coupled constraints on both sides. We address these challenges by designing new techniques combining LP relaxations with contention resolution schemes.
Our first main result proves that the adaptivity gap between one-sided static and one-sided adaptive policies remains bounded by $(e/(e-1))^2$ in the two-way constrained setting, and improves to exactly $e/(e-1)$ for MNL models and for all one-way constrained instances (Theorem~\ref{thm:cardinality_gap_os_oa}). On the algorithmic side,
our main result (Theorem~\ref{thm:approx_oa_constrained_general}) shows that, under the two-way constrained setting and assuming weak substitutability, there exists a polynomial-time algorithm achieving a $(1 - 1/e)^3$-approximation for both one-sided static and one-sided adaptive policies, and a $\frac{1}{2}(1 - 1/e)^3$-approximation for fully adaptive policies. These guarantees improve to $1/2$ and $1/4$, respectively, under MNL choice models or in the one-way constrained setting, where  submodularity and efficient evaluation are preserved (Corollary~\ref{coro:improved_guarantees_constrained}).  Finally, we extend our approximation result for the fully static policy to the constrained setting, showing that a $0.067$-approximation is achievable under MNL preferences (Theorem~\ref{thm:fullystatic_mnl_approx_cardinality}).

\vspace{0.5em}
\noindent
{\bf Empirical Analysis.} In Section~\ref{sec:experiments}, we empirically evaluate the performance of our algorithms and analyze adaptivity gaps across different policy classes using synthetic instances. We compare our methods against both optimal policies (when tractable) and other tractable upper bound  benchmarks. The results show that our algorithms consistently outperform their theoretical worst-case guarantees and remain competitive even when compared to stronger baselines, highlighting their practical effectiveness. \revised{In Appendix~\ref{appendix:sensitivity}, we focus on the empirical performance of the main algorithms on larger instances. We study the sensitivity of their performances with respect to two main aspects that influence the market's outcome: (i) the value of the outside option, and (ii) the agents' heterogeneity.}

\vspace{0.5em}
\noindent\revised{\textbf{Practical Impact.} The model  introduced in this work can be interpreted as a broad abstraction that captures several types of choice-based matching platforms. To the best of our knowledge, our model is the first to include different fundamental features of two-sided markets along with discrete choice modeling: (i) direction of interactions (one-way and two-way); (ii) levels of adaptivity (non-adaptive to adaptive, and their respective one-sided and two-sided variants), and (iii) feasibility constraints. First, the direction of interactions is mainly determined by the specific platform application, for instance, outsourcing and carpooling apps my benefit by either one-way or two-way interactions. Second, depending on the platform's resources and design, certain types of policies might be more appropriate: On the one hand, non-adaptive policies can be an attractive alternative for decision-makers with limited resources as these methods perform relatively well in certain markets. On the other hand, adaptive policies are theoretically and empirically strong, however, re-optimization might be quite limiting for certain applications. Finally, cardinality constraints are useful in overcrowded markets or when dealing with recommendation fatigue, 
pre-selecting a subset of users to be displayed in the first positions can address these issues while optimizing the overall market's outcome. 
Given this, our work reveals not only the theoretical (relative) differences when incorporating these platform features but also provides algorithmic solutions with strong provable guarantees and empirical performance. Finally, we remark that our model can be adapted to incorporate other market features, some of which we discuss in Section~\ref{sec:conclusions} and leave as interesting future directions.
}

%% file: 1-2_related_literature.tex
\color{black}
\subsection{Related Literature}\label{sec:related_literature}
In the following, we present the most relevant streams of literature. \paragraph{One-sided Assortment Optimization.} Traditional assortment optimization has been extensively studied in the one-sided setting where the goal is to select the subset of products that maximizes the expected revenue achieved by the purchase of a single customer. The single-customer assortment optimization problem was introduced in \citet{talluri04}, and, since then, numerous offline and dynamic models, in the unconstrained and constrained settings, have been considered. 
Notably, \citet{talluri04} show that the unconstrained offline problem can be solved in polynomial time when the customer's preferences are governed by the \emph{multinomial choice model}. 
The offline problem has been studied under various choice models such as the Logit-based~\citep{talluri04,rusmevich14,davis14}, the Rank-based~\citep{aouad18} the Markov chain-based \citep{blanchet16} and general choice models~\citep{berbeglia_joret20}. Several of these optimization problems are NP-hard e.g., \citep{rusmevich14,aouad18,desir2022capacitated}, and  a significant portion of research in this area is dedicated to developing approximation algorithms.
The constrained assortment optimization problem has also been studied under several choice models, see e.g. \citep{rusmevich_etal10,sumida2021revenue,desir_etal15,el2023joint,barre2023assortment}. For a survey on assortment optimization, we refer to \citep{kok_etal08}. The online version of the problem has also been studied under several settings and arrival models, see e.g.~\citep{golrezaei2014real,ma2020algorithms,goyal2025asymptotically,feng2022near,gong2022online}. 

\paragraph{Efficiency and Congestion in Two-sided Markets.}
There is an extensive literature on improving efficiency and reducing market congestion starting with the seminal work of \citet{rochet_tirole03}. Since then, the literature has studied several market interventions to improve the final outcome such as matchmaking strategies~\citep{shi2022strategy}, signaling competition levels~\citep{besbes2023signaling}, {notifications~\citep{manshadi22},}   market recommendations~\citep{horton17} and ranking~\citep{fradkin15}, limiting the visibility \citep{halaburda_etal16,arnosti_etal18} and choice~\citep{immorlica2021designing}, controlling who initiates contact \citep{kanoria2021facilitating,ashlagi_etal2020clearing} \revised{and stability constraints~\citep{camelo2025stable}}. We emphasize that most of this literature focuses on analyzing specific design features of the market rather than taking an algorithmic approach as in this work. 

\paragraph{Two-sided Assortment Optimization.} 
One of the first works in this stream of literature was \citep{ashlagi_etal19}. According to the notation we use in this work, \citet{ashlagi_etal19} consider a one-sided static problem where customers initiate the interaction. Moreover, customers follow the same multinomial-logit choice model and suppliers follow a ``uniform'' model. Under this specific setting, they show the problem is strongly NP-hard and give a $10^{-5}$-approximation algorithm. 
Their results were later improved by \citet{torrico21} who show a tight $1-1/e$ approximation factor under a general choice model setting by using an adaptation of the continuous greedy algorithm. \citet{ahmed2022parameterized} consider a fully-static two-sided assortment optimization problem and show parameterized guarantees. We significantly improve upon their results by providing the first constant-factor approximation guarantee for the fully static problem. We emphasize that our framework is a generalization of the settings studied in \citep{ashlagi_etal19,torrico21,ahmed2022parameterized} as we allow a broader space of policies ranging from fully-static to fully-adaptive.

\revised{Motivated by labor markets,} \citet{aouad2020online} introduce an online version of the one-sided static model considered by~\citep{ashlagi_etal19}. Specifically, customers arrive one by one in an online fashion and suppliers' availability also varies from stage to stage. Both arrivals and availability may be adversarially chosen. The authors show that under mild assumptions, an online greedy method (analog of the Arbitrary Order Greedy described in Section~\ref{sec:approx_fully_adaptive})  attains a $1/2$ competitive ratio, which is the best possible. Also, \citet{aouad2020online} show improved guarantees for the i.i.d. arrival settings and specific customers' choice models.
Their online model is different than ours since we allow the platform to make adaptive decisions on not only the assortments, but also on the order that the agents are processed. In this sense, optimal one-sided static and one-sided adaptive policies serve as theoretical bounds for the online setting.
More importantly, their analysis of the greedy algorithm cannot be applied to our setting as they use an order-dependent benchmark. Instead, our primal-dual argument is more general as it is against an order-independent linear programming relaxation. 

\revised{Recently, both \citet{nissim2024revenue} and \citet{ahmadnejadsaein2025adaptive} have studied the two-sided assortment problem with general edge-dependent revenues under the MNL choice model, providing approximations for several policy classes. Their works show that the model can be adapted to more realistic settings while retaining strong theoretical guarantees.}
\revised{From a more practical viewpoint, assortment optimization has also recently been considered for a patient-provider matching problem~\citep{raman2025assortment} where the goal is to reduce choice congestion and improve matches quality.}

Finally, motivated by dating platforms, \citet{rios2022improving,rios2023platform} consider a dynamic setting where the goal of the platform is to show subset of profiles to every agent in each stage with the objective of maximizing the expected total number of matches. Agents' preferences are determined by: (i) pair-dependent probabilities which are independent of the displayed assortments, (ii) agents can like/dislike multiple profiles in each stage and (iii) assortments are restricted by cardinality constraints. 
\citet{rios2022improving} show that this problem is NP-hard and focus mostly on an empirical analysis. Later, \citet{rios2023platform} provide several approximation guarantees under different platform features. 
The most crucial differences with our work are that, in their dynamic model, agents make independent choices and select multiple profiles, while in our setting there is a substitution effect captured by general choice preferences.

\paragraph{Stochastic Matching (Probe-commit Model).}
Motivated by online dating and kidney exchange, \cite{Chen09} introduce the stochastic matching problem in the probe-commit model, also known as stochastic probing. In this problem, every pair of nodes of a given general graph has a probability of being compatible or not. The matchmaker can sequentially probe edges to check the compatibility and, if compatible, then the edge must be included in the final set. The goal is to form a matching of maximum cardinality. They provide the first performance guarantees for different query models and patience constraints (bounds on the number of edges incident to nodes that can be probed). Since then,  several algorithmic improvements have been shown for the (un)weighted versions, levels of patience, and online variants, see e.g.~\citep{adamczyk2011improved, costello2012stochastic,bansal2012lp,gamlath2019beating,brubach2021improved}. For a more detailed discussion, we refer to the recent work of~\citet{derakhshan2023beating}. 
First, we highlight two important aspects that distinguish the stochastic bipartite matching model with our framework: (1) in their model, the decision-maker queries edges, instead our policies must query two nodes to realize an edge; (2) their probabilities are edge-dependent, while in our setting the probabilities depend on the assortments. Finally, we note that our setting generalizes the bipartite case with unit-patience constraints.

\paragraph{Adaptivity Gaps in Stochastic Optimization.} 
The literature of stochastic optimization has extensively studied the concept of adaptivity for different settings and applications such as knapsack constraints~\citep{dean2008approximating,bhalgat2011improved,ma2018improvements}, packing integer program~\citep{dean2005adaptivity},  TSP~\citep{jiang2019algorithms}, covering constraints~\citep{agarwal2019stochastic,ghuge2021power} and hiring~\citep{purohit2019hiring}. Part of this literature closely relates to the stochastic matching probing model described above which was generalized to submodular objectives in~\citep{adamczyk2016submodular}. For monotone submodular functions, we highlight the work of~\citet{asadpour2016maximizing} who show a $e/(e-1)$-gap between non-adaptive and adaptive policies when the elements being probed must satisfy a matroid constraint. 
For the same setting, \citet{hellerstein2015discrete} provides parametric bounds, and the generalization to other constraints is studied in~\citep{gupta2016algorithms,gupta2017adaptivity,bradac2019near}. 
As in the probe-commit matching model above, these settings assume that each probability is element-dependent, while in our setting, instead, each distribution depends on the assortment decided by the platform.
%

%% file: 2_problem_formulation.tex
\color{black}
\section{Model}\label{sec:model}

\subsection{Problem Formulation}\label{sec:problem_formulation}
{In this section, we introduce our general framework for the \emph{unconstrained two-sided assortment optimization problem}.}
Consider a two-sided platform represented by a bipartite graph where, on  one side, we have a set $\C$ with $n$ customers and, on the opposite side, we have a set $\S$ with $m$ suppliers. The entire set of agents in the platform will be denoted as $\A=\C\cup\S$. {
In the following, we formalize the agents' preferences and the platform's objective.} 
\paragraph{Agents' Preferences.}
The preferences of each customer $i\in \C$ are captured by a discrete choice model $\phi_i: (\S\cup\{0\})\times 2^{\S}\to[0,1]$ where $\{0\}$ represents the outside option (e.g., a different platform). For any assortment of suppliers $S\subseteq\S$ presented to customer $i$, the value $\phi_i(j,S)$ corresponds to the probability that customer $i$ chooses option $j\in\S\cup\{0\}$, where these values are such that $\sum_{j\in S\cup\{0\}}\phi_i(j,S)=1$ and $\phi_i(j,S)=0$ for any $j\notin S\cup\{0\}$. In other words, for a given assortment, the customer probabilistically chooses either the outside option or one alternative from the assortment.
We denote the demand function of customer $i$ as $f_i:2^\S\to[0,1]$ that, for any subset $S\subseteq\S$, returns $f_i(S)=\sum_{j\in S}\phi_i(j,S)$, i.e., the probability that customer $i$ chooses someone from $S$. Similarly, the preferences of each supplier $j\in\S$ are determined by a choice model $\phi_j:(\C\cup\{0\})\times 2^{\C}\to[0,1]$ whose demand function will be denoted as $f_j$. 
We emphasize that agents can have heterogeneous preferences, i.e., we do not assume the same choice model for each agent. Throughout this work, we consider the following standard assumptions which have been previously considered in the two-sided assortment optimization literature, e.g., \citep{aouad2020online,torrico21}. 
\begin{assumption}\label{assumption:monotone_submodular}
We consider the following assumptions:
\begin{enumerate}
  \item {\bf Single-agent Assortment Optimization Oracle.} {For any customer $i\in\C$ and non-negative values $\{\theta_{i,j}\}_{j\in \S}$, there is an oracle that solves in polynomial time the \emph{unconstrained single-agent assortment optimization problem} $\max\Big\{\sum_{j\in S}\theta_{i,j}\cdot\phi_i(j,S): \ S\subseteq\S\Big\}$. The assumption is analogous for suppliers.}
 \item {\bf Monotone Submodular Demand Functions.} 
The demand function \revised{$f_a$ of every agent $a\in\C\cup\S$} is a monotone submodular set function.\footnote{A non-negative set function $f:2^\E\to\RR_+$ defined over a set of elements $\E$ is \emph{monotone} if for every pair of subsets $E\subseteq E'\subseteq \E$, we have $f(E)\leq f(E')$. Function $f$ is \emph{submodular} if for every $e\in\E$ and $E\subseteq E'\subseteq \E\setminus\{e\}$, we have $f(E\cup\{e\})-f(E)\geq f(E'\cup\{e\}) - f(E')$.}
\end{enumerate}
\end{assumption}

{We remark that these assumptions are tailored to the unconstrained setting, i.e., the platform is not limited to any assortment structure. One could similarly formalize Assumption~\ref{assumption:monotone_submodular} for the constrained setting which we discuss in Section~\ref{sec:cardinality_constraints}. We also emphasize that Assumption~\ref{assumption:monotone_submodular} is not restrictive. For the first point in Assumption~\ref{assumption:monotone_submodular}, there are polynomial-time algorithms for several choice models, for instance, the logit-based models~\citep{talluri04,sumida2021revenue,davis14} and the Markov-chain model~\citep{blanchet16}. Moreover, this oracle assumption can be easily extended to an approximate oracle; for more details, we refer to Remark~\ref{obs:approximate_oracle} in Section~\ref{sec:approx_fully_adaptive}. For the second point Assumption~\ref{assumption:monotone_submodular}, we note that the demand functions of the wide class of random utility choice models are monotone and submodular~\citep{berbeglia_joret20}.\footnote{The converse is not true: Not every monotone submodular function induces a random utility choice model; for a counterexample, see~\citep{berbeglia_joret20}.}}

\paragraph{The Platform's Objective.}
The platform aims to clear the market by showing a single assortment to each agent (possibly in an adaptive way) such that the {expected total} number of matches is maximized, {where a match is realized when a customer and a supplier mutually select each other. The} assortments and the order in which they are displayed to the agents are determined by a \emph{policy}. 
{Formally,} a feasible policy $\pi$ consists of a sequence of actions such that, in each step, $\pi$ decides to \emph{process} a subset of the agents, specifically: (1) it selects a subset of agents $A\subseteq\A$ and (2) it presents an assortment to each agent $a\in A$. We assume that there are no constraints on the assortments, i.e., any subset of suppliers $S\subseteq\S$ can be displayed to customers and any subset of customers $C\subseteq\C$ can be shown to suppliers.
Once the assortments are shown to the agents in $A$, the policy observes their choices\footnote{The agents' choices are \revised{made} independently, irrevocably and simultaneously.} and {continues} to the next step with the remaining agents $\A\setminus A$, if any. The policy stops when no agents remain to be processed. %
Denote by $\Pi$ the space of all feasible policies, and \revised{any subset of policies $\Pi'\subseteq\Pi$ (possibly satisfying a common property) will be referred to as a policy class}. For a policy $\pi\in \Pi$, let $\Ma_\pi$ be the random variable that indicates the total number of matches obtained under $\pi$. \revised{Given a policy class $\Pi'\subseteq\Pi$,} the platform is interested in solving the following unconstrained two-sided assortment optimization problem:
\begin{equation}\label{eq:general_problem}
\revised{\OPT_{\Pi'}:=} \max\Big\{\EE[\Ma_\pi]: \ \revised{\pi\in\Pi'}\Big\}
\end{equation}
where the expectation is taken over the agents' choices and the possible randomized actions made by the policy.
{As we mentioned earlier, from a practical standpoint, a policy can be interpreted as an abstraction of a platform design. Therefore, depending on the design, the platform may not have access to the whole space of policies $\Pi$ when attempting to solve Problem~\eqref{eq:general_problem}.}
{
Given the vast landscape of policy classes and platform applications, our first goal is to measure how valuable it is for matching platforms to adopt a certain class of policies relative to others. For this, we consider the following metric:
\begin{definition}[Adaptivity Gap]\label{def:adaptivity_gap}
For a given pair of policy classes $\Pi', \Pi''\subseteq \Pi$ such that $\Pi'\subseteq\Pi''$\footnote{$\Pi'\subseteq\Pi''$ means that any policy in $\Pi'$ can be simulated by some policy in $\Pi''$.} we define the adaptivity gap between them as:
\[
\textsc{Gap}(\Pi',\Pi'')\ =\ \max_{\text{all instances } I}\left\{ \frac{\OPT_{\Pi''}(I)}{\OPT_{\Pi'}(I)} \right\}, \ 
\]
where an instance $I$ is determined by $\C$, $\S$ and the choice models $\phi$. In the expression above, $\OPT_{\Pi'}(I)$ indicates the optimal value when restricted to policies in $\Pi'$ and its dependence on the instance $I$. To minimize notation, in the remainder of the paper, we do not write this dependence when it is understood from context.
\end{definition}
Note that since $\Pi'\subseteq\Pi''$, then we have that $\textsc{Gap}(\Pi',\Pi'')\geq 1$, where a gap equal to 1 means that  the amount of information used by optimal policies in $\Pi''$ does not result in a better matching outcome relative to optimal policies in $\Pi'$. In other words, the platform does not benefit from implementing an optimal policy in $\Pi''$ relative to $\Pi'$.}

{Our second goal in this work is to understand the technical and algorithmic limitations within a specific platform design. To study this, we aim to develop efficient policies with provable \emph{approximation guarantees}. Formally, given a policy class $\Pi'\subseteq\Pi$, we say that a feasible policy $\pi\in\Pi'$ guarantees an $\alpha$-approximation factor with $\alpha\in[0,1]$ if {it runs in polynomial-time} and $\EE[\Ma_{\pi}]\geq \alpha\cdot \OPT_{\Pi'}$ where $\OPT_{\Pi'}$ is the optimal objective value restricted to policies in $\Pi'$. Note that adaptivity gaps and approximation guarantees are closely related since the approximation guarantee achieved within one policy class will carry over to a larger class, up to the measured gap. 
}

\subsection{{Our Main Policy Classes}}\label{sec:policy_classification}
{
To capture different platform applications and designs, we classify policies based on two main aspects: (1) the extent to which past information (agents' choices) can be used to design assortments in future steps; and (2) whether the platform allows one-way or two-way interactions, that is, whether one side is fully processed before the other or agents on both sides can be processed in an interleaved manner.}

{First, let us describe our main policy classes with two-way interactions, i.e., where the policy is not required to process an entire side first before processing someone from the opposite side. Practically speaking, these classes corresponds to abstractions of matching platforms that allow any side to initiate an interaction. For example, labor markets
allow suppliers (someone looking to complete a task or project) to reach out customers (freelancers), or vice versa.}

\noindent{\bf Fully static.} We define $\Pi_{\fullyS}$ as the class of \emph{fully static} policies that process all the agents in $\A$ at the same time. Note that these policies are not able to adapt to the agents' choices since all the assortments are simultaneously displayed. {These policies may be suitable for platforms with two-way interactions that compute a family of assortments in advance and are not updated over a given time period (e.g., daily) due to limited resources. We define \fullystatic~as Problem~\eqref{eq:general_problem} restricted to policies in $\Pi_{\fullyS}$ and by $\OPT_{\fullyS}$ its optimal value.}

\noindent{\bf Fully adaptive.} We define $\Pi_{\fullyA}$ as the class of \emph{fully adaptive} policies that process agents in $\A$ one by one. {These policies are the most ``powerful'' since they are allowed to process different types of agents in consecutive iterations, i.e., by alternating sides.} Moreover, in each iteration, these policies can use the information obtained from the choices made by previous processed agents. \revised{Practically speaking, these policies are abstractions of platforms that must update their assortments immediately after an agent makes a choice, for instance, a recently booked freelancer cannot be displayed in subsequent assortments.} We define  \fullyadaptive~as Problem~\eqref{eq:general_problem} restricted to policies in $\Pi_{\fullyA}$ and $\OPT_{\fullyA}$ its optimal value. We observe that \fullyadaptive~can be solved via a dynamic programming formulation with exponentially many states and variables, which we include in Appendix~\ref{sec:DP_fullyadaptive}.
We note that $\Pi_{\fullyA}\supset\Pi_{\fullyS}$ since any policy in $\Pi_{\fullyS}$ can be simulated by a policy in $\Pi_{\fullyA}$ as follows: \revised{the assortments are decided in advance and are not updated while processing each agent one by one.} 

{Secondly, let us describe our main policy classes with one-way interactions, i.e., where the policy must finish processing an entire side first before processing someone from the opposite side. In practice, several matching platforms adopt these type policies for a variety of reasons, such as privacy and security concerns or because they tailor it to a specific market. For example, we find one-way interactions in lodging and carpooling apps where customers (e.g., travelers, passengers) must reach out to suppliers (e.g., hosts, drivers). \revised{However, the opposite one-way interaction is also possible where a supplier (e.g., someone looking for help to complete a task) must initiate the interaction with a customer (e.g., freelancer, handyman).}\vspace{0.5em}

\noindent{\bf One-sided static.} We define $\Pi_{\onesidedS}$ as the class of \emph{one-sided static} policies which {is comprised by two types of policies}: (1) those that simultaneously process all the customers first, it observes their choices and then it processes all the suppliers, {denoted by $\Pi_{\ConesidedS}$}; (2) those that simultaneously process all the suppliers first, it observes their choices and then it processes all the customers, {denoted by $\Pi_{\SonesidedS}$}. {Therefore, we have $\Pi_{\onesidedS} = \Pi_{\ConesidedS}\cup\Pi_{\SonesidedS}$.}
{We emphasize that policies in $\Pi_{\onesidedS}$ cannot randomize between sides, i.e., they must commit to the side which is processed first and then run a policy on that side.
Note that one-sided static policies are at least as good as the fully-static ones because they can observe the choices of one side first, allowing them to adapt the assortments for the opposite side.} Moreover, once one side is entirely processed, then the decision-maker cannot obtain extra information by adaptively processing the opposite side, since our model assumes that their choices are done independently. Therefore, processing the opposite side one by one or all at the same time would achieve the same market outcome. 
{We define \onesidedstatic~as Problem~\eqref{eq:general_problem} restricted to policies in $\Pi_{\onesidedS}$ and $\OPT_{\onesidedS}$ its optimal value. Similarly, we define \Conesidedstatic~and \Sonesidedstatic~as Problem~\eqref{eq:general_problem} restricted to $\Pi_{\ConesidedS}$ and $\Pi_{\SonesidedS}$, respectively. We denote by $\OPT_{\ConesidedS}$ and $\OPT_{\SonesidedS}$ their respective optimal objective values which satisfy
\(
\OPT_{\onesidedS} = \max\Big\{\OPT_{\ConesidedS},\OPT_{\SonesidedS}\Big\}.
\)
Finally, we remark that the problem introduced in \citep{ashlagi_etal19} corresponds to \Conesidedstatic~under specific choice models for customers and suppliers.}\vspace{0.5em}

\noindent{\bf One-sided adaptive.} We define $\Pi_{\onesidedA}$ as the class of \emph{one-sided adaptive} policies which {is comprised by two types of policies}: (1) those that process all the customers one by one (observing each of their choices), and then all the suppliers, {denoted by $\Pi_{\ConesidedA}$}; (2) those that process all the suppliers one by one, and then all the customers, {denoted by $\Pi_{\SonesidedA}$}. 
{
Therefore, we have $\Pi_{\onesidedA}=\Pi_{\ConesidedA}\cup\Pi_{\SonesidedA}$. As before, we remark that policies in $\Pi_{\onesidedA}$ cannot randomize between sides, i.e., they must commit to the side which is processed first and then run an adaptive policy on that side.} 
We similarly define the optimization problems and their optimal values as \onesidedadaptive, \Conesidedadaptive, \Sonesidedadaptive, and $\OPT_{\onesidedA}$, $\OPT_{\ConesidedA}$, $\OPT_{\SonesidedA}$, respectively. Moreover, these values satisfy
\(
\OPT_{\onesidedA} = \max\Big\{\OPT_{\ConesidedA},\OPT_{\SonesidedA}\Big\}.
\)
Finally, as we noted with the fully adaptive and static policies, we have that $\Pi_{\onesidedA}\supset\Pi_{\onesidedS}$.
\begin{observation}\label{obs:policy_inclusion}
All the policy classes defined above are nested in the following way:
$$
\Pi_{\fullyS}\subset\Pi_{\onesidedS}\subset\Pi_{\onesidedA}\subset\Pi_{\fullyA},
$$
which means that $\OPT_{\fullyS}\leq\OPT_{\onesidedS}\leq\OPT_{\onesidedA}\leq \OPT_{\fullyA}$. Moreover, any other feasible policy in the whole space $\Pi$ has an objective value between $\OPT_{\fullyS}$ and $\OPT_{\fullyA}$. %
\end{observation}

%% file: 3_adaptivity_gaps.tex

\section{Adaptivity Gaps}\label{sec:adaptivity_gaps}
Our first worst-case result shows that, in general, fully-static policies can perform arbitrarily worse than one-sided static policies. 
\begin{proposition}\label{prop:gap_static_zero}
There exists an instance such that
\(
\OPT_{\onesidedS}= \Omega \left( n \right)\cdot \OPT_{\fullyS}.
\)
\end{proposition}
{Practically speaking, this result} implies that a platform would clearly benefit by allowing the best possible side to initiate the matching process since this would {reveal} more information about the agents' choices. Observe that, in the example of Proposition~\ref{prop:gap_static_zero}, if we compare \Sonesidedstatic~and \fullystatic, then the gap is actually 1. So, it is important to note that in our definition of \onesidedstatic~we consider the best between \Conesidedstatic~and \Sonesidedstatic. We defer the proof of Proposition~\ref{prop:gap_static_zero} to Appendix~\ref{sec:missing_proofs_adaptivitygaps}.

 Our first main result in this section is given in the following theorem.  
\begin{theorem}\label{thm:gap_onesided_adaptive}
The adaptivity gap between \onesidedstatic~and \onesidedadaptive~is {$e/(e-1)$}.
\end{theorem}
{ Theorem \ref{thm:gap_onesided_adaptive} implies that platforms with one-way interactions can obtain at most approximately $1.58$ times more matches (at optimality) by processing agents adaptively on one side relative to computing those assortments in advance. For instance, consider a carpooling app where customers must reach out to drivers. On the one hand, computing all the customers' assortments in advance can lead to choice congestion, i.e., popular drivers getting more requests than they can actually handle. On the other hand, displaying customers' assortments adaptively allows the platform to shift the attention to other drivers that may be able to accept the requests.}

The proof of Theorem~\ref{thm:gap_onesided_adaptive} is given in Section~\ref{sec:onesided_adaptive}, and is separated in two lemmas. First, in Lemma~\ref{lemma:LB_onesided_adaptive}, we present a novel LP relaxation of \onesidedadaptive, which along with the notion of \emph{correlation gap} introduced in \citep{agrawal2010correlation}, allows us to show that the adaptivity gap is {at most $e/(e-1)$} thanks to the monotonicity and submodularity of the demand functions. In Lemma~\ref{lemma:UB_onesided_adaptive}, we present a non-trivial instance where the gap is arbitrarily close to {$e/(e-1)$} as the size of the market grows.

Our second main result is given by the following theorem.
\begin{theorem}\label{thm:gap_fully_adaptive}
The adaptivity gap between \onesidedadaptive~and \fullyadaptive~is {$2$}.
\end{theorem}
{Theorem \ref{thm:gap_fully_adaptive} implies that platforms can double the expected number of matches if they adopt a two-way interaction model with adaptive policies relative to the equivalent one-way interaction design. For example, an outsourcing app with one-way interactions
would benefit by allowing freelancers to reach out those who need to outsource a task.}
The proof of Theorem~\ref{thm:gap_fully_adaptive} is given in Section~\ref{sec:fully_adaptive}, and as before, it is separated in two lemmas. First, in Lemma~\ref{lemma:LB_fully_adaptive}, we show that the adaptivity gap is {at most $2$} via a coupling argument between fully adaptive and one-sided adaptive policies. Finally, in Lemma~\ref{lemma:UB_fully_adaptive}, we present an instance for which the gap is arbitrarily close to {$2$} as the market's size increases.

{
Note that the previous two theorems ensure that the adaptivity gap between \onesidedstatic\ and \fullyadaptive\ lies between $2$ and $2e/(e-1)\approx 3.163$. Our next result establishes the precise gap.
\begin{theorem}\label{thm:gap_os_fa}
    The adaptivity gap between \onesidedstatic\ and \fullyadaptive\ is $2e/(e-1)$.
\end{theorem}
The proof relies on constructing a non-trivial family of nearly tight instances that we present in Appendix~\ref{apx:os_fa_gap}.
\begin{observation}\label{prop:independent_model}
Note that the adaptivity gap between both extremes, \fullystatic~and \fullyadaptive, is equal to 1 in the special case where  every agent follows the independent demand choice model, which has been widely used in practice, e.g., in the airline industry~\citep{lee1993model,van2005future}. This means that, when the choice probability of the {agents} do not depend on the other offered alternatives, there is no advantage of one policy class relative to other. To gain some intuition, let us formally define, for any customer $i\in\C$, the probability that $i$ chooses $j$ when assortment $S\subseteq\S$ is offered as
\[
\phi_i(j,S)= \left\{\begin{matrix}
    \phi_{i,j} & \text{if}\; j\in S,\\
    0 & \text{otherwise},
\end{matrix}\right. \qquad \text{and} \qquad \phi_i(0,S) = 1-\sum_{j\in S}\phi_{i,j},
\]
 where the values $\phi_{i,j}\in[0,1]$ satisfy $\sum_{j=1}^m \phi_{i,j} \leq  1$. %
Analogously, we  define the independent demand choice model of every supplier $j$ with probability values $\phi_{j,i}$ for any $i\in\C$. Intuitively, this observation follows since \fullystatic~corresponds to the problem of selecting the set of edges that maximizes the weight of the subgraph, where the weight is determined by the choice probabilities. Since the assortments are not restricted by any constraint, then the optimal solution of this problem is to select every edge. On the other hand, for \fullyadaptive, note that displaying everything in each assortment is the optimal action in the DP formulation~\eqref{eq:dp_fullyadaptive}, see Appendix~\ref{sec:DP_fullyadaptive}.
\revised{A key aspect of this observation is that it holds only in the unconstrained setting; once cardinality constraints are imposed on the assortments, one can easily construct examples where the gap
is not 1}. 
\end{observation}}

\subsection{Gap Between \onesidedstatic~and \onesidedadaptive}\label{sec:onesided_adaptive}
First, we focus on the following lower bound of the adaptivity gap between \onesidedstatic~and \onesidedadaptive.

\begin{lemma}\label{lemma:LB_onesided_adaptive}
For every instance, we have that 
{\(
\OPT_{\onesidedA}\leq\left(\frac{e}{e-1}\right)\cdot\OPT_{\onesidedS}.
\)}
\end{lemma}

\proof{Proof.}
{Recall that 
\[
\OPT_{\onesidedS} = \max\Big\{\OPT_{\ConesidedS},\OPT_{\SonesidedS}\Big\} \quad \text{and} \quad \OPT_{\onesidedA} = \max\Big\{\OPT_{\ConesidedA},\OPT_{\SonesidedA}\Big\}.
\] 
Therefore, if we show the inequality for each initiating side, then the result follows.} Fix the initiating side to $\C$. We will show that for any instance, we have 
\begin{equation}\label{eq:c_onesided_gap}
\OPT_{\ConesidedS}\leq \left(\frac{e}{e-1}\right)\cdot\OPT_{\ConesidedA}.
\end{equation}
The same result can be analogously obtained when $\S$ is the initiating side. In this proof, we use that for one-sided policies (either static or adaptive), once a policy processes all the customers, then processing all the suppliers simultaneously or one by one achieves the same objective value, because their choices are independent. In fact, due to the monotonicity of the demand functions (Assumption~\ref{assumption:monotone_submodular}), and since there are no constraints, if $C$ is the set of customers that chose supplier $j$ after all customers were processed, then the best decision is to show exactly $C$ to $j$ because: (i)  showing anything outside $C$ cannot lead to a match (customers outside $C$ did not pick $j$); (ii) showing a smaller set will decrease the probability that $j$ picks someone, i.e., their demand $f_j(\cdot)$.

To prove Inequality~\eqref{eq:c_onesided_gap}, \revised{we will proceed in three steps. First, we introduce a novel LP relaxation, show that it upper bounds $\OPT_{\ConesidedA}$ and call this result Lemma~\ref{lemma:relaxation_onesided_customers}. Then, in Lemma~\ref{thm:correlation_gap_agrawaletal}, we recall a key result from the literature of submodular optimization. Finally, we conclude the proof of Inequality~\eqref{eq:c_onesided_gap} by appropriately linking the two previous results.}

\noindent\textbf{\underline{Step 1}:} We consider the following formulation:
\begin{align}\label{eq:relaxation_onesided_customers}
\max &\quad \sum_{j\in \S}\sum_{C\subseteq \C}f_j(C)\cdot \lambda_{j,C}\\
s.t. &\quad \sum_{C\subseteq \C}\lambda_{j,C} = 1, \hspace{9.8em} \text{for all} \ j\in \S \notag\\
&\quad \sum_{C: C\ni i} \lambda_{j,C} = \sum_{S: S\ni j}\tau_{i,S}\cdot\phi_i(j,S), \hspace{2.3em} \text{for all} \ i\in \C, \ j\in \S \notag\\
&\quad \sum_{S\subseteq \S} \tau_{i,S} =1, \hspace{10em}\text{for all} \ i\in \C, \notag \\
&\quad \lambda_{j,C}, \ \tau_{i,S}\geq 0, \hspace{9.8em} \text{for all} \ j\in \S, \ C\subseteq \C, \ i\in \C, \ S\subseteq \S. \notag
\end{align}
We can interpret $\lambda_{j,C}$ as the probability that supplier $j$ gets chosen by exactly customers in $C$ \revised{(which we refer to as the \emph{backlog} of $j$)}  and $\tau_{i,S}$ as the probability that we show assortment $S$ to customer~$i$. Note that both $C$ and $S$ can potentially be $\emptyset$. The first and third constraints are distribution constraints. The second constraint corresponds to the \revised{probability that a customer $i$ is in the backlog of supplier $j$, i.e., the probability that $i$ sees and chooses $j$.}
We note that a closely related relaxation was studied for the i.i.d. online arrival model introduced by \citet{aouad2020online}, however, our relaxation is more general in that it applies to the adaptive setting that we consider in this paper and to the online model studied in \citep{aouad2020online}. Problem~\eqref{eq:relaxation_onesided_customers} will play a crucial role in our proof and, in fact, it is an upper bound of \Conesidedadaptive:
\begin{lemma}\label{lemma:relaxation_onesided_customers}
Problem~\eqref{eq:relaxation_onesided_customers} is a relaxation of {\Conesidedadaptive}.
\end{lemma}
\revised{In other words, Lemma~\ref{lemma:relaxation_onesided_customers} states that any feasible solution for \Conesidedadaptive\ can be represented by a feasible point in Problem~\eqref{eq:relaxation_onesided_customers}. Moreover, the objective value of  Problem~\eqref{eq:relaxation_onesided_customers} upper bounds that of \Conesidedadaptive.}
We provide the proof of Lemma~\ref{lemma:relaxation_onesided_customers} in Appendix~\ref{sec:proof_relaxation_onesided_customers}.

\noindent \textbf{\underline{Step 2}:} One important aspect of Problem~\eqref{eq:relaxation_onesided_customers} is that it computes the best possible distribution over assortments for each supplier and it correlates that distribution with the distribution over assortments for customers. This will be crucial when we compare it with the solution of  {\sf($\C$-One-sided Static)}. 
Now we focus on {\sf($\C$-One-sided Static)}. Note that this setting corresponds to the model introduced by \citet{ashlagi_etal19}. In this problem, a policy can be simply viewed as simultaneously selecting a family of assortments $S_1,\ldots,S_n$ where subset $S_i$ will be shown to customer $i\in\C$. 
As noted in \citep{torrico21}, 
we can consider a randomized solution in which for each customer $i\in\C$, we have a distribution $\tau_i$ over assortments and the objective value in {\sf($\C$-One-sided Static)}  corresponds to
\begin{equation}\label{eq:objective_c_onesidedstatic}
\sum_{j\in \S}\sum_{C\subseteq \C}f_j(C)\cdot \prod_{i\in C}\Bigg(\sum_{\substack{S\subseteq\S:\\ S\ni j}}\tau_{i,S}\cdot\phi_i(j,S)\Bigg)\prod_{i\in\C\setminus C}\Bigg(1-\sum_{\substack{S\subseteq\S:\\ S\ni j}}\tau_{i,S}\cdot\phi_i(j,S)\Bigg),
\end{equation}
In fact, the following distribution, for all $j\in\C, \ C\subseteq\C$,
\begin{equation}\label{eq:independent_distribution}
\lambda^{\sf ind}_{j,C} = \prod_{i\in C}\left(\sum_{S\subseteq\S: S\ni j}\tau_{i,S}\cdot\phi_i(j,S)\right)\prod_{i\in\C\setminus C}\left(1-\sum_{S\subseteq\S: S\ni j}\tau_{i,S}\cdot\phi_i(j,S)\right), 
\end{equation}
is a feasible solution in Problem~\eqref{eq:relaxation_onesided_customers} which can be easily checked by the reader. A distribution of this form is referred to as the \emph{independent distribution} (also known as product distribution) and we will denote it as $\lambda^{\sf ind}$.
We are now ready to finalize the proof of Lemma~\ref{lemma:LB_onesided_adaptive}. \citet{agrawal2010correlation}  define the \emph{correlation gap} between the ``worst-case'' distribution and the independent distribution. In particular, for any monotone submodular objective function, \citet{agrawal2010correlation} show the following:
\begin{lemma}[\citet{agrawal2010correlation}]\label{thm:correlation_gap_agrawaletal}
For any monotone submodular function $g$,  we have that
\(
\EE_{A\sim D^{\sf ind}}[g(A)] \geq \left(1-\frac{1}{e}\right)\cdot\EE_{A\sim D}[g(A)],
\)
where $D^{\sf ind}$ is the independent distribution and $D$ is any distribution over subsets with the same marginals.
\end{lemma}

\noindent\textbf{\underline{Step 3}:} We conclude the proof of Lemma~\ref{lemma:LB_onesided_adaptive} as follows: let $\tau^\star$ and $\lambda^\star$ be an optimal solution of Problem~\eqref{eq:relaxation_onesided_customers}. Therefore, by Lemma~\ref{lemma:relaxation_onesided_customers} we know that $\OPT_{\sf(\C\text{-}OA)}$ is upper bounded by the optimal value  of Problem~\eqref{eq:relaxation_onesided_customers} which is determined by $\lambda^\star$. On the other hand, we construct the independent distribution $\lambda^{\sf ind}$ as in Equation~\eqref{eq:independent_distribution} with marginal values 
\(
\sum_{S\subseteq\S: S\ni j}\tau^\star_{i,S}\cdot\phi_i(j,S).
\)
Note that $\lambda^\star$ and $\lambda^{\sf ind}$ have the same marginal values. Then, we obtain the following bound for the ratio between the optimal values of {\sf($\C$-One-sided Static)} and {\sf($\C$-One-sided Adaptive)},
\begin{equation}\label{eq:aux_correlation_gap}
\frac{\OPT_{\ConesidedS}}{\OPT_{\ConesidedA}} \geq \frac{\sum_{j\in\S}\sum_{C\subseteq\C}f_j(C)\cdot\lambda_{j,C}^{\sf ind}}{\sum_{j\in\S}\sum_{C\subseteq\C}f_j(C)\cdot\lambda_{j,C}^{\star}} \geq 1-\frac{1}{e},
\end{equation}
where the first inequality follows by noting that $\tau^\star$ can be used as a randomized feasible solution in {\sf($\C$-One-sided Static)} whose objective value is equal to Expression~\eqref{eq:objective_c_onesidedstatic} evaluated {at} $\tau^\star$ and is exactly the expected value with respect to the independent distribution $\lambda^{\sf ind}$. Also note that $\OPT_{\ConesidedA}$ is upper bounded by the optimal objective value of Problem~\eqref{eq:relaxation_onesided_customers}. The last inequality is due to Lemma~\ref{thm:correlation_gap_agrawaletal}. Finally, an analogous relaxation can be constructed when $\S$ is the initiating side and the analysis follows to show that
{
\(
\OPT_{\SonesidedA}\leq \left(\frac{e}{e-1}\right)\cdot\OPT_{\SonesidedS}.
\)}
\endproof

The second part of the proof of Theorem~\ref{thm:gap_onesided_adaptive} consists of lower bounding the  gap between \onesidedstatic~and \onesidedadaptive.
\begin{lemma}\label{lemma:UB_onesided_adaptive}
There exists an instance such that
{
\(
\OPT_{\onesidedA} \geq \left(1-\mathcal{O}\left(\frac{\log n}{n}\right)\right)\cdot \frac{e}{e-1}\cdot\OPT_{\onesidedS}.
\)}
In particular, the adaptivity gap tends to {$e/(e-1)$} as $n$ goes to $\infty$.
\end{lemma}
{The instance in Lemma~\ref{lemma:UB_onesided_adaptive} consists of all agents having a demand function that tends to 1 as the size of the offered set increases.}
We provide the proof of this lemma in Appendix~\ref{sec:missing_proofs_adaptivitygaps}.

\subsection{Gap Between \onesidedadaptive~and \fullyadaptive}\label{sec:fully_adaptive}
First, we focus on the following upper bound of the adaptivity gap.
\begin{lemma}\label{lemma:LB_fully_adaptive}
For every instance, we have that 
{
\(
\OPT_{\fullyA}\leq2\cdot\OPT_{\onesidedA}.
\)}
\end{lemma}

\proof{Proof.}
Consider an optimal policy $\pi^*$  for \fullyadaptive. For each agent $a \in \A$, we denote as $X_{a}$ the random variable which takes value $1$ if $a$ is matched in $\pi^*$ and  $0$ otherwise.  We denote by $Y_{a}$ the random variable  which takes value $1$ if $a$ is matched to an opposite agent that was processed before $a$ in $\pi^*$ and  $0$ otherwise. We denote by $Z_{a}$ the random variable  which takes value $1$ if $a$ is matched to an opposite agent that was processed after $a$ in $\pi^*$ and  $0$ otherwise. Clearly, we have $X_{a} =Y_{a}+Z_{a}$. 
Given this notation, observe that $\sum_{i\in \C}Y_{i}=\sum_{j\in \S}Z_{j}$ and $\sum_{j\in \S}Y_{j}=\sum_{i\in \C}Z_{i}$. 
    Let $\Ma_{\pi^*}$ be the number of matches obtained by $\pi^*$ which is given by
$$      
        \Ma_{\pi^*}=\sum_{i\in \C}X_{i}= \sum_{i\in \C} Y_{i}+ \sum_{i\in \C}  Z_{i} =  \sum_{i\in \C}  Y_{i} + \sum_{j\in \S
} Y_{j}.
$$   
Let us define a feasible (randomized) one sided adaptive policy $\pi_{\C}$ that processes customers first. We will couple this policy with $\pi^*$. In particular, $\pi_{\C}$ offers exactly the same assortment as $\pi^*$ whenever $\pi^*$ processes a customer. If $\pi^*$ processes a supplier before finishing all customers at some time, then $\pi_{\C}$ does not process any agent at that time but just simulates the choice of the supplier so that $\pi^*$ and $\pi_{\C}$ follow the same decision tree.\footnote{\revised{Recall that every policy can be described by a decision tree where each node corresponds to the \emph{current state} and the edges represent each possible \emph{action}. A policy then follows this decision tree (starting at the initial state, i.e., the root of the tree) and may either deterministically or stochastically decide the next action at each node.}} Note that $\pi_{\C}$ is a randomized policy because of the simulation. After processing all customers, $\pi_{\C}$ offers to each supplier the same subset of customers as in $\pi^*$. Let $\Ma_{\C}$ denote the number of matches obtained by $\pi_{\C}$, we will show that for each supplier $j \in \S$, the probability of matching $j$ in  the policy $\pi_{\C}$ is greater than $\EE [Y_j]$. Recall that $\EE [Y_j]$ corresponds to the probability of matching $j$ in the optimal policy $\pi^*$ to a customer that was processed before $j$.   
In fact, suppose $\pi^*$ is processing supplier $j$ at some point (a node in the decision tree of $\pi^*$) and offers her assortment $C\subseteq\C$. We can write $C$ as $C= C_o \cup C_f$ where $C_o$ is the \revised{backlog of supplier $j$} 
(note that it is always optimal to include this subset in the assortment for supplier $j$) and $C_f$ is a subset from the set of customers that have not been processed yet by $\pi^*$. These are the scenarios  right after processing $j$: 

- Supplier $j$  chooses an agent in $C_f$ or the outside option, this happens with a total probability of $\sum_{i \in C_f \cup\{0\} }  \phi_j(i,C)$. In that case,  we have $Y_j=0$.

- Supplier $j$  chooses an agent in $C_o$, this happens with a total probability of $\sum_{i \in C_o  }  \phi_j(i,C)$. In that case,  we have $Y_j=1$.

By taking expectation over the scenarios described above we get $\EE [Y_j]=\sum_{i \in C_o} \phi_j(i,C)$. On the other hand, in all the scenarios  above, the policy $\pi_{\C}$  offers at the end the set $C=C_o \cup C_f$ to supplier $j$. From the set $C$, we know that all customers in $C_o$ chose $j$, but customers in $C_f$ might or might not have chosen $j$ depending on each  scenario. Therefore,  the probability of matching $j$ in the policy $\pi_{\C}$  is greater or equal than $\sum_{i \in C_o} \phi_j(i,C) $.  This implies that the probability of matching $j$ in  the policy $\pi_{\C}$ is greater than $\EE [Y_j]$. By taking the sum over all $j$ and the expectation over all sample paths, we get
$\EE [\Ma_{\C}] \geq \EE  [\sum_{j\in \S} Y_{j}]$. Therefore, for \Conesidedadaptive~we have,
$  \OPT_{\ConesidedA} \geq    \EE  [\sum_{j\in \S} Y_{j}].    $
A similar argument can be used to show that
$  \OPT_{\SonesidedA} \geq    \EE  [ \sum_{i\in \C} Y_{i}   ]      $. Finally, we conclude that
  \begin{align}  
  \OPT_{\onesidedA} = \max \Big\{\OPT_{\ConesidedA},\OPT_{\SonesidedA} \Big\} & \geq \frac{1}{2} \cdot(\OPT_{\ConesidedA}+\OPT_{\SonesidedA}) \notag \\
  & \geq \frac{1}{2}\cdot
  \EE\Bigg[ \sum_{j\in \S} Y_{j}\Bigg]+  \frac{1}{2}\cdot \EE  \Bigg[\sum_{i\in \C} Y_{i}\Bigg]  \notag\\
  &= \frac{1}{2} \cdot\EE[\Ma_{\pi^*}]= \frac{1}{2}\cdot \OPT_{\fullyA}.  \label{eq:aux_ineq_fa}
  \end{align}
\endproof
We conclude this section by presenting an instance in which the adaptivity gap is arbitrarily close to 2.
\begin{lemma}\label{lemma:UB_fully_adaptive}
There exists an instance such that
{\(
\OPT_{\fullyA}= 2\cdot\left(1 - \mathcal{O}\left(\frac{1}{n}\right)\right)\cdot \OPT_{\onesidedA}.
\)}
In particular, the adaptivity gap tends to {$2$} as $n$ goes to $\infty$.
\end{lemma}
{The  instance in Lemma~\ref{lemma:UB_fully_adaptive} features two distinct sub-markets: one with 1 supplier and $n-1$ customers, and the other with 1 customer and $n-1$ suppliers. We designed a choice structure where each sub-market operates independently, ensuring zero cross-selection probability between agents of different sub-markets, which enables us to establish the gap. }
We defer the proof of this result to Appendix~\ref{sec:missing_proofs_adaptivitygaps}.

%% file: 4_approx_fully_adaptive.tex

\section{Approximation Guarantees for \fullyadaptive}\label{sec:approx_fully_adaptive}
{The LP upper bound~\eqref{eq:relaxation_onesided_customers} was instrumental in establishing the adaptivity gap between \onesidedstatic~and \onesidedadaptive. In this section, we further leverage this relaxation as a central tool in analyzing the approximation guarantee for \onesidedadaptive, and consequently, for \fullyadaptive.}

First, we know that \citet{torrico21} show a $(1-\frac{1}{e})$-approximation for \Conesidedstatic. Therefore, by applying this result to both sides of the market and using our adaptivity gaps obtained in Section~\ref{sec:adaptivity_gaps}, we can immediately conclude the following constant guarantees:
\begin{corollary}\label{coro:simple_approx_factor}
There exists a one-sided static policy that runs in polynomial-time and guarantees a $(1-1/e)^2\approx 0.399$ approximation factor for \onesidedadaptive~and a $\frac{1}{2}\cdot(1-1/e)^2\approx 0.199$ for \fullyadaptive.
\end{corollary}
We defer the details of the proof of Corollary  \ref{coro:simple_approx_factor} to Appendix~\ref{sec:proof_corollary_simple_approx_factor}. Our main goal in this section is to obtain better guarantees for \fullyadaptive\ than in Corollary~\ref{coro:simple_approx_factor}.
{\color{black}
To establish our approximation guarantees, the main building block is provided in the following theorem:
\begin{theorem}\label{thm:approxfactor_greedy}
There exists a polynomial-time algorithm that achieves, in expectation, a $1/2$ approximation factor for \Conesidedadaptive. Moreover, the same approximation factor holds for its adaptation to \Sonesidedadaptive.
\end{theorem}

From a practical perspective, Theorem~\ref{thm:approxfactor_greedy} can be interpreted as follows: if the platform is committed to a one-sided adaptive policy, i.e., it has predetermined whether customers or suppliers initiate the matching, then there exists a policy that achieves at least half of the optimal number of matches within the corresponding restricted policy space (i.e., $\OPT_{\ConesidedA}$ or $\OPT_{\SonesidedA}$, respectively).
To prove Theorem~\ref{thm:approxfactor_greedy}, we design a greedy method that we name {\em Arbitrary Order Greedy}. The analysis of this algorithm relies on the LP relaxation introduced in Section~\ref{eq:relaxation_onesided_customers} and draws inspiration from the randomized primal-dual technique initially introduced by \citet{devanur2013randomized} for online matching.

We leverage Theorem~\ref{thm:approxfactor_greedy} to establish our next main result in Theorem \ref{coro:approx_fully_adaptive}, which provides a $1/4$-approximation guarantee for \fullyadaptive. The corresponding algorithm is randomized: with equal probability, it selects which side initiates the matching process and then runs the Arbitrary Order Greedy on that side. We present detailed analyses of Theorems~\ref{thm:approxfactor_greedy} and \ref{coro:approx_fully_adaptive} in Sections~\ref{sec:thm4} and \ref{sec:thm5}, respectively.
\begin{theorem}\label{coro:approx_fully_adaptive}
There exists a randomized polynomial-time policy that achieves, in expectation, a $1/4$ approximation factor for \fullyadaptive.
\end{theorem}

Finally, we address a natural question for platforms that have committed to one-sided adaptive policies but have not yet decided which side should initiate. In this case, we propose a pseudo-polynomial time sampling-based method that identifies, with high probability, the better initiating side and applies the greedy algorithm accordingly. This leads to the following result:
\begin{corollary}\label{coro:approx_onesided_adaptive}
There exists a pseudo-polynomial-time sampling policy  that achieves (with high probability) a $1/2-\epsilon$ approximation factor for \onesidedadaptive, for any $\epsilon>0$.
\end{corollary}
The policy consists of a sampling subroutine that estimates, through repeated simulations, the expected number of matches resulting from applying the greedy algorithm on each side. It then deterministically selects the side with the higher estimated value. This policy, which belongs to the class $\Pi_{\onesidedA} = \Pi_{\ConesidedA} \cup \Pi_{\SonesidedA}$, thus guarantees a constant-factor approximation relative to the optimal one-sided adaptive benchmark $\OPT_{\onesidedA} = \max \{\OPT_{\ConesidedA}, \OPT_{\SonesidedA}\}$.
Compared to the result in Corollary~\ref{coro:simple_approx_factor}, the approximation factor in Corollary~\ref{coro:approx_onesided_adaptive} is better, and the policy is adaptive rather than static. However, the sampling procedure in Corollary~\ref{coro:approx_onesided_adaptive} runs in a pseudo-polynomial-time, in contrast to the polynomial-time algorithm in Corollary~\ref{coro:simple_approx_factor}.
The complete analysis of Corollary~\ref{coro:approx_onesided_adaptive} is presented in Section~\ref{sec:thm6}.
 }

\subsection{Proof of Theorem~\ref{thm:approxfactor_greedy}} \label{sec:thm4}

{Our main algorithm to prove Theorem~\ref{thm:approxfactor_greedy}~is formalized in Algorithm~\ref{alg:random_order_greedy}. } 
\begin{algorithm}[htpb]
 \caption{(Arbitrary Order) Greedy for Customers}\label{alg:random_order_greedy}
 \begin{algorithmic}[1]
 \Require  Choice models $\phi$ {and an arbitrary order $i_1,\ldots,i_n$ over $\C$.}
 \Ensure Assortments $S_i$ for all $i\in\C$.
\State For each $j\in \S$, set $C^0_j=\emptyset$.
\For{$t=1,\ldots,n$} 
\State $S_{i_t} \in \argmax \left\{\sum_{j\in S}\left[f_j(C^{t-1}_j \cup \{i_t\})-f_j(C^{t-1}_j)\right]\cdot\phi_{i_t}(j,S): \ S \subseteq \S \right\}$
\State Observe the choice $\ell\in S_{i_t}\cup\{0\}$ of customer $i_t$ which happens w.p. $\phi_{i_t}(\ell,S_{i_t})$.
\State If $\ell\neq 0$, update $C_{\ell}^t \leftarrow C_{\ell}^{t-1}\cup\{i_t\}$ and $C_j^{t} \leftarrow C_j^{t-1}$ for all $j\in\S\setminus\{\ell\}$. Otherwise, $C_j^{t} \leftarrow C_j^{t-1}$ for all $j\in\S$.
\EndFor
 \end{algorithmic}
 \end{algorithm}
\citet{torrico21} show that the greedy algorithm (Algorithm~\ref{alg:random_order_greedy}) achieves a $1/2$ factor for the \Conesidedstatic~problem {and \citet{aouad2020online} show the same guarantee for the online model.} We will show, via a novel primal-dual analysis in the context of two-sided assortment optimization, that this algorithm actually achieves a $1/2$ approximation guarantee for \Conesidedadaptive~{with respect to an order-independent linear relaxation.} Roughly speaking, this algorithm considers an arbitrary permutation over the set of customers and treats that order as a stream of customers arrivals. For each arriving customer, the algorithm greedily offers the assortment of suppliers that maximizes the \emph{weighted} probability of choosing a supplier. For completeness, we formalize this method in Algorithm~\ref{alg:random_order_greedy} for \Conesidedadaptive. \citet{torrico21} considers this algorithm for a revenue maximizing objective, however, we are interested only {in} the expected size of the matching.
Note that the subproblem in Line 4 can be efficiently solved thanks to the first point in Assumption~\ref{assumption:monotone_submodular}. We remark that in Line 5, the method in \citep{torrico21} writes ``simulate the choice of'' instead of ``observe the choice of''. We make this distinction here because for \Conesidedstatic~the platform cannot access the choices of each customer but it can rather simulate them. Instead, for \Conesidedadaptive, the platform can actually observe these choices as they are made and, subsequently, design the future assortments. These observations serve to construct the backlog sets $C_j^n$ which, for each supplier $j$, corresponds to the set of customers that chose $j$. In this way, we are able to control the demand of each supplier by using the marginals $f_j(C^{t-1}_j\cup\{i_t\})-f_j(C^{t-1}_j)$. {These values roughly correspond to the marginal probability that supplier $j$ chooses customer $i_t$ back}. Finally, note that Algorithm~\ref{alg:random_order_greedy} can be easily adapted for \Sonesidedadaptive~{by simply exchanging} the roles of customers and suppliers (Assumption~\ref{assumption:monotone_submodular} applies to all agents).
 \begin{proof}{Proof of Theorem~\ref{thm:approxfactor_greedy}.} 
 We show that Algorithm~\ref{alg:random_order_greedy} achieves a $1/2$ approximation factor for \Conesidedadaptive~via a randomized primal-dual analysis. 
Consider an arbitrary permutation $i_1,\ldots,i_n$ over $\C$ where each customer is treated as an arrival.
 The dual program of Problem~\eqref{eq:relaxation_onesided_customers} is:
 \begin{subequations}\label{eq:dual_relaxation_onesided_customers}
\begin{align}
 \min &\quad \sum_{j\in \S}\beta_j + \sum_{i\in \C} \alpha_i \label{eq:dualobjective_unconstrained}\\
 s.t. &\quad  \beta_j +\sum_{i\in C}\gamma_{i,j} \geq f_j(C), \hspace{1.4em} \qquad \text{for all} \ C\subseteq\C, \ j\in \S, \label{eq:demand_dualconstraint_unconstrained}\\
 &\quad \alpha_i - \sum_{j\in S}\gamma_{i,j}\cdot \phi_{i}(j,S) \geq 0, \qquad \text{for all} \ S\subseteq\S, \ i\in \C,\label{eq:choicemodel_dualconstraint_unconstrained}\\
 &\quad \alpha_i, \ \beta_j, \ \gamma_{i,j}\in\RR, \hspace{6em} \text{for all}\; i\in \C, \ j\in \S.\notag
 \end{align}
 \end{subequations}
We now construct random variables $\tilde{\alpha}_i$, $\tilde{\gamma}_{i,j}$ and $\tilde{\beta}_j$ from the actions taken by Algorithm~\ref{alg:random_order_greedy} according to the customer arrivals. As observed in \citep{devanur2013randomized}, this solution does not necessarily need to be feasible in the dual, but its expected values must be feasible. Moreover, we need to show that the dual objective value achieved by these random variables is a constant factor away from the primal objective value. 

{\bf Construction of the Variables.} For each arriving customer $i=i_t$ at step $t$ and assortment $S_i$ output by Algorithm~\ref{alg:random_order_greedy},  set $\tilde{\alpha_i}$ to
 { \(
 \Big[f_{\xi}(C^{t-1}_\xi\cup\{i\})-f_\xi(C^{t-1}_\xi)\Big]
 \)
 where $\xi:=\xi_i(S_i)$ is the random variable that indicates the supplier that $i$ chooses when assortment $S_i$ is displayed (if $i$ chooses the outside option, set $\tilde{\alpha}_i=0$). So for $j\in S_i\cup\{0\}$, we have
 \(
 \PP(\xi_i(S_i) = j) = \phi_i(j,S_i)
 \)
 and zero otherwise.} Note that $\tilde{\alpha}_i$ is, in fact, a random variable that depends on the past customers' choices and the current choice.
Set $\tilde{\gamma}_{i,j}$ for every $j\in \S$ to
 \(
 \Big[f_j(C^{t-1}_j\cup\{i\})-f_j(C^{t-1}_j)\Big]
 \).
 Note that $\tilde{\gamma}_{i,j}$ is also a random variable since it depends on the set $C^{t-1}_j$. Finally, set $\tilde{\beta}_j = f_j(C^n_j)$, which is the random value obtained from supplier $j$ at the end of the process. 
 
 {\bf Feasibility and Objective Value.}
 Now, we show that the expected value of these variables form a feasible solution in the dual and that their dual objective value is a constant factor away from the primal objective value. 
 First, let us compute their objective value. {Denote by $\xi_t$ the supplier (if any) chosen by customer $i_t$.} Clearly,
 {\begin{align*}
\sum_{t=1}^n\tilde{\alpha}_{i_t} &= \sum_{t=1}^n\Big[f_{\xi_t}(C^{t-1}_{\xi_t}\cup\{i_t\})-f_{\xi_t}(C^{t-1}_{\xi_t})\Big] =\sum_{j\in\S}\sum_{t=1}^n\Big[f_j(C^{t}_j)-f_j(C_j^{t-1})\Big]\\
 &=\sum_{j\in\S}f_j(C^n_j) - f_j(\emptyset) = \sum_{j\in\S}f_j(C^n_j),
 \end{align*}}
 where in the second equality we use that for every $j\in\S$ and $t\in[n]$, the set $C_j^t=C_j^{t-1}$ if $i_t$ does not choose $j$ and  $C_j^t=C_j^{t-1}\cup\{i_t\}$ if $i_t$ does  choose $j$. In that equality, we also observe that $C_1^n,\ldots,C_j^n$ are disjoint sets (each customer chooses only one option). The following equality corresponds to a telescopic sum. Therefore, the dual objective value is
 \[
 \sum_{i\in\C}\tilde{\alpha}_i + \sum_{j\in\S}\tilde{\beta}_j = 2\cdot\sum_{j\in\S}f_j(C^n_j),
 \]
 and the primal value is exactly $1/2$ of that expression, i.e., $\sum_{j\in\S}f_j(C^n_j)$ since we get a value $f_j(C^n_j)$ for every supplier.

 Now, we prove that the expected value of the variables are feasible in the dual program. Recall that  Constraint~\eqref{eq:demand_dualconstraint_unconstrained} is: For any $C\subseteq\C$ and $j\in\S$, the variables must satisfy
 \[
 \EE[\tilde{\beta}_j] + \sum_{i\in C}\EE[\tilde{\gamma}_{i,j}] \geq f_j(C),
 \]
 where the expectation is over the customers' choices.
 The left hand side of the inequality satisfies
\begin{align*}
 \EE[\tilde{\beta}_j] + \sum_{i\in C}\EE[\tilde{\gamma}_{i,j}] 
 &=\EE\left[f_j(C^n_j)+\sum_{t=1}^n[f_j(C^{t-1}_j\cup\{i_t\})-f_j(C^{t-1}_j)]\cdot \one_{\{i_t\in C\}}\right]\\
 &\geq \EE\left[f_j(C^n_j)+\sum_{t=1}^n[f_j(C^{n}_j\cup\{i_t\})-f_j(C^{n}_j)]\cdot \one_{\{i_t\in C\}}\right]\\
 &\geq \EE\left[f_j(C^n_j)+f_j(C^{n}_j\cup C)-f_j(C^{n}_j)]\right]\\
 & \geq f_j(C),
 \end{align*}
 where the first equality follows from observing that the sum is over $i\in C$
 the following two inequalities are due to submodularity of $f_j$. 
 On the other hand, the expected value of the variables must also satisfy Constraint~\eqref{eq:choicemodel_dualconstraint_unconstrained}: For any $i\in\C$ and $S\subseteq \S$
 \[
 \EE[\tilde{\alpha}_i] - \sum_{j\in S}\EE[\tilde{\gamma}_{i,j}]\cdot\phi_i(j,S) \geq 0.
 \]
 Let us condition first on $C_j^{t-1}$ for all $j$, where $t$ is the time period that customer $i$ arrives. Given this, variables $\tilde{\gamma}_{i,j}$ are deterministic and the following expectation is only over the random choice $\xi$ of $i$ which we denote by $\EE_\xi[\cdot]$:
 {\begin{equation}\label{eq:aux_primal_dual1}
 \EE_\xi[\tilde{\alpha}_i] = \EE_\xi\left[f_\xi(C^{t-1}_\xi\cup\{i\}) - f_\xi(C^{t-1}_\xi)\right] = \sum_{j\in S_i}[f_j(C^{t-1}_j\cup\{i\}) - f_j(C^{t-1}_j)]\cdot\phi_i(j,S_i).  
 \end{equation}
Due to the optimality of $S_i$ in Step 3 of Algorithm~\ref{alg:random_order_greedy}, we obtain that for any subset $S\subseteq \S$
\begin{equation*}
 \sum_{j\in S_i}[f_j(C^{t-1}_j\cup\{i\}) - f_j(C^{t-1}_j)]\cdot\phi_i(j,S_i) \geq \sum_{j\in S}[f_j(C^{t-1}_j\cup\{i\}) - f_j(C^{t-1}_j)]\cdot\phi_i(j,S)
= \sum_{j\in S}\tilde{\gamma}_{i,j}\cdot\phi_i(j,S),
\end{equation*}
where in the equality we use the definition of $\tilde{\gamma}_{i,j}$.}
 Constraint~\eqref{eq:choicemodel_dualconstraint_unconstrained} then follows by taking expectation over $C_j^{t-1}$. 
Since this analysis does not depend on the specific order in which customers arrive, we conclude that Algorithm~\ref{alg:random_order_greedy} achieves a $1/2$ approximation factor for \Conesidedadaptive. 
\end{proof}
 
 {
 \begin{observation}
For general interest, the proof above recovers the 1/2-approximation for the setting with online (customers) arrivals in the adversarial model initially studied in \citep{aouad2020online}. This is because their clairvoyant benchmark, which we denote by $\OPT_{\textsf{On}}$, knows the order of arrival but cannot alter it. Therefore, it is clear that 
\(
\OPT_{\ConesidedS}\leq \OPT_{\textsf{On}}\leq \OPT_{\ConesidedA}\leq \OPT_{\eqref{eq:relaxation_onesided_customers}},
\)
where $\OPT_{\eqref{eq:relaxation_onesided_customers}}$ is the optimal value of the \emph{order-independent} relaxation~\eqref{eq:relaxation_onesided_customers}. Moreover, thanks to \eqref{eq:aux_correlation_gap} we know that $(1-1/e)\cdot\OPT_{\eqref{eq:relaxation_onesided_customers}}\leq \OPT_{\ConesidedS}$.
 This implies that $\OPT_{\eqref{eq:relaxation_onesided_customers}}\leq \frac{e}{e-1}\cdot \OPT_{\textsf{On}}$, i.e.,  the gap between both benchmarks is at most $e/(e-1)$. We leave as an open question whether the gap is tight or not.
\end{observation}
\begin{observation}\label{obs:approximate_oracle}
A natural question is how much we lose in the approximation factor under an approximate oracle in Assumption~\ref{assumption:monotone_submodular}. Let us assume that, for each agent $a\in\A$, the decision-maker has access to an oracle that achieves a $\rho_a$ approximation factor for the \emph{unconstrained} single-agent assortment optimization problem. Then, in the primal-dual analysis above, we re-define the random variables as $\tilde{\alpha}_i = \frac{1}{\rho}_i\cdot \Big[f_{\xi}(C^{t-1}_\xi\cup\{i\})-f_\xi(C^{t-1}_\xi)\Big]$ for each $i\in\C$, while variables $\tilde{\beta}_j$ and $\tilde{\gamma}_{i,j}$ remain the same.\\
{\bf Feasibility.} The proof of the feasibility of their expected values in Constraint~\eqref{eq:demand_dualconstraint_unconstrained} does not change. Now, let us discuss Constraint~\eqref{eq:choicemodel_dualconstraint_unconstrained}. We note that, conditioned on the backlogs $C_j^{t-1}$, we have
\[
\EE_\xi[\tilde{\alpha}_i] = \frac{1}{\rho_i}\cdot\sum_{j\in S_i}[f_j(C^{t-1}_j\cup\{i\}) - f_j(C^{t-1}_j)]\cdot\phi_i(j,S_i).
\]
On the other hand, thanks to the $\rho_i$-approximate oracle, we know that the solution set $S_i$ found in Step 3 of Algorithm~\ref{alg:random_order_greedy} satisfies
\[
\sum_{j\in S_i}[f_j(C^{t-1}_j\cup\{i\}) - f_j(C^{t-1}_j)]\cdot\phi_i(j,S_i)\geq \rho_i\cdot\sum_{j\in S^\star}[f_j(C^{t-1}_j\cup\{i\}) - f_j(C^{t-1}_j)]\cdot\phi_i(j,S^\star),
\]
where $S^\star$ is an optimal set for the maximization problem in Step~3. Putting these observations together, and since $S^\star$ is an optimal solution, we get for all $S\subseteq\S$
\begin{align*}
\EE_\xi[\tilde{\alpha}_i]\geq \sum_{j\in S^\star}[f_j(C^{t-1}_j\cup\{i\}) - f_j(C^{t-1}_j)]\cdot\phi_i(j,S) &\geq \sum_{j\in S}[f_j(C^{t-1}_j\cup\{i\}) - f_j(C^{t-1}_j)]\cdot\phi_i(j,S) \\&= \sum_{j\in S}\tilde{\gamma}_{i,j}\cdot\phi_i(j,S). 
\end{align*}
{\bf Objective.} As in the proof of Theorem~\ref{thm:approxfactor_greedy}, we can easily show that 
 \(
 \sum_{i\in\C}\tilde{\alpha}_i \leq \frac{1}{\rho}\cdot\sum_{j\in\S}f_j(C^n_j),
 \)
 where $\rho = \min_{i\in\C}\rho_i$. Therefore,
 \[
 \sum_{i\in\C}\tilde{\alpha}_i + \sum_{j\in\S}\tilde{\beta}_j \leq \left(\frac{1}{\rho}+1\right)\cdot\sum_{j\in\S}f_j(C^n_j),
 \]
 which leads to a $\rho/(\rho+1)$ approximation factor.
 \end{observation}
}

{\color{black}

\subsection{Proof of Theorem~\ref{coro:approx_fully_adaptive}}  \label{sec:thm5}
Denote by $\EE[\Ma_\C],\; \EE[\Ma_\S]$ the expected total number of matches achieved by Algorithm~\ref{alg:random_order_greedy} on each side. Our algorithm for \fullyadaptive\ flips a fair coin to select the side that is processed first. Let \alg~be the expected value of this randomized policy.
Then, we have
\[
\alg = \frac{1}{2}\cdot\left(\EE[\Ma_\C] +\EE[\Ma_\S]\right) \geq \frac{1}{4}\cdot\left(\OPT_{\ConesidedA} +\OPT_{\SonesidedA}\right)\geq\frac{1}{4}\cdot\OPT_{\fullyA},
\]
where the first inequality follows from Theorem~\ref{thm:approxfactor_greedy} and the last inequality is due to $\OPT_{\ConesidedA} +\OPT_{\SonesidedA}\geq \OPT_{\fullyA}$, which was proved in \eqref{eq:aux_ineq_fa}.
}

\color{black}
\subsection{Proof of Corollary~\ref{coro:approx_onesided_adaptive}}  \label{sec:thm6}
Our  policy for \onesidedadaptive~is formalized in Algorithm~\ref{alg:onesidedadaptive_policy}. This procedure consists of running Algorithm~\ref{alg:random_order_greedy} sufficiently enough times to get an estimate of the expected total number of matches achieved independently for each side. Then, it selects the side with the highest estimate. When the customers are processed first, Lemma~\ref{lemma:sampling_greedy} establishes the precise number of sample runs to obtain an approximately good estimate. A similar result can be obtained when suppliers are processed first. Specifically, the number of runs is pseudo-polynomial in the size of the input as it depends on the following parameter:
\begin{equation}\label{eq:phi_min}
\phi_{\min} = \min\left\{\max_{S\subseteq\S}\phi_i(j,S):\; i\in\C,\; j\in\S,\; \max_{S\subseteq\S}\phi_i(j,S)>0\right\}.
\end{equation}
This parameter represents the smallest non-zero choice probability (given a maximizing assortment), among all customers' choice models. 
Moreover, $\phi_{\min}$ corresponds to 
\[
\min\left\{\phi_i(j,\{j\}):\; i\in\C,\; j\in\S,\; \phi_i(j,\{j\})>0\right\}
\] when the choice models are weak-substitutable (a standard assumption in the assortment optimization literature), i.e., for all $j\in S\subseteq S'\subseteq\S$, $\phi_i(j,S)\geq \phi_i(j,S')$. For example, if all customers follow MNL choice models, $\phi_{\min}$ would corresponds to $v_{\min}/(1+v_{\min})$ where $v_{\min}$ is the smallest non-zero preference weight over suppliers, among all customers.
\begin{lemma}\label{lemma:sampling_greedy}
Consider any order in $\C$ and any $\epsilon>0$, $\delta\in(0,1)$. Let $\EE[\Ma_{\pi_\C}]$ be the expected total number of matches achieved by Algorithm~\ref{alg:random_order_greedy} when customers are processed first and $\phi_{\min}$ be the parameter defined in \eqref{eq:phi_min}. Then, the sample average $\hat{\Ma}_{\pi_\C}$ over $\mathcal{O}\left(\frac{1}{\epsilon^2 \phi_{\min}}\cdot\log \frac{1}{\delta}\right)$ independent runs of Algorithm~\ref{alg:random_order_greedy} satisfies, with probability $1-\delta$, 
\[
|\hat{\Ma}_{\pi_\C} - \EE[\Ma_{\pi_\C}]| \leq \epsilon\cdot\OPT_{\ConesidedA},
\] 
This can be shown analogously for Algorithm~\ref{alg:random_order_greedy} adapted to \Sonesidedadaptive.
\end{lemma}
\begin{algorithm}[htpb]
 \caption{Policy for \onesidedadaptive}\label{alg:onesidedadaptive_policy}
 \begin{algorithmic}[1]
 \Require  Choice models $\phi$ and $T=\mathcal{O}\left(\frac{1}{\epsilon^2 \phi_{\min}}\cdot\log \frac{1}{\delta}\right)$
 \Ensure Assortments $S_i$ for all $i\in\C$ and $C_j$ for all $j\in\S$.
 \State {Consider arbitrary orders $\sigma_\C$ and $\sigma_\S$ over $\C$ and $\S$, respectively.}
 \State {Make $T$ independent runs of Algorithm~\ref{alg:random_order_greedy} over $\sigma_{\C}$} to approximately solve \Conesidedadaptive. {Let $\hat\Ma_{\pi_\C}$ be the estimate of the expected number of matches $\EE[\Ma_{\pi_\C}]$}.
 \State {Make $T$ independent runs of the adaptation of Algorithm~\ref{alg:random_order_greedy} over $\sigma_{\S}$} to approximately solve \Sonesidedadaptive. {Let $\hat\Ma_{\pi_\S}$ be the estimate of the expected number of matches $\EE[\Ma_{\pi_\S}]$}.
 \State Select the side, either $\C$ or $\S$, with the highest estimate between $\hat\Ma_{\pi_\C}$ and $\hat\Ma_{\pi_\S}$. 
 \State Run Algorithm~\ref{alg:random_order_greedy} on the selected side.
 \end{algorithmic}
 \end{algorithm}
We defer the proof of this lemma to Appendix~\ref{appendix:proof_lemma_sampling}. We are now ready to prove Corollary~\ref{coro:approx_onesided_adaptive}.
\begin{proof}{Proof of Corollary~\ref{coro:approx_onesided_adaptive}.}
 Let $\hat\Ma_{\pi_\C}$ and $\hat\Ma_{\pi_\S}$ be the estimate obtained in Lines 2 and 3 for the expected total number of matches  $\EE[\Ma_{\pi_\C}]$ and $\EE[\Ma_{\pi_\C}]$, respectively, achieved by Algorithm~\ref{alg:random_order_greedy} for each side.
 First, due to Lemma~\ref{lemma:sampling_greedy} we know that, with probability $1-\delta$,
 \[
 \hat{\Ma}_{\pi_\C}\geq \EE[\Ma_{\pi_\C}]-\epsilon\cdot\OPT_{\ConesidedA}\geq \left(\frac{1}{2}-\epsilon\right)\cdot\OPT_{\ConesidedA},
 \]
 where in the second inequality we use Theorem~\ref{thm:approxfactor_greedy}, i.e., $\EE[\Ma_{\pi_\C}]\geq (1/2)\cdot\OPT_{\ConesidedA}$.

Second, since each estimate was computed with independent runs, the estimated total number of matches achieved by  Algorithm~\ref{alg:onesidedadaptive_policy} on the selected side (Step 4) can be simply computed as $\alg = \max\left\{\hat\Ma_{\pi_\C},\hat\Ma_{\pi_\S}\right\}$. We conclude the proof by noting that with probability $(1-\delta)^2\geq 1-2\delta$, we have
\[
\alg =\max\left\{\hat\Ma_{\pi_\C},\hat\Ma_{\pi_\S}\right\} \geq\left(\frac{1}{2}-\epsilon\right)\cdot\max\{\OPT_{\ConesidedA},\OPT_{\SonesidedA}\}\geq
\left(\frac{1}{2}-\epsilon\right)\cdot\OPT_{\onesidedA}.
\]
 
We remark that the pseudo-polynomial number of sample runs can be improved to a polynomial number under certain conditions discussed in Remark~\ref{remark:polynumber_samples} in Appendix~\ref{appendix:proof_lemma_sampling}. 
\end{proof}
  

%% file: 5_approx_fully_static.tex

\section{Approximation Guarantees for \fullystatic\ under MNL}\label{sec:approx_fully_static}
In this section, we focus on the \fullystatic\ problem under the multinomial logit model (MNL). Specifically, each agent (i.e., supplier or customer) makes choices according to an MNL model. We allow agents to have different MNL models. Specifically, the MNL model for a customer $i \in \C$ is defined as follows: We use $v_{ij}$ to denote the preference weight of customer $i$ for supplier $j \in \S$. The preference weight for the outside option is normalized to 1. Therefore, for any assortment of suppliers $S \subseteq \S$, the probability that customer $i$ chooses supplier $j \in S$ under the MNL model is given by $\phi_i(j,S)= \frac{v_{ij}}{1+\sum_{\ell \in S} v_{i\ell}}$. The outside option is chosen with probability $\phi_i(0,S)= \frac{1}{1+\sum_{\ell \in S} v_{i\ell}}$. Similarly, the MNL model for a supplier $j \in \S$ is defined as follows:  We use $w_{ji}$ to denote the preference weight of supplier $j$ for customer $i \in \C$. The preference weight for the outside option is normalized to 1. Hence, for any assortment of customers $C \subseteq \C$, the probability that supplier $j$ chooses customer $i \in C$ under the MNL model is given by $\phi_j(i,C)= \frac{w_{ji}}{1+\sum_{k \in C} w_{jk}}$. The outside option is chosen with probability $\phi_j(0,C)= \frac{1}{1+\sum_{k \in C} w_{jk}}$.

In the \fullystatic\ problem, the decision maker must decide the assortment to offer each agent, with all assortments displayed simultaneously. It is clear that in an optimal solution, if a customer $i$ is shown to a supplier $j$, then supplier $j$ should also be shown to customer $i$. This is because a match can only occur if both agents choose each other, making it useless to offer an agent to another without reciprocation. Therefore, the \fullystatic\ problem is reduced to selecting edges $(i,j)$ where $i \in \C$ and $j \in \S$. Selecting an edge $(i,j)$ implies including customer $i$ in the assortment offered to supplier $j$ and vice versa. For each pair $(i,j) \in \C \times \S$, let the binary decision variable $x_{ij}$ indicate whether $(i,j)$ is selected, i.e., $x_{ij}=1$ if and only if $(i,j)$ is selected. Recall that $m$ is the number of suppliers and $n$ is the number of customers. Using this notation, the probability that customer $i$ chooses supplier $j$ given the offered assortment is $\frac{v_{ij} x_{ij}}{1+ \sum_{\ell=1}^m v_{i\ell} x_{i\ell}}$, and the probability that supplier $j$ chooses customer $i$ given the offered assortment is $\frac{w_{ji} x_{ij}}{1 + \sum_{k=1}^n w_{jk} x_{kj}}$.
Therefore, the \fullystatic\ problem can be formulated as
\begin{align}
    \label{eq:main_static}
    \OPT_{\fullyS}=\max_{\mathbf{x}\in\{0,1\}^{n\times m}} \;\;\; \sum_{i=1}^n \sum_{j=1}^m \frac{w_{ji}}{1+ \sum_{k=1}^n w_{jk} x_{kj}}\cdot\frac{v_{ij}}{1 + \sum_{\ell=1}^m v_{i\ell} x_{i\ell}} \cdot x_{ij}.
\end{align}

First, we observe that the \fullystatic\ problem is strongly NP-hard, a conclusion that we can derive directly from the analysis presented in Proposition 1 of \citep{ashlagi_etal19}. Specifically, \citet{ashlagi_etal19} demonstrates the strong NP-hardness of the \Conesidedstatic\  problem using a reduction from the three-partition problem. To prove this result, they construct an instance where the optimal solution involves offering disjoint subsets of suppliers to individual customers (see their Lemma B2). Hence, in this instance, the \Conesidedstatic\ and \fullystatic\ problems are  equivalent. This equivalence arises because, in this instance of the \Conesidedstatic\ problem, each supplier is offered to precisely one customer; thus, irrespective of the customer's decision to select the supplier or not, it is always optimal to offer that customer to the supplier. In other words, adapting to the decisions of customers does not change the expected number of matches. {Consequently, the  \fullystatic\ is also strongly NP-hard.} Given this hardness result, we focus on approximation algorithms for \fullystatic.

{Second, we note that Problem~\eqref{eq:main_static} has an equivalent bilinear formulation with continuous variables (instead of binary) and packing constraints; for more details, we refer to Appendix~\ref{app:fs_bilinear}. As shown in~\citep{el2024lp}, this class of problems is significantly more challenging to approximate within any constant ratio, in general. Despite this hardness of approximation for general bilinear problems, we show the following constant guarantee for \fullystatic:
}
\begin{theorem}\label{thm:fullystatic_mnl_approx}
There exists a {$0.067$}-approximation algorithm for the \fullystatic\ problem when customers and suppliers make choices according to MNL models.
\end{theorem}

To prove Theorem \ref{thm:fullystatic_mnl_approx}, we fix a parameter $\alpha > 0$. In Section \ref{sec:lowlow}, we first consider cases where both suppliers and customers have low values, i.e., $w_{ji} \leq \alpha$ and $v_{ij} \leq \alpha$. In this scenario, we develop a randomized polynomial-time algorithm that yields a $\frac{1}{(2+\alpha)^2}$-approximation for \fullystatic. Then, in Section \ref{sec:high-supplier}, we address cases where either all suppliers have high values ($v_{ij} \geq \alpha$) or all customers have high values ($w_{ji} \geq \alpha$). 
\revised{Since the same argument will apply to either side,}
we focus on the case with high-valued customers ($w_{ji} \geq \alpha$) and provide a randomized algorithm that achieves {$\left(\frac{e-1}{e}\cdot\frac{\alpha}{1+\alpha}\right)$}-approximation for \fullystatic. Building on this analysis, we combine our results to present an algorithm that achieves  {a $0.067$}-approximation for the general case, as demonstrated in Section~\ref{sec:general-proof-static}.

\subsection{Low Value Customers and Low Value Suppliers} \label{sec:lowlow}

In this section, we assume that $w_{ji}\leq \alpha$ and $v_{ij}\leq \alpha$ for all $(i,j)\in \C \times \S $ for some given constant $\alpha>0$ that \revised{will be optimized later.}
We refer to this case as low value customers and low value suppliers. Under this assumption, 
we have the following {lemma}.
\begin{lemma} \label{thm:low_low}
Suppose $w_{ji}\leq \alpha$ and $v_{ij}\leq \alpha$ for all $(i,j)\in \C \times \S $ for some $\alpha>0$. There exists a randomized polynomial time algorithm that gives a  $\frac{1}{(2+\alpha)^2}$-approximation to \fullystatic.
\end{lemma}
To show this result, we present a linear programming (LP)  relaxation \revised{of} Problem \eqref{eq:main_static}. Then, we use the optimal solution of the LP to construct a randomized solution for Problem \eqref{eq:main_static} whose expected objective value is within $\frac{1}{(2+\alpha)^2}$ of the optimal objective value of \eqref{eq:main_static}. 
In particular, consider the following linear program:
\begin{equation} \label{eq:relax_static_simple}
z_{\sf LP}= \max_{\mathbf{y}\geq 0}\left\{ \sum_{i=1}^n  \sum_{j=1}^m  v_{ij}w_{ji} y_{ij}: \
 y_{ij}+\sum_{\ell=1}^m v_{i\ell}y_{i\ell} \leq 1, 
\quad y_{ij}+\sum_{k=1}^n w_{jk}y_{kj} \leq 1, 
\  \forall i\in[n], j \in [m] \right\}.
\end{equation}
\begin{lemma}\label{lemma:upper_bound_static_lowlow}
$  \OPT_{\fullyS} \leq z_{\sf LP}$.
\end{lemma}
The proof of Lemma \ref{lemma:upper_bound_static_lowlow} is given in Appendix \ref{apx:lemma:upper_bound_static_lowlow}.
Now, we give our randomized algorithm for \fullystatic~for the low value setting. Our algorithm is based on a randomized rounding of the optimal solution of the LP given in \eqref{eq:relax_static_simple}. {Consider ${\bf y}^* \in {\mathbb{R}}^{n\times m}$ an optimal solution of \eqref{eq:relax_static_simple}. Note that $0\le y^*_{ij}\le1$.}

\noindent
{\bf Independent Randomized Rounding.} Our proposed randomized solution for \fullystatic\ is the following. For $i\in [n]$ and $j\in [m]$, we define 
\begin{equation} \label{eq:rando}
\tilde{x}_{ij} = \left\{
    \begin{array}{ll}
        1 & \; \; \; \mbox{with probability  } y^*_{ij} \\
        0 & \; \; \; \mbox{with probability  } 1-y^*_{ij}.
    \end{array}
\right.    
\end{equation}
All the variables $\tilde{x}_{ij}$ are independent and well defined because  $ 0 \leq y^*_{ij} \leq 1.$

We show that, in expectation, the randomized solution \eqref{eq:rando} gives $\frac{1}{(2+\alpha)^2}$-approximation to \fullystatic, which enables us to prove Lemma \ref{thm:low_low}. We defer the remaining details of Lemma  \ref{thm:low_low} to Appendix~\ref{apx:low_low_approx}. 

\subsection{High Value Customers} \label{sec:high-supplier}

In this section, we assume that $w_{ji} \geq \alpha$ for all $(i,j) \in \C \times \S$. We refer to this case as suppliers having high values. Under this assumption, we establish the following theorem.

\begin{lemma} \label{thm:high_value_suppliers}
Suppose $w_{ji} \geq \alpha$ for all $(i,j) \in \C \times \S$. There exists a randomized polynomial time algorithm that provides a {$\left(\frac{e-1}{e}\cdot\frac{\alpha}{1+\alpha}\right)$}-approximation to \fullystatic.
\end{lemma}

{To prove Lemma~\ref{thm:high_value_suppliers}, we introduce the following formulation 
\begin{equation} \label{eq:high_value_supp_optimist}
z_{\sf C} = \max \left\{\sum_{i=1}^n \; \frac{\sum_{j=1}^m v_{ij} y_{ij}}{1+\sum_{j=1}^m v_{ij} y_{ij}}: \ 
 \sum_{i=1}^n y_{ij} \leq 1, \ \forall j \in [m], \ 
\;\; y_{ij} \in\{0,1\}, \ \forall i \in [n], \forall j \in [m] \right\},
\end{equation}
which, as the following lemma states, upper bounds Problem \eqref{eq:main_static}. 
\begin{lemma}\label{lemma:upper_bound_static_concave}
$ \OPT_{\fullyS} \leq z_{\sf C}$.
\end{lemma}
The proof of this lemma can be found in Appendix \ref{apx:lemma:upper_bound_static_concave}.
More importantly, Problem~\eqref{eq:high_value_supp_optimist} corresponds to a submodular optimization problem\footnote{The demand function of the MNL choice model is monotone and submodular.} under a partition constraint which can be approximated within $1-1/e$ factor by the well-known \emph{continuous greedy algorithm}~\citep{vondrak2008optimal}. Finally, for high value customers, we can show that any solution of Problem~\eqref{eq:high_value_supp_optimist} with objective value {\sf val} is also a solution of Problem~\eqref{eq:main_static} with objective value at least $\frac{\alpha}{1+\alpha}\cdot {\sf val}$. Therefore, by taking the solution provided by approximately solving~\eqref{eq:high_value_supp_optimist}, we get a $\frac{e-1}{e}\cdot\frac{\alpha}{1+\alpha}$ approximation factor of the main problem, which proves Lemma~\ref{thm:high_value_suppliers}.  The detailed proof is provided in Appendix~\ref{proof:high_value_approx}.}

\subsection{Proof of Theorem \ref{thm:fullystatic_mnl_approx}} \label{sec:general-proof-static}
In this section, we address the general case for any values of $w_{ji}$ and $v_{ij}$ without imposing restrictive assumptions and present the proof of our main result as stated in Theorem \ref{thm:fullystatic_mnl_approx}. We partition the pairs in $\C \times \S$ as follows:
\begin{align*}
    \E^1 &= \{(i,j) \in \C \times \S \mid w_{ji} \geq \alpha\}, \quad
    \E^2 = \{(i,j) \in \C \times \S \mid v_{ij} \geq \alpha, w_{ji} < \alpha\}, \\
    \E^3 &= \{(i,j) \in \C \times \S \mid v_{ij} < \alpha, w_{ji} < \alpha\}.
\end{align*}
The triplet $(\E^1, \E^2, \E^3)$ forms a partition of $\C \times \S$. For each $t \in \{1, 2, 3\}$, we define
\begin{align}
    \label{eq:substatic}
    \OPT^t_{\fullyS} = \max_{\mathbf{x} \in \{0,1\}^{|\E^t|}} \;\;\; \sum_{(i,j) \in \E^t} \frac{w_{ji}}{1 + \sum_{k: (k,j) \in \E^t} w_{jk} x_{kj}}\cdot\frac{v_{ij}}{1 + \sum_{\ell: (i,\ell) \in \E^t} v_{i\ell} x_{i\ell}} \cdot x_{ij}.
\end{align}
{We note that these subproblems satisfy the following property:}
\begin{lemma} \label{lemma:subadditivity}
    $\OPT_{\fullyS} \leq \OPT^1_{\fullyS} + \OPT^2_{\fullyS} + \OPT^3_{\fullyS}$.
\end{lemma}
{We defer the proof of this lemma to Appendix~\ref{proof:subadditivity}.}
{
For Problem \eqref{eq:substatic}, when $t=1$ (case $w_{ji} \geq \alpha$),  Lemma~\ref{thm:high_value_suppliers} provides a randomized  algorithm that achieves a $\frac{e-1}{e}\frac{\alpha}{1+\alpha}$-approximation for $\OPT^1_{\fullyS}$.  By symmetry, the case $t=2$ ($v_{ij} \geq \alpha$) mirrors the case of high value customers and by inverting the roles of suppliers and customers we apply Lemma~\ref{thm:high_value_suppliers} to obtain a $\frac{e-1}{e}\frac{\alpha}{1+\alpha}$-approximation for $\OPT^2_{\fullyS}$. For $t=3$ ($w_{ji} \leq \alpha$ and $v_{ij} \leq \alpha$), we use Lemma~\ref{thm:low_low} to obtain a $\frac{1}{(2+\alpha)^2}$-approximation of $\OPT^3_{\fullyS}$. }

{\noindent {\bf Main Algorithm.} Our method for Problem \eqref{eq:main_static} operates as follows: For each $t \in \{1,2,3\}$, we apply our algorithms to solve Problem~\eqref{eq:substatic} approximately. Let $\Ma^t$ denote the random objective value of our approximate solution for Problem~\eqref{eq:substatic} and $\EE[\Ma^t]$ its expected value. Lemmas~\ref{thm:low_low} and~\ref{thm:high_value_suppliers} guarantee that $\EE[\Ma^t]\geq\frac{e-1}{e}\frac{\alpha}{1+\alpha}\cdot\OPT^t_{\fullyS}$ for $t=1,2$ and $\EE[\Ma^3]\geq\frac{1}{(2+\alpha)^2}\cdot\OPT^3_{\fullyS}$. Our algorithm then selects the solution $t$ with the highest $\Ma^t$ among $t \in \{1,2,3\}$ which we denote by $t^\star$. 
We obtain
\begin{align*}
\OPT_{\fullyS}&\leq \OPT^1_{\fullyS}+\OPT^2_{\fullyS}+\OPT^3_{\fullyS}\\
&\leq \frac{e}{e-1}\cdot\frac{1+\alpha}{\alpha}\cdot\EE[\Ma^1] + \frac{e}{e-1}\cdot\frac{1+\alpha}{\alpha}\cdot\EE[\Ma^2]+(2+\alpha)^2\cdot\EE[\Ma^3]\\
&\leq \left(\frac{2e}{e-1}+\frac{2e}{e-1}\cdot\frac{1}{\alpha}+(2+\alpha)^2\right)\cdot\max\left\{\EE[\Ma^1],\EE[\Ma^2],\EE[\Ma^3]\right\}\\
&\leq \left(\frac{2e}{e-1}+\frac{2e}{e-1}\cdot\frac{1}{\alpha}+(2+\alpha)^2\right)\cdot\EE[\max\left\{\Ma^1,\Ma^2,\Ma^3\right\}]
\end{align*}
where the last step is a direct application of Jensen's inequality. We choose $\alpha >0$ that minimizes  $\left(\frac{2e}{e-1}+\frac{2e}{e-1}\frac{1}{\alpha}+(2+\alpha)^2\right)$, which gives $\alpha \approx 0.7574$. This leads to
\(
\EE[\max\left\{\Ma^1,\Ma^2,\Ma^3\right\}]\geq 0.067\cdot\OPT_{\fullyS}.
\)
We remark that our $0.067$-approximate solution for $\OPT_{\fullyS}$ is constructed by setting $x_{ij}$ equal to our approximate solution for $\OPT^{t^\star}_{\fullyS}$ for $(i,j) \in \E^{t^\star}$ and $x_{ij}=0$ otherwise.
}

%% file: 6_extensions_cardinality.tex
\section{Extensions to Cardinality Constraints}\label{sec:cardinality_constraints}
In this section, we extend Problem~\eqref{eq:general_problem} to the cardinality constrained setting.
Before discussing our results, we first need to address the following question: 
To which side should we constraint the assortment size? This is particularly relevant for platform designs with one-way interactions where constraints can be imposed to both sides or only to the side who gets processed first. For the latter, for example, a carpooling app may restrict the assortment size to passengers looking for a ride, but not limit the requests received by drivers. 
In the remainder of this section, for one-sided policies, we refer to the side who is processed first as the \emph{initiating side} and the one processed second as the \emph{responding side}. 
To the best of our knowledge, the two-sided assortment optimization literature has only focused on one-way interactions where the assortment size is restricted only for the initiating side, either in the online model~\citep{aouad2020online} or offline model~\citep{torrico21}. In the context of dating platforms, \citet{rios2023platform} study settings with constraints on both sides, however, we clarify that they do not consider discrete choice models.

In our constrained framework, each policy must show at most $K_a\in\ZZ_+$ profiles to each processed agent $a\in\A$.
To avoid introducing additional notation, we preserve the names of each policy class, optimization problems and optimal values, but it is understood that they refer to the constrained setting. We distinguish two models:
\begin{itemize}
\item \emph{Two-way constrained:} This setting applies to all four policy classes, where $K_a$ is an arbitrary positive integer for all $a \in \A$. Figure~\ref{fig:adaptivity_gap_constrained} summarizes our main findings—adaptivity gaps and approximation algorithms—for general choice models under this setting. Some of these results improve further when all agents follow MNL choice models. We discuss these results in more detail in the following sections.

\item \emph{One-way constrained:} This setting applies only to one-sided policies, where we assume that on the initiating side, $K_a$ is an arbitrary positive integer, while on the responding side, there is no cardinality constraint, i.e., $K_a = +\infty$. For example, any policy $\pi \in \Pi_{\ConesidedS}$ imposes a cardinality constraint for each customer but not for suppliers. In the following sections, we show how the results from the unconstrained setting extend to this constrained model, under an analogous version of Assumption~\ref{assumption:monotone_submodular}. 
\end{itemize}

\begin{figure}[htpb]
\centering
\scalebox{1}{
 \begin{tikzpicture}[
BC/.style args = {#1/#2/#3}{
        decorate,
        decoration={calligraphic brace, amplitude=6pt,
        pre =moveto, pre  length=1pt,
        post=moveto, post length=1pt,
        raise=#1,
              #2},
         thick,
        pen colour={#3}, text=#3
        },              ]
\draw[thick,-] (0,0) -- (3,0)  {};
\draw[thick,dotted] (3,0) -- (4,0)  {};
\draw[thick,-] (4,0) -- (15,0)  {};
\draw[-] (1,-0.2) -- node[below=2ex] {\footnotesize\fullystatic} (1,0.2);
\draw[-] (6,-0.2) -- node[above=2ex] {\footnotesize\onesidedstatic} (6,0.2);
\draw[-] (10,-0.2) -- node[below=2ex] {\footnotesize\onesidedadaptive} (10,0.2);
\draw[-] (14,-0.2) -- node[above=2ex] {\footnotesize\fullyadaptive} (14,0.2);
 \draw[arrows = {Stealth[scale=1.1]-|},black] (1,0.5) -- node[above=2ex] {\footnotesize Theorem~\ref{thm:fullystatic_mnl_approx_cardinality}: $0.067$ (for MNL)} (1,1.3);
 \draw[arrows = {Stealth[scale=1.1]-|},black] (14,1) -- node[above=2ex] {
 \footnotesize Theorem~\ref{thm:approx_oa_constrained_general}: $\frac{1}{2}(1-1/e)^3$} (14,1.8);
 \draw[arrows = {Stealth[scale=1.1]-|},black] (10,0.5) -- node[above=2ex] {\footnotesize Theorem~\ref{thm:approx_oa_constrained_general}: $(1-1/e)^3$
 } (10,1.3);
 \draw[arrows = {Stealth[scale=1]-|},black] (6,0.9) -- node[above=2ex] {\footnotesize Theorem~\ref{thm:approx_oa_constrained_general}: $(1-1/e)^3$} (6,1.6);
\draw [BC=1mm/mirror/black] 
    (1,-1.1) -- node[below=2ex] {\footnotesize Proposition~\ref{prop:gap_static_zero}: ${\Omega}(n)$} (6,-1.1);
\draw [BC=1mm/mirror/black] 
    (6,-1.1) -- node[below=2ex] {\footnotesize Theorem~\ref{thm:cardinality_gap_os_oa}: $\text{\sc Gap}\in\left[\frac{e}{e-1},\frac{e^2}{(e-1)^2}\right]$} (10,-1.1);
\draw[BC=1mm/mirror/black]    
(10,-1.1) -- node[below=2ex] {\footnotesize Theorem~\ref{thm:gap_fully_adaptive}: $2$} (14,-1.1);
 \end{tikzpicture}}
 \caption{Main results for the two-way constrained setting. The adaptivity gaps are presented below the brackets. On the top, we present our guarantees for each class of problems.
 }
\label{fig:adaptivity_gap_constrained}
\end{figure}
\subsection{Adaptivity Gaps}

We begin by studying the adaptivity gap between \onesidedstatic~and \onesidedadaptive~in the constrained setting. As a first step, we redefine the demand function for each agent. Formally, for each supplier $j \in \S$, we define the \emph{constrained demand function} as:
\begin{equation}\label{def:constrained_supplier_demand}
f^K_j(C) := \max\left\{\sum_{i \in C'} \phi_j(i, C') :\; |C'| \leq K_j,\; C' \subseteq C\right\}, \qquad \forall\; C \subseteq \C.
\end{equation}
The constrained demand function for customers is defined analogously.  
Note that when $K_j = +\infty$, Equation~\eqref{def:constrained_supplier_demand} reduces to the \emph{unconstrained} demand function:
$
f^\infty_j(C) = f_j(C) = \sum_{i \in C} \phi_j(i, C),
$
which, under Assumption~\ref{assumption:monotone_submodular}, is monotone and submodular. This is particularly useful for analyzing the one-way constrained setting, as we can invoke similar arguments to those in Theorem~\ref{thm:gap_onesided_adaptive}.  

In contrast, the two-way constrained model introduces significantly more complexity. In general, the constrained demand function $f^K_j$ remains monotone but is \emph{not submodular} (see a counterexample in Appendix~\ref{app:not_submodular}). Furthermore, it is not clear whether this function can be evaluated efficiently for general choice models.  
Despite these challenges, we show that the adaptivity gap between \onesidedstatic~and \onesidedadaptive~remains bounded by a constant, even in the two-way constrained setting. Specifically:
 For general choice models under the two-way constrained setting, the adaptivity gap lies between $e/(e-1)$ and $(e/(e-1))^2$. When both customers and suppliers follow MNL choice models, the gap improves to exactly $e/(e-1)$ under the two-way constrained setting. In the one-way constrained setting, the gap is also exactly $e/(e-1)$. We summarize these results in the following theorem:
\begin{theorem}\label{thm:cardinality_gap_os_oa}
Under the two-way constrained setting, the adaptivity gap between \onesidedstatic~and \onesidedadaptive~lies between $e/(e-1) \approx 1.582$ and $(e/(e-1))^2 \approx 2.502$. The upper bound improves to $e/(e-1)$ when both customers and suppliers follow MNL choice models. Under the one-way constrained setting, the adaptivity gap is exactly $e/(e-1)$, matching the unconstrained case.
\end{theorem}

To prove the first part of the theorem, we (i) adapt the LP relaxation to incorporate cardinality constraints on both sides, and (ii) leverage a novel connection between \emph{Contention Resolution Schemes}~\citep{chekuri2011submodular,kashaev2023simple} and the \emph{correlation gap} for monotone submodular functions. For the result regarding MNL choice models, we show that the constrained demand function in~\eqref{def:constrained_supplier_demand} is submodular and efficiently computable in polynomial time. Full details of the proof of Theorem \ref{thm:cardinality_gap_os_oa} are deferred to Appendix~\ref{app:proof_cardinality_gap_os_oa}.

\begin{observation}
The proof of the first part of Theorem~\ref{thm:cardinality_gap_os_oa} may be of independent interest. In particular, we show that the adaptivity gap for functions of the form $g(C) = \max\{f(C') : C' \subseteq C,\; |C'| \leq K\}$, where $f$ is monotone and submodular, is upper bounded by $(e/(e-1))^2$. In Appendix~\ref{app:proof_cardinality_gap_os_oa}, we also discuss how this bound extends to more general constraints such as matroids (e.g., partition constraints), and how the bound approaches $e/(e-1)$ in large-market regimes with \revised{large} 
assortments.
\end{observation}

The adaptivity gaps between \fullystatic~and \onesidedstatic, and between \onesidedadaptive~and \fullyadaptive, remain the same as in the unconstrained setting, and follow from essentially identical analyses. In particular, we have the following:

\begin{observation}\label{obs:constrained_gap_fs_fa}
Under the two-way constrained setting, the adaptivity gap between \fullystatic~and \onesidedstatic~is ${\Omega}(n)$, and between \onesidedadaptive~and \fullyadaptive~is 2.
\end{observation}

The results in Remark \ref{obs:constrained_gap_fs_fa} hold because: The lower bound in Proposition~\ref{prop:gap_static_zero} uses assortments of size 1, and thus still applies under cardinality constraints. The result in Theorem~\ref{thm:gap_fully_adaptive} does not rely on Assumption~\ref{assumption:monotone_submodular}. We note that Remark~\ref{obs:constrained_gap_fs_fa} does not extend to the one-way constrained setting, as that setting is only well-defined for \emph{one-sided} policies.

\subsection{Approximation Algorithms}
Before presenting our approximation guarantees, we adapt Assumption~\ref{assumption:monotone_submodular} to the constrained setting:
\begin{assumption}\label{assumption:constrained_oa_oneway}
Under any cardinality constrained setting, we assume:
\begin{enumerate}
  \item {\bf Single-agent Constrained Assortment Optimization Oracle.} For any customer $i\in\C$ and non-negative values $\{\theta_{i,j}\}_{j\in \S}$, there exists polynomial-time  oracle that solves  the \emph{cardinality constrained single-agent assortment optimization problem} $\max\Big\{\sum_{j\in S}\theta_{i,j}\cdot\phi_i(j,S): \ S\subseteq\S, \ |S|\leq K_i\Big\}$. The assumption is analogous for suppliers.
 \item {\bf Monotone Submodular Demand Functions.} 
The \emph{unconstrained} demand function of every agent is a monotone submodular set function.
\end{enumerate}
\end{assumption}

Note that our results in this section extend to settings where the oracle is approximate, as discussed in Remark~\ref{obs:approximate_oracle}. In Theorem \ref{thm:approx_oa_constrained_general}, we present our constant approximation guarantee for one-sided adaptive and fully adaptive policies. To obtain these guarantees, we additionally assume that agents follow weakly substitutable choice models, a standard class in the assortment optimization literature which is verified by all random utility choice models.  Formally, a choice model $\phi$ is \emph{weakly substitutable} if for any subsets $A \subseteq A'$ and any $a \in A \cup \{0\}$, it holds that $\phi(a, A) \geq \phi(a, A')$. This condition is stronger than monotonicity of the unconstrained demand function, since $\phi(0, A) \geq \phi(0, A')$ implies $\sum_{b \in A} \phi(b, A) \leq \sum_{b \in A'} \phi(b, A')$ via the identity $\phi(0, A) = 1 - \sum_{b \in A} \phi(b, A)$.

\begin{theorem}\label{thm:approx_oa_constrained_general}
Under the two-way constrained setting, there exists a polynomial-time algorithm that achieves, in expectation, a $(1-1/e)^3\approx 0.252$ approximation factor for \Conesidedadaptive~and \Conesidedstatic,~when customers and suppliers follow general choice models that satisfy weak-substitutability. Moreover, the same factor holds for its adaptation to \Sonesidedadaptive~and \Sonesidedstatic.
Furthermore, there exists a polynomial-time algorithm that achieves a $\frac{1}{2}\cdot(1-1/e)^3$-approximation for \fullyadaptive. 
\end{theorem}

Theorem~\ref{thm:approx_oa_constrained_general} is the only result in this section that requires the weak substitutability assumption. All the other results below rely only on Assumption~\ref{assumption:constrained_oa_oneway}.
To prove Theorem~\ref{thm:approx_oa_constrained_general}, we adapt the LP relaxation~\eqref{eq:relaxation_onesided_customers} to the two-way constrained setting (assuming customers are the initiating side), and solve it via the ellipsoid method with an approximate separation oracle. Based on the LP solution, we sample an assortment for each customer. After observing customer choices, we obtain an optimal assortment of customers for each supplier by solving $f_j^K(C)$ thanks to the oracle in Assumption~\ref{assumption:constrained_oa_oneway}. The guarantee for \fullyadaptive~follows from arguments similar to Theorem~\ref{coro:approx_fully_adaptive}. The full proof of Theorem~\ref{thm:approx_oa_constrained_general} is provided in Appendix~\ref{app:proof_cardinality_approximation_oa}.

We can further improve the results in Theorem \ref{thm:approx_oa_constrained_general} in special cases as outlined in the following corollary.

\begin{corollary}\label{coro:improved_guarantees_constrained}
The guarantees in Theorem~\ref{thm:approx_oa_constrained_general} improve to $1/2$ for \Conesidedadaptive~and \Sonesidedadaptive, and to $1/4$ for \fullyadaptive~in the following cases:
(1) Under the two-way constrained setting, when customers and suppliers follow MNL choice models.
(2) Under the one-way constrained setting, when agents satisfy Assumption~\ref{assumption:constrained_oa_oneway}, regardless of whether their choice models are weakly substitutable.
\end{corollary}

The key technical ingredient in Corollary~\ref{coro:improved_guarantees_constrained} is the submodularity of the constrained demand function on the responding side. In Case (1), the constrained MNL demand function is submodular (see Lemma~\ref{lemma:MNL_constrained_submodularity} in Appendix~\ref{app:twowayconstrained_mnl}) and efficiently computable. In Case (2), Assumption~\ref{assumption:constrained_oa_oneway} suffices to carry out the analysis. Full details are provided in Appendix~\ref{app:proof_corollary_improved_approx_constrained}.

Our final result in this section extends the approximation guarantee for \fullystatic~to the cardinality-constrained setting:

\begin{theorem}\label{thm:fullystatic_mnl_approx_cardinality}
Under the two-way constrained setting, there exists a polynomial-time algorithm that achieves a $0.067$-approximation for the \fullystatic~problem when customers and suppliers follow MNL choice models.
\end{theorem}

This result follows by adapting the relaxations for the low-value and high-value regimes and applying more sophisticated rounding techniques. See Appendix~\ref{app:proof_cardinality_approximation_fs} for the complete proof.

%% file: 7_numerical_experiments.tex

\section{Numerical Experiments}\label{sec:experiments}
In this section, we conduct a computational study on synthetic instances to empirically evaluate our algorithms and  adaptivity gaps. Our goal is to assess how closely these empirical results align with our worst-case guarantees.

\vspace{1mm}
\noindent
{\bf Experimental Setup.} We consider the unconstrained setting with both sides having the same size, i.e., $m=n$. Customers and suppliers make choices according to MNL choice models. Specifically, the value that a customer $i\in\C$ has over a supplier $j\in\S$ is sampled from the uniform distribution on $[0,1]$, i.e., $v_{ij}\sim~U[0,1]$.  On the other hand, the value that a supplier $j\in\S$ has over a customer $i\in\C$ is sampled from the exponential distribution of parameter 1, i.e., $w_{ji}\sim \text{Exp}(1)$. For each problem size, we sample 20 instances and we report the average result across all solved instances within a time limit of 2 hours, along with the minimum and maximum values. All the experiments where performed on a machine with an Intel(R) i5, 4 cores clocked at $2.40$GHz, and $8$GB of RAM. To solve linear and bilinear programs we use Gurobi $11.0.3$ and to solve convex programs we use CVXPY $1.5.3$.

 \vspace{1mm}
\noindent
{\bf Algorithms and Benchmarks.} We now describe the methods used for each problem.
\begin{itemize}
\item For \fullystatic, we implement our algorithm from Section~\ref{sec:approx_fully_static}, which we denote by $\alg_{\fullyS}$. To compute the optimal value (when computationally possible), we solve an equivalent bilinear reformulation~\eqref{eq:bilinear} for \fullystatic, which is included in Appendix~\ref{app:fs_bilinear}.
\item For \onesidedstatic, we implement a policy that uses the greedy method analyzed in \citep{torrico21}. Namely, we simulate Algorithm~\ref{alg:random_order_greedy} for each side, and select the assortments that lead to the best outcome. In expectation, this leads to a $1/2$ worst-case approximation factor against $\OPT_{\onesidedS}$, as observed by \citet{torrico21}. We do not implement their continuous greedy approach because it is less efficient and provides a marginal objective value improvement.
\item For \onesidedadaptive, we implement Algorithm~\ref{alg:onesidedadaptive_policy}, denoted by $\alg_{\onesidedA}$. We perform 100 simulation runs to estimate the expected total number of matches for each initiating side. Additionally, we compute $\OPT_{\onesidedA}$ using a one-sided dynamic program, analogous to the one described in Appendix~\ref{sec:DP_fullyadaptive}. We also consider a concave relaxation~\eqref{eq:ub_oa},  detailed in Appendix~\ref{app:ub_oa}, which provides an upper bound for \onesidedadaptive\ and is denoted by $\UB_{\onesidedA}$.
\item For \fullyadaptive, we compute the optimal value using the dynamic program described in Appendix~\ref{sec:DP_fullyadaptive}. We denote our algorithm from Theorem~\ref{coro:approx_fully_adaptive} as $\alg_{\fullyA}$. We also design a linear programming relaxation~\eqref{eq:ub_fa}, detailed in Appendix~\ref{app:ub_fa}, which provides an upper bound for \fullyadaptive\ and is denoted by $\UB_{\fullyA}$.
\end{itemize}

\vspace{1mm}
\noindent
{\bf Results.}
In Table~\ref{table:fs}, we report the performance of our algorithm for \fullystatic. We first observe that, on average, the empirical performance of $\alg_\fullyS$ is significantly better in practice. In particular, the performance ratio $\alg_\fullyS/\OPT_{\fullyS}$ averages around 0.80 across all instances, which is substantially higher than the worst-case guarantee of 0.067. Moreover, this ratio appears to improve as the market size increases.
Second, computing the optimal solution for \fullystatic\ becomes increasingly challenging as the instance size grows. Even when using the bilinear reformulation for $\OPT_{\fullyS}$, computational time deteriorates with the size of the market. Specifically, for $n = m = 11$, we observe that 6 out of 20 instances were not solved within the two-hour time limit. In contrast, our algorithm is significantly faster, requiring less than $10^{-2}$ seconds for all instances. This highlights the necessity of approximation algorithms as the problem dimension increases.

\begin{table}[htpb]
{
    \caption{\bf Performance of our algorithm for \fullystatic}\label{table:fs}
    \centerline{
    \scalebox{0.9}{
    \begin{tabular}{cccc|cc}
        \toprule
        \multicolumn{1}{c}{} & \multicolumn{3}{c}{$\alg_\fullyS/\OPT_\fullyS$} & \multicolumn{2}{c}{time $\OPT_\fullyS$} \\ \cmidrule{2-4} \cmidrule{5-6}
                            $m=n$    & min & mean & max &   mean & max  \\ 
                                 \midrule
        2 & 0.60 & 0.83 & 1  & $0.032$ & $0.078$ \\ 
        3 & 0.59 & 0.79 & 1 & $0.049$ & $0.094$\\ 
        4 & 0.52 & 0.76 & 1  & $0.093$ & $0.220$\\ 
        5 & 0.62 & 0.80 & 1  & $0.271$ & $3.483$\\ 
        6 & 0.54 & 0.77 & 0.97  & $0.675$ & 2.848\\ 
        8 & 0.66 & 0.79 & 0.92  & 140.8 & 1930\\ 
        9 & 0.70 & 0.82 & 0.92  & 121.6 & 650\\ 
        10 & 0.73 & 0.85 & 0.92  & 407.1 & 1610\\ 
        11 & 0.72 & 0.85 & 0.92  & 1587 & -\\ 
       \bottomrule
    \end{tabular}}}
    \vspace{1em}
    \centerline{\begin{minipage}[t]{13cm}
    \footnotesize{\emph{Note.} 
    Time performance (in seconds), and ratio between our algorithm and of the optimal solution which was solved up to a MipGap tolerance of $1\%$. We use `-' to indicate that the execution time exceeds the 2 hours limit. For the last row, only $14/20$ computations executed in less than 2 hours; all statistics are computed with respect to this sub-sample.
    }
    \end{minipage}}}
\end{table}

\begin{table}[htpb]
{
    \caption{\bf Results for the adaptivity gaps}\label{table:gaps}
    \centerline{
    \scalebox{0.8}{
    \begin{tabular}{cccc|ccc|ccc|ccc|ccc|ccc}
        \toprule
        \multicolumn{1}{c}{} & \multicolumn{3}{c}{$\OPT_{\onesidedS}/\OPT_{\fullyS}$} & \multicolumn{3}{c}{$\UB_{\onesidedA}/\alg_{\fullyS}$} & \multicolumn{3}{c}{$\OPT_{\onesidedA}/\alg_{\onesidedS}$} & \multicolumn{3}{c}{$\UB_{\onesidedA}/\alg_{\onesidedS}$} & \multicolumn{3}{c}{$\OPT_{\fullyA}/\OPT_{\onesidedA}$} & \multicolumn{3}{c}{$\UB_{\fullyA}/\alg_{\onesidedA}$} \\ \cmidrule{2-19} 
        $m=n$ & min & mean & max & min & mean & max & min & mean & max & min & mean & max & min & mean & max & min & mean & max \\
        \midrule 
        2 & 1.00 & 1.12 & 1.33 & 1.41 & 2.20 & 3.85 & 1.00 & 1.02 & 1.14 & 1.30 & 1.56 & 2.04 & 1.00 & 1.07 & 1.14 & 1.45 & 2.10 & 2.59\\ 
        3 & 1.03 & 1.19 & 1.31 & 1.74 & 2.41 & 3.19 & 1.02 & 1.05 & 1.11 & 1.41 & 1.61 & 1.94 & 1.07 & 1.10 & 1.17 & 1.75 & 2.14 & 2.73\\ 
        4 & 1.01 & 1.25 & 1.38 & 1.72 & 2.65 & 4.29 & 1.03 & 1.10 & 1.16 & 1.36 & 1.62 & 1.86 & 1.06 & 1.11 & 1.19 & 1.63 & 2.17 & 2.58\\ 
        5 & - & - & - & 1.92 & 2.45 & 3.26 & 1.05 & 1.11 & 1.18 & 1.51 & 1.67 & 1.83 & 1.06 & 1.11 & 1.17 & 1.79 & 2.16 & 2.57\\ 
        6 & - & - & - & 1.79 & 2.65 & 4.02 & 1.07 & 1.14 & 1.19 & 1.61 & 1.71 & 1.86 & 1.06 &  1.09 & 1.13 & 1.94 & 2.19 & 2.54\\ 
        7 & - & - & - & 1.90 & 2.34 & 3.12 & 1.12 & 1.17 & 1.23 & 1.59 & 1.69 & 1.83 & - & - & - & 1.94 & 2.16 & 2.49\\ 
        8 & - & - & - & 2.15 & 2.65 & 3.33 & - & - & - & 1.61 & 1.72 & 1.83 & - & - & - & 1.93 & 2.23 & 2.50\\ 
        10 & - & - & - & 2.09 & 2.34 & 2.94 & - & - & - & 1.63 & 1.74 & 1.84 & - & - & - & 1.94 & 2.13 & 2.41 \\ 
        15 & - & - & - & 2.08 & 2.32 & 2.55 & - & - & - & 1.75 & 1.79 & 1.83 & - & - & - & 1.97 & 2.05 & 2.20\\ 
        20 & - & - & - & 2.16 & 2.34 & 2.49 & - & - & - & 1.79 & 1.82 & 1.87 & - & - & - & 1.95 & 2.01 & 2.07\\ 
        \bottomrule
    \end{tabular}}}
     \vspace{1em}
                \centerline{\begin{minipage}[t]{13cm}
                \footnotesize{\emph{Note.} 
                We use `-' to indicate that the execution time exceeds the 2 hours limit. We report the minimum, average and maximum values among the instances that were solved within the time limit.
                }
                \end{minipage}}
                }
\end{table}

In Table~\ref{table:gaps}, we present different ratios to analyze the adaptivity gap between all our policy classes.
We remark that all the optimal values can be computed only for small instances due to the exponentially-sized dynamic programs. Because of this, to analyze larger instances, we use the corresponding upper bounds. Here is a summary of our results for Table~\ref{table:gaps}:

\begin{itemize}
    \item The first column of Table~\ref{table:gaps} reports the gap between \onesidedstatic\ and \fullystatic, given by the ratio $\OPT_{\onesidedS}/\OPT_{\fullyS}$. The average values of this gap ranges between 1.12 and 1.25. 
    This provides empirical evidence that fully static policies can perform well in practice as compared to one-sided static. This observation is further supported by the second column, which reports $\UB_{\onesidedA}/\alg_{\fullyS}$, an upper bound on $\OPT_{\onesidedS}/\OPT_{\fullyS}$. For larger instances, this average upper bound does not exceed 2.7 (the worst-case single instance was for $n=4$ where the bound reaches $4.29$). In particular, these empirical gaps are significantly better than the worst-case bound of $\Omega(n)$ demonstrated in Proposition~\ref{prop:gap_static_zero}. Furthermore,  these  upper bounds that remain tractable for large instances, may be loose; hence, the actual gap should be even smaller than the values in the second column suggest.

    \item The third column, $\OPT_{\onesidedA}/\alg_{\onesidedS}$, provides an upper bound on the gap between \onesidedadaptive\ and \onesidedstatic. We observe that these values are close to 1, substantially improving upon the worst-case bound of $e/(e-1) \approx 1.58$. For larger instances, computing $\OPT_{\onesidedA}$ becomes intractable, so we instead report a more tractable upper bound in the fourth column, $\UB_{\onesidedA}/\alg_{\onesidedS}$. While this bound exceeds the worst-case value, it remains fairly tight, with average values between 1.56 and 1.82. Despite being coarse, these results are still useful for matching platforms where computing optimal values is infeasible, as they offer informative benchmarks.

    \item Finally, the fifth column, $\OPT_{\fullyA}/\OPT_{\onesidedA}$, presents the exact empirical gap between \onesidedadaptive\ and \fullyadaptive\ for small instances. The observed values (averaging between 1.07 and 1.11) are well below the worst-case bound of 2. For larger instances, the sixth column reports the upper bound $\UB_{\fullyA}/\alg_{\onesidedA}$, which exceeds the theoretical bound but remains close to 2. \revised{These empirical findings suggest that even  upper bounds, which might be loose with respect to the optimal values, can offer meaningful empirical adaptivity gaps for larger instances.}
\end{itemize}

\begin{table}[htpb]
{
    \caption{\bf Performance of our algorithm for \fullyadaptive~and \onesidedadaptive}\label{table:adaptiveperf}
    \centerline{
    \scalebox{0.9}{
    \begin{tabular}{cccc|ccc|ccc|ccc}
        \toprule
        \multicolumn{1}{c}{}  & \multicolumn{3}{c}{$\alg_{\onesidedA}/\OPT_{\onesidedA}$} &  \multicolumn{3}{c}{$\alg_{\onesidedA}/\UB_{\onesidedA}$} & \multicolumn{3}{c}{$\alg_{\fullyA}/\OPT_{\fullyA}$} & \multicolumn{3}{c}{$\alg_{\fullyA}/\UB_{\fullyA}$} \\ \cmidrule{2-13} 
        $m=n$ & min & mean & max & min & mean & max & min & mean & max & min & mean & max \\
        \midrule 
        2 & 0.98 & 0.99 & 1 & 0.50 & 0.66 & 0.77 & 0.86 & 0.93 & 1 & 0.39 & 0.49 & 0.69 \\ 
        3 & 0.97 & 0.99 & 1 & 0.51 & 0.65 & 0.74 & 0.84 & 0.90 & 0.93 & 0.37 & 0.47 & 0.57 \\ 
        4 & 0.97 & 0.99 & 0.99 & 0.59 & 0.67 & 0.76 & 0.84 & 0.89 & 0.94 & 0.39 & 0.47 & 0.61\\ 
        5 & 0.97 & 0.98 & 0.99 & 0.59 & 0.65 & 0.70  & 0.83 & 0.89 & 0.93 & 0.39 & 0.47 & 0.56\\ 
        6 & 0.97 & 0.98 & 0.99 & 0.59 & 0.66 & 0.71 & 0.86 & 0.90 & 0.92 & 0.39 & 0.46 & 0.52 \\ 
        7 & 0.97 & 0.98 & 0.99 & 0.62 & 0.67 & 0.71 & - & - & - & 0.40 & 0.46 & 0.51 \\ 
        8 & - & - & - & 0.60 & 0.67 & 0.71 & - & - & - &  0.40 & 0.45 & 0.52\\ 
        10 & - & - & - & 0.63 & 0.68 & 0.71 & - & - & - & 0.41 & 0.47 & 0.52 \\ 
        15 & - & - & - & 0.68 & 0.70 & 0.73     & - & - & - &  0.46 & 0.49 & 0.51\\ 
        20 & - & - & - & 0.70 & 0.71 & 0.73 & - & - & - & 0.48 & 0.50 & 0.51 \\ 
        \bottomrule
    \end{tabular}}}
     \vspace{1em}
                \centerline{\begin{minipage}[t]{13cm}
                \footnotesize{\emph{Note.} 
                We use `-' to indicate that the execution time exceeds the 2 hours limit. We report the minimum, average and maximum values among the instances that were solved within the time limit.
                }
                \end{minipage}}
                }
\end{table}

In Table~\ref{table:adaptiveperf}, we present the performance of our algorithm and benchmarks concerning  \onesidedadaptive~and \fullyadaptive.
As before, we remark that $\OPT_{\onesidedA}$ and $\OPT_{\fullyA}$ can be computed only for small instances due to the exponentially-sized dynamic programs, which is reflected on their time performance that we report in Table~\ref{table:adaptivetimes} (see Appendix~\ref{app:additional_tables}). Therefore, to analyze the performance of our method for larger instances, we use the corresponding relaxations $\UB_{\onesidedA}$ and $\UB_{\fullyA}$. Both of these bounds can be solved efficiently for all the market sizes that we tested (see Table~\ref{table:adaptivetimes} in Appendix~\ref{app:additional_tables}). In that regard, we remark as well that $\alg_{\onesidedA}$ and $\alg_{\fullyA}$ take less than a few seconds to run.
We note that, for smaller instances and against the true optimal values, our method empirically outperforms its worst-case guarantees of $1/2-\epsilon$ and $1/4$. Specifically, the empirical ratio for $\alg_{\onesidedA}$ are (on average) at least 98\% for \onesidedadaptive~when compared to $\OPT_{\onesidedA}$ and at least 65\% relative to $\UB_{\onesidedA}$. Moreover, $\alg_{\fullyA}$'s performance achieves at least 89\% when compared to $\OPT_{\fullyA}$ and at least 45\% relative to $\UB_{\fullyA}$. 
In conclusion, our methods perform well even against these coarser upper bounds.

\revised{We also perform additional experiments to study the performance of our algorithms in large instances and their sensitivity to key features and parameters. We refer the reader to Appendix~\ref{appendix:sensitivity} for all the details.}

%% file: 8_conclusions.tex
\section{Conclusions}\label{sec:conclusions}
In this work, we propose a toolkit to address efficiency and design challenges faced by choice-based matching platforms. Specifically, we introduce a two-sided assortment optimization framework  under general choice models in which a wide range of policies are permitted, ranging from \fullystatic~to \fullyadaptive. To compare different platform designs, we studied the adaptivity gaps between different policy classes and, to address efficiency questions, we provided strong approximation guarantees for \fullyadaptive~and \fullystatic. To the best of our knowledge, our work is the first, from a theoretical standpoint, to precisely measure one-sided and two-sided adaptivity gaps which present a technical departure from the related literature. Practically speaking, our work reveals a natural tension between not only one-sided and two-sided policies but also between non-adaptive and adaptive ones. Specifically, a challenging trade-off is revealed for matching platforms: On the one hand, adaptive policies enjoy stronger guarantees, but might be infeasible to implement because of limited resources, specific markets or it could be expensive to wait for the decision of certain agents. On the other hand, non-adaptive policies might be easier to implement but have worse guarantees; although, they show promising empirical performance.

Several research questions remain open. 
First, our model focuses on optimizing the expected size of the matching between suppliers and customers. It would be intriguing to explore a generalization of the objective function, where the goal is to maximize the total expected revenue of a matching, with a defined revenue value for each pair of customers and suppliers. Also, it would be interesting to explore other operational constraints and forms of interactions \revised{such as repeated interactions or online arrivals}.
Finally, our framework focuses on a bipartite graph of customers and suppliers. Generalizing this framework to non-bipartite graphs, where a match can occur between arbitrary nodes (e.g., for kidney exchange), represents an interesting future direction.

%% file: 999_appendix.tex

\section{DP Formulation for \fullyadaptive}\label{sec:DP_fullyadaptive}
In \fullyadaptive, the number of steps that a policy makes is equal to the number of agents, i.e., $|\A| = n+m$. Let us denote by $T=n+m$ the total number of steps and $t\in[T]$ a step in which a policy decides to process an agent.

Fix a step $t\in[T]$. Let $N_a$ be the agents in the opposite side of $a\in\A$, i.e., for $i\in\C$ we have $N_i=\S$ and for $j\in\S$ we have $N_j = \C$. Denote by $B_a\subseteq N_a$ the subset of agents that have chosen agent $a\in\A$ in a stage before $a$ is processed, which we call backlogs. 
For any subset $A\subseteq \A$ of agents that remain to be processed and their current backlogs $\{B_k\}_{k\in A}$, define $V_t(A,\{B_k\}_{k\in A})$ as the maximum expected number of matches we can obtain from stage $t$ forward (value to-go function). Then, in each time period $t\in[T]$ we have the following Bellman equation in which we decide the agent to be processed:
\begin{align}\label{eq:dp_fullyadaptive}
V_t(A,&\{B_k\}_{k\in A}) = \max_{\substack{a\in A, \\ S\subseteq N_a}}\Bigg\{\sum_{a'\in S\cap B_a}\phi_{a}(a',S)\cdot\Big(1+V_{t+1}(A\setminus\{a\},\{B_k\}_{k\in A\setminus \{a\}})\Big) \notag\\
&\hspace{2em}\quad +\sum_{a'\in S\cap B^c_a}\phi_a(a',S)\cdot\one_{\{a'\in A\}}\cdot V_{t+1}\Big(A\setminus\{a\},\Big\{\{B_k\}_{k\in A\setminus\{a,a'\}}, B_{a'}\cup\{a\}\Big\}\Big) \notag\\
&\hspace{2em}\quad +\phi_a(0,S)\cdot V_{t+1}(A\setminus\{a\},\{B_k\}_{k\in A\setminus \{a\}})
\Bigg\}, 
\end{align}
where the first sum corresponds to $a\in A$ selecting someone in the backlog; in that case, the matches increase by one and we just remove $a$ from non-processed agents. The second sum corresponds to the case when $a$ selects someone that is not in the backlog. In this sum, the only term that could generate a match is when $a$ selects someone that will be processed in the future (in this case we update the backlog of $a'$). The last term corresponds to the case in which $a$ chooses the outside option in which case we do not gain anything.

\section{Missing Proofs in Section~\ref{sec:adaptivity_gaps}}\label{sec:missing_proofs_adaptivitygaps}

\subsection{Proof of Proposition~\ref{prop:gap_static_zero}}
Consider an instance with $n$ customers and $1$ supplier.
For each customer $i \in \{1,\ldots,n\}$, the choice model of $i$ is defined such that, $ \phi_i( 1 , \{1\}) = 1/n $ and $ \phi_i( 0 , \{1\}) = 1-1/n $. The choice model of supplier $1$ is given such that, for  any $C \subseteq\C $, $ \phi_{1}( i , C) = \frac{1}{|C|} $ if  $i \in C$ and $ \phi_{1}( i , C) = 0 $ if $i \notin C$ or $i=0$. It is easy to check that these choice models {are particular cases of MNL choice models which satisfy Assumption~\ref{assumption:monotone_submodular}.}

Let us study  \onesidedstatic. First, for \Conesidedstatic, the only option is to offer supplier $1$ to every customer. Then, supplier $1$ will get chosen by at least one customer with probability $1- (1-1/n)^n$
and supplier $1$ will choose one customer back with probability $1$. For \Sonesidedstatic, we note that the outside option probability is always equal to $0$ for supplier $1$ independently of the assortment $C$ that we offer. So, regardless of the assortment, supplier 1 chooses some customer $i\in C$, and the best decision to do after that is to offer supplier $1$ to customer $i$. Hence we will have a match with probability $1/n$.  Therefore,
$\OPT_{\onesidedS}=\max\{1- (1-1/n)^n,1/n\} = 1- (1-1/n)^n\xrightarrow{n \xrightarrow{}\infty}  1-1/e$. 

Now consider \fullystatic. Suppose we offer an assortment $C$ to supplier $1$ and show supplier $1$ to every customer in the assortment $C$. The probability that we get a match between  a customer in $C$ and supplier $1$ is $\frac{1}{n} \cdot \frac{1}{|C|}$. The probability that we get a match over all is $\sum_{i \in C} \frac{1}{n} \cdot \frac{1}{|C|}=1/n$ for any $C$. Therefore, $\OPT_{\fullyS}=1/n$. We conclude that {\color{black} $\OPT_{\onesidedS} = \Omega (n) \cdot \OPT_{\fullyS}.$ }
\subsection{Proof of Lemma~\ref{lemma:relaxation_onesided_customers}}\label{sec:proof_relaxation_onesided_customers}
We prove Lemma~\ref{lemma:relaxation_onesided_customers} by constructing a feasible solution $\lambda,\tau$ of Problem~\eqref{eq:relaxation_onesided_customers} from the sample paths produced by an optimal policy in \Conesidedadaptive. Let $\pi^\star$ be an optimal policy in \Conesidedadaptive~ and denote by $\Omega$ the set of sample paths of this policy, where a particular sample path is denoted by $\omega$. A sample path can be viewed as a path, from the root (first customer is processed) to a leaf (last customer is processed and makes its choice), in the decision tree defined by $\pi^\star$. Recall that in \Conesidedadaptive\ we are processing customers first one by one, and then we process all of the suppliers. As we discussed in Section~\ref{sec:onesided_adaptive}, in any sample path of $\pi^\star$, the assortments for suppliers will be the set of customers that choose each of them.

For a given sample path $\omega$, let $S_i^\omega$ be the assortment shown to customer $i\in\C$ and $\xi_i^\omega\in S_i^\omega\cup\{0\}$ be the alternative that $i$ chose. Note that when fixing $\omega$, $\xi_i^\omega$ is deterministic. On the other hand, for a given $\omega$, we denote as $C^\omega_j$ the set of customers that chose supplier $j\in\S$.
For $\omega\in\Omega$, we denote as $\PP(\omega)$ the probability of that path, which is equal to
\[
\PP(\omega) = \prod_{i\in\C}\phi_i(\xi_i^\omega,S_i^w).
\]
With this, we define our variables as follows: for all $j\in\S, \ C\subseteq\C$,
\begin{align*}
\lambda_{j,C} = \sum_{\omega\in\Omega: C^\omega_j=C}\PP(\omega) &=    \sum_{\omega\in\Omega: C^\omega_j=C}\prod_{i\in\C}\phi_i(\xi_i^\omega,S_i^w)   \\&= \sum_{\substack{\omega\in\Omega:\\ C^\omega_j=C}}\prod_{i\in C}\phi_i(j,S^\omega_i)\prod_{i\in\C\setminus C}\phi_i(\xi_i^\omega,S^\omega_i),
\end{align*}
and for all $i\in\C, \ S\subseteq\S$,
\[
\tau_{i,S} = \sum_{\omega\in\Omega: S^\omega_i=S}\PP(\omega_{\text{-}i}) = \sum_{\substack{\omega\in\Omega:\\ S^\omega_i=S}}\prod_{k\in\C\setminus\{i\}}\phi_k(\xi_k^\omega,S^\omega_k), 
\]
where $\PP(\omega_{\text{-}i})=\prod_{k\in\C\setminus\{i\}}\phi_k(\xi_k^\omega,S^\omega_k)$ denotes the probability of the  path $\omega$ but excluding the choice of customer $i$ in the decision tree.  This is because the probability that $\pi^\star$ shows $S$ to $i$ does not depend on the choice of $i$.
Let us prove  these variables satisfy the constraints: First, for any $j\in\S$ we have
\begin{align*}
\sum_{C\subseteq\C}\lambda_{j,C} = \sum_{C\subseteq\C}\sum_{\substack{\omega\in\Omega: \\ C^\omega_j=C}}\prod_{i\in\C}\phi_i(\xi_i^\omega,S_i^w) &=\sum_{\omega\in\Omega}\left(\sum_{C\subseteq\C}\one_{\{C^\omega_j=C\}}\right)\prod_{i\in\C}\phi_i(\xi_i^\omega,S_i^w)\\
&=\sum_{\omega\in\Omega}\prod_{i\in\C}\phi_i(\xi_i^\omega,S_i^w) = \sum_{\omega\in\Omega} \PP(\omega)= 1,
\end{align*}
where $\one_{\{C^\omega_j=C\}}$ is the indicator that equals 1 when $C^\omega_j=C$ and zero otherwise. The second equality is a change of the sum operands, the next equality holds because in each path there exists exactly one set of customers equal to $C_j^\omega$. The last equality holds because we are summing over all possible paths. 

We can prove analogously that for any $i\in\C$ we have 
$$
\sum_{S\subseteq\S}\tau_{i,S} = 1.
$$
Now we prove the final constraint that correlates $\lambda$ and $\tau$. For any $i\in \C$ and $j\in\S$, 
\begin{align*}
\sum_{C\subseteq\C: C\ni i}\lambda_{j,C}&=\sum_{\substack{C\subseteq\C: \\C\ni i}}\sum_{\substack{\omega\in\Omega:\\ C^\omega_j=C}}\prod_{k\in\C}\phi_k(\xi_k^\omega,S_k^\omega)\\
&=\sum_{\substack{C\subseteq\C: \\C\ni i}}\sum_{\substack{\omega\in\Omega:\\ C^\omega_j=C}}\phi_i(j,S^\omega_i)\cdot\prod_{k\neq i}\phi_k(\xi_k^\omega,S_k^\omega)\\
&=\sum_{\substack{C\subseteq\C: \\C\ni i}}\sum_{\substack{\omega\in\Omega:\\ C^\omega_j=C}}\left(\sum_{\substack{S\subseteq\S:\\S\ni j}}\phi_i(j,S)\cdot\one_{\{S=S_i^\omega\}}\right)\cdot\prod_{k\neq i}\phi_k(\xi_k^\omega,S_k^\omega)\\
&=\sum_{\omega\in\Omega}\sum_{\substack{C\subseteq\C: \\C\ni i}}\sum_{\substack{S\subseteq\S:\\S\ni j}}\phi_i(j,S)\cdot\one_{\{S=S_i^\omega\}}\cdot\one_{\{C=C_j^\omega\}}\cdot\prod_{k\neq i}\phi_k(\xi_k^\omega,S_k^\omega)\\
&=\sum_{\omega\in\Omega}\sum_{\substack{S\subseteq\S:\\S\ni j}}\phi_i(j,S)\cdot\one_{\{S=S_i^\omega\}}\left(\sum_{C\subseteq\C\setminus\{i\}}\one_{\{C_j^\omega=C\cup\{i\}\}}\right)\cdot\prod_{k\neq i}\phi_k(\xi_k^\omega,S_k^\omega)\\
&=\sum_{\substack{S\subseteq\S:\\S\ni j}}\sum_{\substack{\omega\in\Omega:\\ S_i^\omega = S}}\phi_i(j,S)\cdot\PP(\omega_{\text{-}i})\\
&=\sum_{S\subseteq\S: S\ni j}\phi_i(j,S)\cdot\tau_{i,S},
\end{align*}
where in the second to last equality we used that in each path there is be exactly one set of customers equal to $C^\omega_j$.

For a fix path $\omega\in\Omega$, the probability that supplier $j\in\S$ chooses back one of the customers that chose her initially is $f_j(C_j^\omega)$. Therefore, the objective value achieved by $\pi^\star$ corresponds to
\begin{align*}
\sum_{\omega\in\Omega}\left(\sum_{j\in\S}f_j(C_j^\omega)\right)\cdot\PP(\omega) &= \sum_{\omega\in\Omega}\left(\sum_{j\in\S}\sum_{C\subseteq\C}f_j(C)\cdot\one_{\{C=C_j^\omega\}}\right)\cdot\PP(\omega)\\
&=\sum_{j\in\S}\sum_{C\subseteq\C}f_j(C)\sum_{\omega\in\Omega}\one_{\{C=C_j^\omega\}}\cdot\PP(\omega)\\
&=\sum_{j\in\S}\sum_{C\subseteq\C}f_j(C)\sum_{\substack{\omega\in\Omega:\\ C_j^\omega =C}}\PP(\omega)\\
&=\sum_{j\in\S}\sum_{C\subseteq\C}f_j(C)\cdot\lambda_{j,C}.
\end{align*}

\subsection{Proof of Lemma~\ref{lemma:UB_onesided_adaptive}}

Consider an instance with $n$ customers and $n$ suppliers. On the one hand, all customers have the same choice model, such that for any $i \in \C$,  any subset of suppliers $S \subseteq\S$, $ \phi_i( j , S) = 1/|S| $ for all $j \in S$ and $ \phi( 0 , S) = 0 $. For any $k \in [n]$, define the sequence $\beta_k=k\cdot(1-e^{-\frac{1}{k}})$ and $\beta_0=0$. Note that $0 \leq \beta_k \leq 1$ .
On the other hand, all suppliers have the same choice model: for any $j \in \S$,  any subset of customers $C \subseteq\C$ with $|C|=k$, we consider  $ \phi_j( i , C) = \beta_k/k $ for all $i \in C$ and $ \phi( 0 , C) = 1-\beta_k $. It is easy to see that all these choice models are monotone and submodular, satisfying the second point in Assumption~\ref{assumption:monotone_submodular}. 

First, let us compute the optimal solution of \onesidedstatic. For this, consider the \Conesidedstatic~problem and denote by $\Ma_1$ the number of matches obtained by showing a family of assortments $S_1,\ldots,S_n\subseteq\S$ where for each $i \in \C $, $S_i$ is a feasible assortment offered to customer  $i$. After customers make their choices, for each $j \in \S$, let $\alpha_j \in \{0,1,\ldots,n\}$ be the random number of customers who choose supplier $j$. Since the remaining unmatched probability in the customers' choice models is $0$, then  $\sum_{j=1}^n \alpha_j=n$. Therefore, conditioned on $\alpha_1,\ldots,\alpha_n$:
\[\EE[\Ma_1|\alpha_1,\ldots,\alpha_n]  =  \sum_{j=1}^n \beta_{\alpha_j} 
=     \sum_{j=1}^n     \alpha_j\cdot(1-e^{-\frac{1}{\alpha_j}}) = n - \sum_{j=1}^n     \alpha_j\cdot e^{-\frac{1}{\alpha_j}},
\]
where in the first equality we use that supplier $j$ sees $\alpha_j$ customers and the demand equals $f_j(C_j)=\beta_{\alpha_j}$.
Thus, by applying Jensen's inequality, we get
\[
\frac{1}{n} \sum_{j=1}^n     \alpha_j\cdot e^{-\frac{1}{\alpha_j}} \geq  e^{ - \frac{1}{n} \sum_{j=1}^n     \alpha_j\cdot\frac{1}{\alpha_j}}  = e^{-1}, 
\]
which implies that
$\EE[\Ma_1|\alpha_1,\ldots,\alpha_n] \leq  n- n e^{-1} $. Therefore, by taking expectation over $\alpha_1,\ldots,\alpha_n$ 
we conclude that $\EE[\Ma_1]\leq n\cdot(1-1/e)$.
Now, consider \Sonesidedstatic~and denote by $\Ma_2$ the number of matches obtained under a family of assortments $C_1,\ldots,C_n\subseteq\C$ where for each $j \in S $, $C_j$ is the feasible assortment offered to supplier $j$ with $|C_j|=k_j$. Since the probability of remaining unmatched in the choice model of customers (the responding side) is always $0$, then the number of matches is equal to the number of customers who are chosen after suppliers are processed, i.e.,
$\Ma_2 = \sum_{i=1}^n \mathbbm{1}_{\{i \text{ is chosen}\}}$. Therefore, $\EE[\Ma_2] = n- \sum_{i=1}^n  \mathbbm{P}( i \text{ is not chosen}) $. Moreover, we have for any $i \in \C$,
\begin{align*}
    \frac{1}{n}\sum_{i=1}^n\mathbbm{P}( i \text{ is not chosen}) =  \frac{1}{n}\sum_{i=1}^n\prod_{\substack{j \in \S:\\ i \in C_j}}  \left( 1- \frac{ \beta_{k_j}}{k_j}    \right)   
     =  \frac{1}{n}\sum_{i=1}^n\exp \Bigg(  - \sum_{\substack{j \in \S:\\ i \in C_j}}  \frac{1}{k_j} \Bigg), 
\end{align*}
where the first equality is because $\mathbbm{P}( i \text{ is not chosen})$ is equal to the product of the probabilities that each $j$ did not choose $i$, for each $j\in\S$ that sees $i$. The following equality is due to the definition of $\beta_{k_j}$. Finally, by Jensen's inequality, we can prove that 
 \[
 \frac{1}{n} \sum_{i=1}^n   \exp \left(  - \sum_{ j \in \S : i \in C_j}  \frac{1}{k_j} \right) \geq 
    \exp \left(  -  \frac{1}{n} \sum_{i=1}^n \sum_{ j \in \S : i \in C_j}  \frac{1}{k_j} \right)
 =  \exp \left(  -  \frac{1}{n} \sum_{j=1}^n \sum_{ i =1}^n \mathbbm{1}_{\{i \in C_j\}}\cdot  \frac{1}{k_j} \right)
 = e^{-1},
 \]
which implies that $\EE[\Ma_2] \leq  n- n e^{-1}$.
Given this, we conclude $\OPT_{\onesidedS} \leq  n\cdot(1-1/e)$. 
To show that is equal to $n\cdot(1-1/e)$, we now provide a feasible solution.  Consider the following solution for \Conesidedstatic: each customer is offered exactly and exclusively one supplier. Then, each customer will choose the supplier with probability $1$, and the supplier will choose back that customer with probability $\beta_1=1-1/e$. Hence, the expected number of matches is $n\cdot(1-1/e)$ which implies that $\OPT_{\onesidedS} =  n\cdot(1-1/e).$  

To get an upper bound of the adaptivity gap we just need a feasible policy for the \onesidedadaptive.~Consider the following: We start processing customers in any random order one by one. For each customer, we offer all suppliers who are not chosen yet. In particular, the first customer is offered all suppliers in $\S$ and makes a choice. If the customer chooses a supplier $j$, we update the list of suppliers who are not chosen yet to $\S\setminus\{j\}$ and proceed to the next customer. Recall that all customers and suppliers in our instance are similar. Let $\Ma_{p,q}$ be the number of matches under our policy when we have a sub-instance with $p$
 customers and $q$ suppliers. By conditioning on the choice of first customer, we get
\[
   \EE [\Ma_{n,n}] = (1+ \EE [\Ma_{n-1,n-1}] )\cdot \beta_n +  \EE [\Ma_{n,n-1}] \cdot (1- \beta_n ).
\]
Since adding one customer can only make the matching better then we have $\EE[\Ma_{n,n-1}] \geq \EE [\Ma_{n-1,n-1}]$, which implies that
\( \EE[\Ma_{n,n}] \geq \beta_n + \EE[\Ma_{n-1,n-1}].\)
Hence, by induction, we get 
\[
\OPT_{\onesidedA} \geq \EE[\Ma_{n,n}] \geq  \sum_{k=1}^n \beta_k = \sum_{k=1}^n k\cdot(1-e^{-1/k}) \geq \sum_{k=1}^n k\cdot \left( \frac{1}{k} - \frac{1}{2k^2}\right)=n-\frac{1}{2}\cdot \sum_{k=1}^n \frac{1}{k},
\]
where in the third inequality we use $e^{-x} \leq 1-x+x^2/2$. In other words, we have that $ \OPT_{\onesidedA} \geq  n-\mathcal{O}(\log n)$. Finally, we conclude that 
{\color{black}
\[
\frac{\OPT_{\onesidedA} }{\OPT_{\onesidedS}} \geq \frac{e}{e-1}\cdot(1-  \mathcal{O}(\log n /n)   ). 
\]  }
\subsection{Proof of Lemma~\ref{lemma:UB_fully_adaptive}}
The proof considers 2 ``inverted'' copies of the instance given in the proof of Proposition~\ref{prop:gap_static_zero}. Formally, consider an instance with $n$ customers and $n$ suppliers. The choice models of customers are defined as follows.
For any customer $i \in \{2,\ldots,n\}$, for any subset of suppliers $S \subseteq\S$, if $S$ contains supplier $1$, then $ \phi_i( 1 , S) = 1/n$, $\phi_i(0,S)=1-1/n$ and $ \phi_i( j , S) = 0 $ for any $j \neq 1$. If $S$ does not contain supplier $1$, then $ \phi_i( j , S) = 0 $ for any $j \in S$. In other words, customers $i \in \{2,\ldots,n\}$ choose supplier $1$ with probability $1/n$ if it is offered, otherwise they remain unmatched. The choice model of customer $1$ is such that: 
\begin{itemize}
\item For  any $S \subseteq\S$ such that $S\neq\emptyset$ and $ 1 \notin S$ (supplier $1$ is not in $S$), $ \phi_{1}( j , S) = 1/|S| $  for all $j\in S$.
\item For  any $S \subseteq\S$ such that $S\neq\emptyset$ and $1\in S$, $ \phi_{1}( j , S) = 1/(|S|-1) $  for all $j\in S\setminus\{1\}$ and $ \phi_{1}( 1 , S) = 0 $.
\item For $S=\emptyset$, $\phi_1(0,S) = 1$.
\end{itemize}
Similarly, the choice models of suppliers are defined as follows.
For any supplier $j \in \{2,\ldots,n\}$, for  any subset of customers $C \subseteq\C$, if $C$ contains customer $1$, then $ \phi_j( 1 , C) = 1/n$, $\phi_j(0,C)=1-1/n$ and $ \phi_j( i , C) = 0 $ for any $i \neq 1$. If $C$ does not contain customer $1$, $ \phi_j( i, C) = 0 $ for any $i \in C$. The choice model of supplier $1$, we define it analogously to the one of customer 1, i.e.:
\begin{itemize}
\item For  any $C \subseteq\C$ such that $C\neq\emptyset$ and $ 1 \notin C$ (customer $1$ is not in $C$), $ \phi_{1}( i , C) = 1/|C| $  for all $i\in C$.
\item For  any $C \subseteq\C$ such that $C\neq\emptyset$ and $1\in C$, $ \phi_{1}( i , C) = 1/(|C|-1) $  for all $i\in C\setminus\{1\}$ and $ \phi_{1}( 1 , C) = 0 $.
\item For $C=\emptyset$, $\phi_1(0,C) = 1$.
\end{itemize}
Clearly, these choice models have a monotone submodular demand function. 

First, let us compute the optimal solution of \onesidedadaptive. For this, we study first \Sonesidedadaptive. Since the set of customers feasible to match with supplier $1$ does not intersect with the set of customers feasible to match with  suppliers $2,\ldots,n$, it doesn't matter which supplier we process first. Assume, without loss of generality, that we start with supplier $1$.
We note that, for this supplier, the outside option probability is always equal to $0$ independently of the assortment $C$ that we offer. So, regardless of the assortment, supplier 1 chooses some customer $i'\in C$. For suppliers $j\in\{2,\ldots,n\}$, the only option is to offer customer $1$. Then, the probability that at least one supplier in $\{2,\ldots,n\}$ choose customer $1$ is $1-(1-1/n)^{n-1}$. Now, we process the customers' side. Show supplier $1$ to customer $i'$ which results in a match with probability $1/n$. For customer $1$, show some supplier $j\in\{2,\ldots,n\}$ that initially chose customer $1$, then a match occurs with
probability $1-(1-1/n)^{n-1}$. Therefore, the total number of matches is $1-(1-\frac{1}{n})^{n-1}+\frac{1}{n}$ which tends to $1-1/e$ as $n$ goes to infinity. A similar analysis can be done for \Conesidedadaptive.

Now let us present a feasible policy of \fullyadaptive. Show customer 1 to every supplier $j\in\{2,\ldots,n\}$ and supplier $1$ to every customer $i\in\{2,\ldots,n\}$. Customer 1 gets chosen by at least one supplier with probability $1-(1-1/n)^{n-1}$ and, similarly, for supplier 1. Show the correspond supplier to customer 1 and the corresponding customer to supplier 1. The expected number of matches is $2\cdot\left[1-(1-1/n)^{n-1}\right]$ which tends to $2\cdot(1-1/e)$ as $n$ increases.
\color{black}
\subsection{Proof of Theorem~\ref{thm:gap_os_fa}}\label{apx:os_fa_gap}

We know that the adaptivity gap cannot be more than $2e/(e-1)$ due to Theorems~\ref{thm:gap_onesided_adaptive} and~\ref{thm:gap_fully_adaptive}. We now provide a family of tight examples.
To ease the analysis, we first analyze a sub-instance with $n$ customers and $n$ suppliers. Our final example (with $2n$ customers and $2n$ suppliers) will be comprised by two ``copies'' of this sub-instance where in the second copy we exchange the roles of customers and suppliers.  

In the sub-instance, all $n$ customers have the same choice model, which is an MNL choice model where every supplier has a preference weight $1/\sqrt{n}$. On the other hand, all $n$ suppliers have the same choice model: Assuming they are offered at least one customer, they choose uniformly at random among the displayed customers.
First, let us analyze \onesidedstatic.
For \Conesidedstatic, we denote as $\Ma_\C$ the random variable counting the number of matches in that setting. Notice that every supplier that gets at least one offer will be matched. We denote as $r_i$ the number of suppliers offered to customer $i$. We can write
\begin{align*}
    \EE[\Ma_\C]=n-\sum_{j\in\S}\PP(j\mbox{ is not chosen})=n\left[1-\frac{1}{n}\sum_{j\in\S}\PP(j\mbox{ is not chosen})\right].
\end{align*}
We now have
\begin{align}
    \frac{1}{n}\sum_{j\in\S}&\PP(j\mbox{ is not chosen})=\frac{1}{n}\sum_{j\in\S}\prod_{\substack{i\in\C\\j\in S_i}}\left(1-\frac{1}{\sqrt{n}+r_i}\right)\notag=\frac{1}{n}\sum_{j\in\S}\prod_{\substack{i\in\C\\j\in S_i}}\exp\left[\log\left(1-\frac{1}{\sqrt{n}+r_i}\right)\right]\notag\\
    &=\frac{1}{n}\sum_{j\in\S}\exp\left[\sum_{\substack{i\in\C\\j\in S_i}}\log\left(1-\frac{1}{\sqrt{n}+r_i}\right)\right]\ge\exp\left[\frac{1}{n}\sum_{j\in\S}\sum_{\substack{i\in\C\\j\in S_i}}\log\left(1-\frac{1}{\sqrt{n}+r_i}\right)\right]\label{eq:jensen2}\\
    &=\exp\left[\frac{1}{n}\sum_{i\in\C}r_i\log\left(1-\frac{1}{\sqrt{n}+r_i}\right)\right]\ge\exp\left[-\frac{1}{n}\sum_{i\in\C}\frac{r_i/(\sqrt{n}+r_i)}{1-1/(\sqrt{n}+r_i)}\right]\label{eq:ineq_log}\\
    &=\exp\left[-\frac{1}{n}\sum_{i\in\C}\frac{r_i}{\sqrt{n}+r_i-1}\right]\ge\exp\left[-\frac{1}{n}\sum_{i\in\C}\frac{n}{\sqrt{n}+n-1}\right]\label{eq:increasing}\\
    &=\exp\left[-\frac{n}{n+\sqrt{n}-1}\right]=\frac{1}{e}-\mathcal{O}\left(\frac{1}{\sqrt{n}}\right).   \label{eq:linear approx}
\end{align}
Inequality~\eqref{eq:jensen2} is due to Jensen's inequality on the exponential function. The next equality is because each customer $i$ is offered exactly $r_i$ suppliers. Inequality~\eqref{eq:ineq_log} is a consequence of $\log(1-x)\ge-x/(1-x)$. Inequality~\eqref{eq:increasing} follows since the function $x\mapsto x/(\sqrt{n}+x-1)$ is increasing and that the $r_i$ are bounded by $n$. Finally, Equality~\eqref{eq:linear approx} is the linear approximation for the exponential function. Therefore, we conclude that $\EE[\Ma_\C]\leq \left(1-\frac{1}{e}\right)\cdot n+\mathcal{O}(\sqrt{n}).$

Now, let us look at \Sonesidedstatic. We denote as $\Ma_\S$ the random variable counting the number of matches in that setting, and as $k_i$ the random variable counting the number of proposals obtained by customer $i$. Because there are only $n$ suppliers, we have that $k_1+...+k_n\le n$. Taking the expected value (conditioned on $k_i$'s) with the information over the choices of the suppliers, we have
\begin{align*}
    \EE[\Ma_\S|k_1,\ldots,k_n]=\sum_{i\in\C}\frac{k_i/\sqrt{n}}{1+k_i/\sqrt{n}}\le\sum_{i\in\C}\frac{k_i}{\sqrt{n}}\le\sqrt{n}.
\end{align*}
We conclude that $\EE[\Ma_\S]\leq \sqrt{n}$ by taking the expectation over $(k_1,...,k_n)$. Therefore, we have that $\OPT_{\onesidedS}=\max\{\OPT_{\ConesidedS},\OPT_{\SonesidedS}\}\le(1-\frac{1}{e})\cdot n+\mathcal{O}(\sqrt{n})$.

For the sub-instance, we would like to get a lower bound on the adaptivity gap between \fullyadaptive~and \onesidedstatic, so we will design a feasible one-sided adaptive policy which, consequently, is feasible for \fullyadaptive. Consider the following: we start processing customers in order from $1$ to $n$. To each customer, we offer all remaining suppliers. We denote as $\Ma_{p,q}$ the expected number of matches obtained with this policy with $p$ customers and $q$ suppliers; we see that $\Ma_{p,q}$ equals the number of picked suppliers. Let us denote $\beta_i=i/(\sqrt{n}+i)$. By conditioning on the choice of the first customer, we get
\begin{align*}
    \EE[M_{n,n}]=\beta_n(1+\EE[M_{n-1,n-1}])+(1-\beta_n)\EE[M_{n-1,n}].
\end{align*}
Since adding a supplier can only make the matching better, we have $\EE[M_{n-1,n-1}]\le \EE[M_{n-1,n}]$. Hence, by induction we get
\begin{align*}
    \OPT_{\fullyA}\ge\EE[M_{n,n}]\ge\sum_{i=1}^n \beta_i=\sum_{i=1}^n \frac{i}{\sqrt{n}+i}= n-\sqrt{n}\sum_{i=1}^n\frac{1}{\sqrt{n}+i}= n-\mathcal{O}(\log(n)\sqrt{n}).
\end{align*}
We are now ready to analyze our final example. For this, we consider an instance with $2n$ customers and $2n$ suppliers. Customers in $\{1,...,n\}$ can only match with suppliers in $\{1,...,n\}$; they have the same choice models as explained in the sub-instance above. For customers and suppliers in $\{n+1,...,2n\}$, we proceed similarly as above, but we exchange their choice models, i.e., suppliers now follow the MNL choice model with weights $1/\sqrt{n}$.
For \onesidedstatic, no matter which side you start from, it is impossible to achieve an outcome better than $\sqrt{n}+(1-1/e)\cdot n+\mathcal{O}(\sqrt{n})=(1-1/e)\cdot n+\mathcal{O}(\sqrt{n})$. For \fullyadaptive, by first processing customers in $\{1,...,n\}$ and then suppliers in $\{n+1,...,2n\}$, we can achieve $2n-\mathcal{O}(\log(n)\sqrt{n})$.
We conclude that
$$
    \frac{\OPT_\fullyA}{\OPT_\onesidedS}\geq\frac{2e}{e-1}-\mathcal{O}\left(\frac{\log(n)\sqrt{n}}{n}\right).  
$$

\color{black}
\section{Missing Proofs in Section~\ref{sec:approx_fully_adaptive}}
\subsection{Proof of Corollary~\ref{coro:simple_approx_factor}}\label{sec:proof_corollary_simple_approx_factor}
\citet{torrico21} show that a variant of the continuous greedy algorithm achieves a $1-1/e$ factor for \Conesidedstatic. Note that this algorithm can be easily adapted to \Sonesidedstatic, we just have to assume that the matching process is initiated by suppliers. This is possible since our setting optimizes the expected size of the matching and not the revenue. Therefore, the following policy guarantees a $1-1/e$ factor to \onesidedstatic: Run the algorithm in \citep{torrico21} with $\C$ being the initiating side and, separately, with $\S$ being the initiating side, then output the solution that achieves the highest objective value between the two. Let us call those random variables as $\Ma_\C$ and $\Ma_\S$. Let us denote the expected value of this method by $\textsc{ALG} = \EE[\max\{\Ma_\C,\Ma_\S\}]$. By using the guarantee in \citep{torrico21} and Jensen's inequality, we obtain:
\begin{align*}
\textsc{ALG}&\geq\max\{\EE[\Ma_\C],\EE[\Ma_\S]\}\geq(1-1/e)\cdot\max\left\{\OPT_{\ConesidedS},\OPT_{\SonesidedS}\right\} = (1-1/e)\cdot\OPT_{\onesidedS} \\ &\geq (1-1/e)^2\cdot \OPT_{\onesidedA},
\end{align*}
where the last inequality follows from Lemma \ref{lemma:LB_onesided_adaptive}.
Finally, the guarantee for \fullyadaptive~follows from Lemma \ref{lemma:LB_fully_adaptive} which states that $\OPT_{\onesidedA} \geq \frac{1}{2}\cdot\OPT_{\fullyA}$. 
\color{black}
\subsection{Proof of Lemma~\ref{lemma:sampling_greedy}}\label{appendix:proof_lemma_sampling}
We start by defining for all $i\in\C$ and $j\in\S$,
\begin{equation*}p_{ij}\coloneqq\max_{S\subseteq\S} \; \; \phi_i(j,S).
\end{equation*}
Note that this definition simplifies to $p_{ij}=\phi_i(j,\{j\})$ when the choice model $\phi_i$ is weak-substitutable, i.e., for all $j\in A\subseteq B\subseteq\S$, we have $\phi_i(j,A)\ge\phi_i(j,B)$. 
Using notation $p_{ij}$, recall the definition of $\phi_{\min}$ as 
\[
\phi_{\min} = \min\left\{p_{ij}:\; i\in\C,\; j\in\S,\; p_{ij}>0\right\},
\]
i.e., the smallest non-zero $p_{ij}$ over all $i$ and $j$. We also define
\begin{equation*}
M^\star := \max\left\{\sum_{j\in\S} f_j(C_j):\ \{C_j\}_{j\in\S}\text{ form a partition of } \C\mbox{ and }\forall j\in\S,\; C_j\subseteq\{i\in\C:\; p_{ij}>0\}\right\},
\end{equation*}
We will use $M^\star$ as an upper bound on the number of matches achieved by Algorithm~\ref{alg:random_order_greedy} in each sample run. 

Fix an arbitrary order over the customers, $\epsilon>0$, $\delta\in(0,1)$ and $T$ independent runs of Algorithm~\ref{alg:random_order_greedy}, which we will determine later. Let us denote by $C_j^t$ the observed backlog of supplier $j\in\S$ and by $M_t$ the number of matches in the $t$-th run.
We would like to estimate the expected total number of matches, 
which we denote by $\mu = \EE_\pi\left[\sum_{j\in\S} f_j(C^\pi_j)\right]$ where the expectation is over the customers' choices which are dependent on the assortments crafted by  Algorithm~\ref{alg:random_order_greedy} denoted by $\pi$. Because of feasibility we have that $\mu\leq \OPT_{\ConesidedA}$. The sample average $\hat\mu$ over $T$ sample runs is defined as
\[
\hat\mu\coloneqq\frac{1}{T}\sum_{t=1}^T M_t=\frac{1}{T}\sum_{t=1}^T\sum_{j\in\S}f_j(C_j^t).
\]
We use the notation $[T]=\{1,\ldots,T\}$. Observe that for any $t\in[T]$, if $i\in C_j^t$, then this means that $p_{ij}>0$, otherwise $i$ would not choose $j$. This implies that $M_t\leq M^\star$. To prove our bound on $T$, we use the following multiplicative form of the Chernoff bound:
\begin{lemma}[see e.g., \citet{dubhashi2009concentration}]\label{lemma:hoeffding}
Let $X_1,\ldots, X_T$ be independent random variables such that $X_t\in[0,1]$ almost surely and $\theta= \EE\Big[\sum_{t\in[T]}X_t\Big]$. Then, for all $\epsilon\in(0,1)$,
\[
\PP\Big(\Big|\sum_{t\in[T]}X_t - \theta\Big|\geq \epsilon\cdot\theta \Big)\leq 2\exp\left(-\theta\epsilon^2/3\right)
\]
\end{lemma}
Consider $X_t = M_t/M^\star$ which  is in $[0,1]$ as we observed before. Also, note that $\EE[X_t] = \mu/M^\star$ and $\EE\Big[\sum_{t\in[T]}X_t\Big] = T\mu/M^\star$. Then, in Lemma~\ref{lemma:hoeffding}, we get
\[
\PP\Big(\Big|\frac{1}{T}\sum_{t\in[T]}\frac{M_t}{M^\star} - \frac{\mu}{M^\star}\Big|\geq \epsilon\cdot\frac{\mu}{M^\star}\Big)\leq 2\exp\left(-\frac{2T\epsilon^2}{3}\cdot\frac{\mu}{M^\star}\right),
\]
which is equivalent to
\[
\PP\Big(\big|\hat\mu - \mu\big|\geq \epsilon\cdot\mu\Big)\leq 2\exp\left(-\frac{2T\epsilon^2}{3}\cdot\frac{\mu}{M^\star}\right).
\]
Since $\mu\leq\OPT_{\ConesidedA}$, the above inequality implies
\[
\PP\Big(\big|\hat\mu - \mu\big|\geq \epsilon\cdot\OPT_{\ConesidedA}\Big)\leq 2\exp\left(-\frac{2T\epsilon^2}{3}\cdot\frac{\mu}{M^\star}\right).
\]
Now, to impose the right hand side to be less or equal than $\delta\in(0,1)$, we need
\[
T \geq \frac{M^\star}{\mu}\cdot \frac{3}{2\epsilon^2}\cdot\log\frac{2}{\delta}.
\]
Our goal now is to find an upper bound for $M^\star/\mu$.
\begin{claim}\label{claim:bound_phi_min}
We have that
$\mu\geq \frac{\phi_{\min}}{2}\cdot M^\star$.
\end{claim}
Before proving this claim, note that the proof of Lemma~\ref{lemma:sampling_greedy} follows immediately. Namely, by taking $T = \frac{3}{\epsilon^2\phi_{\min}}\cdot\log\frac{2}{\delta}$, we get
\[
\PP\Big(\big|\hat\mu - \mu\big|\leq \epsilon\cdot\OPT_{\ConesidedA}\Big)\geq 1-\delta.
\]

\begin{observation}\label{remark:polynumber_samples}
Note that instead of the bound stated in Claim~\ref{claim:bound_phi_min}, we could use that $M^\star \leq \min\{m,n\}$, i.e., the number of matches is at most the minimum size between both sides. Therefore, we can take $T = \frac{\min\{m,n\}}{\mu}\cdot\frac{1}{2\epsilon^2}\cdot\log\frac{2}{\delta}$, so any kind of practical lower bound on $\mu$ (for example $\mu\geq1$ if there exists a pair of agents that are ``easy to match'') would ensure a polynomial number of samples instead of pseudo-polynomial.
\end{observation}

\begin{proof}[Proof of Claim~\ref{claim:bound_phi_min}.]
To prove the claim, we construct a feasible policy using an optimal solution of $M^\star$. Denote by $C_1^\star,\ldots, C_{|\S|}^\star$ an optimal partition for $M^\star$ and, for every $i\in\C$, denote $j_i$ the supplier that verifies $i\in C_{j_i}$. Define the one-sided static policy that shows $S_i = \argmax_{S\subseteq\S}\phi_i(j_i,S)$ for each $i\in\C$. The objective value of this policy satisfies
\[
\OPT_{\ConesidedA}\geq\sum_{j\in\S}\sum_{C\subseteq C^\star_j} f_j(C)\prod_{i\in C}\phi_i(j,S_i)\prod_{i\notin C}(1-\phi_i(j,S_i))\geq \sum_{j\in\S}\sum_{C\subseteq C^\star_j} f_j(C)\cdot \phi_{\min}^{|C|}\cdot(1-\phi_{\min})^{|\C\setminus C|},
\]
where the second inequality is due to the monotonicity of $f_j$ and  $\phi_i(j_i,S_i))\geq \phi_{\min}$ for all $i\in\C$. To conclude the claim, we use the following lemma:
\begin{lemma}[\citet{feige_etal2011}]\label{lemma:aux_simulation}
Let $g:2^\E\to\RR_+$ be a submodular function over a set of elements $\E$. Denote by $E(p)$ a random subset of $E\subseteq\E$ where each element appears with probability $p$. Then, 
\[
\EE[g(E(p))] = \sum_{A\subseteq E}g(A)\cdot p^{|A|}\cdot (1-p)^{|E\setminus A|}\geq (1-p)g(\emptyset)+pg(E).
\]
\end{lemma}
This lemma implies that 
\[
\OPT_{\ConesidedA}\geq \sum_{j\in\S}\sum_{C\subseteq C^\star_j} f_j(C)\cdot \phi_{\min}^{|C|}\cdot(1-\phi_{\min})^{|\C\setminus C|}\geq \phi_{\min}\cdot \sum_{j\in\S}f_j(C^\star_j) = \phi_{\min}\cdot M^\star.
\]
We conclude the proof of the claim by using  from Theorem~\ref{thm:approxfactor_greedy} that $\mu\geq \frac{1}{2}\cdot\OPT_{\ConesidedA}$. Therefore,
\[
\mu\geq \frac{1}{2}\cdot\OPT_{\ConesidedA}\geq \frac{\phi_{\min}}{2}\cdot M^\star.
\]
\end{proof}

\section{Missing Proofs in Section~\ref{sec:approx_fully_static}}\label{sec:missing_proofs_approx_fullystatic}

\color{black}

\subsection{Proof of Lemma~\ref{thm:low_low}.}\label{apx:low_low_approx}
Consider $\tilde{\bf x}$ the randomized solution defined in \eqref{eq:rando}.
Let $i \in [n]$ and $j \in [m]$. We have
\begin{align*}
    \mathbb{E} \biggr[     \frac{w_{ji}}{1+   \sum_{k=1}^n  w_{jk} \tilde{x}_{kj}    }&\cdot\frac{v_{ij}}{1  +  \sum_{\ell=1}^m  v_{i\ell}  \tilde{x}_{i\ell}   } \cdot  \tilde{x}_{ij}\biggr] \\
    &=     \mathbb{E} \biggr[     \frac{w_{ji}}{1+   \sum_{k=1}^n  w_{jk} \tilde{x}_{kj}    }\frac{v_{ij}}{1  +  \sum_{\ell=1}^m  v_{i\ell}  \tilde{x}_{i\ell}   } \cdot  \tilde{x}_{ij} \biggr| \tilde{x}_{ij} =1    \biggr]  \mathbbm{P} ( \tilde{x}_{ij} =1) \\
    &=     \mathbb{E} \biggr[     \frac{w_{ji}}{1+ w_{ji} +  \sum_{k=1, k\neq i}^n  w_{jk} \tilde{x}_{kj}    }\frac{v_{ij}}{1  +  v_{ij} + \sum_{\ell=1, \ell \neq j}^m  v_{i\ell}  \tilde{x}_{i\ell}   }     \biggr] \cdot y^*_{ij} \\
    &=     \mathbb{E} \biggr[     \frac{w_{ji}}{1+ w_{ji} +  \sum_{k=1, k\neq i}^n  w_{jk} \tilde{x}_{kj}    }  \biggr] \mathbb{E} \biggr[  \frac{v_{ij}}{1  +  v_{ij} + \sum_{\ell=1, \ell \neq j}^m  v_{i\ell}  \tilde{x}_{i\ell}   }   \biggr]   \cdot y^*_{ij} \\
        & \geq        \frac{w_{ji}}{1+ w_{ji} +  \sum_{k=1, k\neq i}^n  w_{jk}   \mathbb{E} [\tilde{x}_{kj }]    }    \frac{v_{ij}}{1  +  v_{ij} + \sum_{\ell=1, \ell \neq j}^m  v_{i\ell}    \mathbb{E} [\tilde{x}_{i\ell} ]  }     \cdot y^*_{ij} \\
                & =        \frac{w_{ji}}{1+ w_{ji} +  \sum_{k=1, k\neq i}^n  w_{jk}  y^*_{kj}   }    \frac{v_{ij}}{1  +  v_{ij} + \sum_{\ell=1, \ell \neq j}^m  v_{i\ell}   {y}^*_{i\ell}   }     \cdot y^*_{ij} \\
                & \geq \frac{ w_{ji}}{1+\alpha+1} \cdot \frac{ v_{ij}}{1+\alpha+1}  y^*_{ij}\\
                & = \frac{1}{(2+\alpha)^2} \cdot w_{ji} v_{ij} y^*_{ij}.
\end{align*}
The first equality holds because $\tilde{x}_{ij}$ is a Bernoulli random variable. The third equality holds because the random variables $\tilde{x}_{kj}$ for $k \neq i$ and $\tilde{x}_{i\ell}$ for $\ell \neq j$ are independent. The first inequality follows from Jensen's inequality (i.e., $\mathbb{E}[1/X] \geq 1 / \mathbb{E}[X])$. The second inequality holds because $w_{ji} \leq \alpha$  and $v_{ij} \leq \alpha$  according to the assumption in this regime. Moreover, we have  $\sum_{\ell=1, \ell \neq j}^m  v_{i\ell}   {y}^*_{i\ell} \leq 1$ and
$\sum_{k=1, k\neq i}^n  w_{jk}  y^*_{kj} \leq 1$ which follow from the first and second constraints in Problem \ref{eq:relax_static_simple}. Finally, by talking the sum over all $i$ and $j$, we get
\begin{align*}
\mathbb{E} \biggr[   \sum_{i=1}^n  \sum_{j=1}^m  \frac{w_{ji}}{1+   \sum_{k=1}^n  w_{jk} \tilde{x}_{kj}    }\frac{v_{ij}}{1  +  \sum_{\ell=1}^m  v_{i\ell}  \tilde{x}_{i\ell}   } \cdot  \tilde{x}_{ij}\biggr] 
                & \geq \frac{1}{(2+\alpha)^2} \cdot   \sum_{i=1}^n  \sum_{j=1}^m w_{ji} v_{ij} y^*_{ij} \\
                &=  \frac{1}{(2+\alpha)^2}\cdot z_{\sf LP} \\
                &\geq \frac{1}{(2+\alpha)^2}\cdot \OPT_{\fullyS},    
\end{align*}
where the last inequality follows from Lemma \ref{lemma:upper_bound_static_lowlow}. This concludes our proof.

\color{black}
\subsection{Proof of Lemma \ref{lemma:upper_bound_static_lowlow}} \label{apx:lemma:upper_bound_static_lowlow}

    Consider an optimal solution ${\bf x}^* \in \{0,1\}^{n \times m}$ of Problem   \eqref{eq:main_static}. We will construct a feasible solution $\mathbf{y}\in\RR^{(n+1)\times(m+1)}$ for the linear program \eqref{eq:relax_static_simple}. For any $i \in [n]$ and $j \in [m]$, let 
    $$ y_{ij}=           \frac{1}{1+   \sum_{k=1}^n  w_{jk} x^*_{kj}    }\frac{1}{1  +  \sum_{\ell=1}^m  v_{i\ell} x^*_{i\ell}   } \cdot x^*_{ij}.         $$
 This solution satisfies the constraints of Problem \eqref{eq:relax_static_simple}. In fact, for any $i\in[n],j \in [m]$, we have
\begin{align*}
     y_{ij}+\sum_{\ell=1}^m v_{i\ell}y_{i\ell} & \leq  \frac{1}{1  +  \sum_{\ell=1}^m  v_{i\ell} x^*_{i\ell}   }+
 \sum_{\ell=1}^m \frac{v_{i\ell} x^*_{i\ell}   }{1  +  \sum_{\ell=1}^m  v_{i\ell} x^*_{i\ell}   }\\
 &=1.
\end{align*}
Similarly, we can show that
\(
y_{ij}+\sum_{k=1}^n w_{jk}y_{kj}\leq 1.
\)
This means that solution $\mathbf{y}$ is feasible for the LP \eqref{eq:relax_static_simple}, thus 
\(z_{\sf LP} \geq   \sum_{i=1}^n  \sum_{j=1}^m v_{ij}  w_{ji} y_{ij} =   \OPT_{\fullyS}.\)

\color{black}
\subsection{Proof of Lemma \ref{thm:high_value_suppliers}}\label{proof:high_value_approx}
Consider $\mathbf{y}\in\{0,1\}^{n\times m}$ a solution to Problem~\eqref{eq:high_value_supp_optimist} with objective value {\sf val}. Assume that suppliers are high value, meaning that $w_{ji}\ge\alpha$. First, note that for every
 $i \in [n]$ and $j \in [m]$, we have
\begin{equation} \label{eq:pales}
     \frac{w_{ji}y_{ij}}{1 + \sum_{k=1}^n w_{j k} y_{kj}} \geq \frac{\alpha}{1+\alpha} \cdot y_{ij}.
\end{equation}
Indeed, if $y_{ij} = 0$, then both sides of the inequality are zero, and if $y_{ij} = 1$, then $y_{kj} = 0$ for all $k \neq i$ because of constraint $\sum_{k=1}^n y_{kj} \leq 1$. Thus,
\[
\frac{w_{ji}y_{ij}}{1 + \sum_{k=1}^n w_{jk} y_{kj}} = \frac{w_{ji}}{1 + w_{ji}} \geq \frac{\alpha}{1+\alpha} = \frac{\alpha}{1+\alpha} \cdot y_{ij},
\]
where the inequality holds because $w_{ji} \geq \alpha$ according to the assumption in this section and the function $\frac{x}{1+x}$ is increasing. This proves the inequality \eqref{eq:pales}.
Then
\begin{align*}
    \sum_{i=1}^n\sum_{j=1}^m \frac{w_{ji}y_{ij}}{1 + \sum_{k=1}^n w_{jk} y_{kj}}\frac{v_{ij}}{1+\sum_{\ell=1}^m v_{i\ell} y_{i\ell}}\ge\frac{\alpha}{1+\alpha}\cdot\sum_{i=1}^n 
    \frac{\sum_{j=1}^m v_{ij}y_{ij}}{1+\sum_{j=1}^m v_{ij}y_{ij}}\ge\frac{\alpha}{1+\alpha}\cdot \text{\sf val}.
\end{align*}

\color{black}

\subsection{Proof of Lemma \ref{lemma:upper_bound_static_concave}} \label{apx:lemma:upper_bound_static_concave}

    Notice that
\begin{equation} \label{eq:disjoint}
z_{\sf C} = \max \left\{\sum_{i\in\C} \; f_i(S_i): \text{ subsets }S_i\subseteq\S 
 \text{ for all $i\in\C$ are disjoint} \right\}
\end{equation}
     is a reformulation of Problem~\eqref{eq:high_value_supp_optimist}. For every $D\subseteq S\subseteq\S$, \revised{let $g_i(D,S)=\sum_{j\in D} \phi_i(j,S)$ be} the relative demand function of customer $i$. Since the MNL choice model satisfies $\phi_i(j,S)\leq \phi_i(j,D)$ for all $D\subseteq S$, then $g_i(D,S)\leq g_i(D,D)=f_i(D)$, \revised{i.e., the relative demand function coincides with the demand function when considering the same sets.}
    Consider an optimal solution of~\fullystatic. We denote as $S^\star_1,\ldots,S_n^\star$ the subsets offered to customers and as $D_1,...,D_n$ the subsets of suppliers that chose each customer. Note that $D_1,\ldots,D_n$ are disjoint. Therefore, we have
    \begin{align*}
        \OPT_{\fullyS}&=\EE\left[\sum_{i=1}^n\one_{\left\{i\text{ is matched}\right\}}\right]=\EE\left[\sum_{i=1}^n g_i(D_i,S_i^\star)\right]\\
        &\le\EE\left[\sum_{i=1}^n g_i(D_i,D_i)\right]=\EE\left[\sum_{i=1}^n f_i(D_i)\right]\le z_{\sf C},
    \end{align*}
    where the last inequality follows from the fact that the maximum of a random variable is bigger than its expected value. 

\color{black}

\subsection{Proof of Lemma \ref{lemma:subadditivity}}\label{proof:subadditivity}
Consider $\mathbf{x}^*$, an optimal solution for Problem \eqref{eq:main_static}. We observe that
\begin{align*}
    \OPT_{\fullyS} &= \sum_{i=1}^n \sum_{j=1}^m \frac{w_{ji}}{1 + \sum_{k=1}^n w_{jk} x^*_{kj}}\cdot\frac{v_{ij}}{1 + \sum_{\ell=1}^m v_{i\ell} x^*_{i\ell}} \cdot x^*_{ij} \\
    &= \sum_{t=1}^3 \sum_{(i,j) \in \E^t} \frac{w_{ji}}{1 + \sum_{k=1}^n w_{jk} x^*_{kj}}\cdot\frac{v_{ij}}{1 + \sum_{\ell=1}^m v_{i\ell} x^*_{i\ell}} \cdot x^*_{ij} \\
    &\leq \sum_{t=1}^3 \sum_{(i,j) \in \E_t} \frac{w_{ji}}{1 + \sum_{k: (k,j) \in \E^t} w_{jk} x^*_{kj}}\cdot\frac{v_{ij}}{1 + \sum_{\ell: (i,\ell) \in \E_t} v_{i\ell} x^*_{i\ell}} \cdot x^*_{ij} \leq \sum_{t=1}^3 \OPT^t_{\fullyS},
\end{align*}
where the first inequality arises because denominators in the fractions are reduced, and the second inequality follows directly from feasibility.

\color{black}

\section{Missing Proofs in Section~\ref{sec:cardinality_constraints}}\label{app:proofs_cardinality_constraints}
For this section, recall that we define the \emph{constrained demand function} of supplier $j\in\S$ as
\begin{equation} \label{eq:cdf}
   f^K_j(C) := \max\left\{\sum_{i\in C'} \phi_{j}(i,C'):\; |C'|\leq K_j, \ C'\subseteq C\right\}, \qquad \forall\; C\subseteq\C. 
\end{equation}
Similarly, we define $f^K_i(\cdot)$ for every customer $i\in\C$. We emphasize that we still assume that the \emph{unconstrained demand function} $f_j^\infty\equiv f_j$ is monotone and submodular. It is easy to see that $f^K_j$ is monotone, however, the submodularity of $f_j$ does not guarantee that $f^K_j$ is submodular as well; see counterexample in Appendix~\ref{app:not_submodular}.

\subsection{Counterexample to the Submodularity of the Constrained Demand Function}\label{app:not_submodular}
To ease the exposition, we do not write the index for the supplier, i.e. $f_j=f$.
Consider a choice model over a set of customers $\{1,2,3,4\}$ determined by a mixture of two MNL choice models with weights $w_{1i}$ and $w_{2i}$ for $i\in \{1,2,3,4\}$ and arrival probabilities $\theta_1=\theta_2= \frac{1}{2}$. Formally, consider coefficients $w_{11}=1,w_{12}=1,w_{13}=0,w_{14}=\infty$ and $w_{21}=0,w_{22}=0,w_{23}=1,w_{24}=0$. Then, the unconstrained demand function $f$ is defined as
\begin{align*}
    f(S)=\frac{1}{2}\left(\frac{\sum_{i\in S}w_{1i}}{1+\sum_{i\in S}w_{1i}}+\frac{\sum_{i\in S}w_{2i}}{1+\sum_{i\in S}w_{2i}}\right), \qquad \forall\; S\subseteq \{1,2,3,4\}.
\end{align*}
Clearly, $f$ is submodular and monotone  increasing. However, for $K =2$, one can easily check that
\begin{align*}
    f^K(\{1,2\})=\frac{1}{3},\quad f^K(\{1,2,4\})=\frac{1}{2},\quad f^K(\{1,2,3\})=\frac{1}{2},\quad f^K(\{1,2,3,4\})=\frac{3}{4}.
\end{align*}
Hence, $f^K(\{1,2,3,4\})-f^K(\{1,2,3\})=1/4>1/6=f^K(\{1,2,4\})-f^K(\{1,2\})$, which means that $f^K$ is not submodular.
\subsection{Proof of Theorem~\ref{thm:cardinality_gap_os_oa}}\label{app:proof_cardinality_gap_os_oa}
In the remainder of this section, we assume that the initiating side is the customers' side; the proofs are analogous for the other case. We divide the proof in three parts:
\begin{enumerate}
    \item[I.] Under the two-way constrained setting and general choice models, the adaptivity gap between \onesidedstatic~and \onesidedadaptive~lies between $e/(e-1)$ and $(e/(e-1))^2$.
    \item[II.] Under the two-way constrained setting, the adaptivity gap between \onesidedstatic~and \onesidedadaptive~is $e/(e-1)$, when customers and suppliers follow MNL choice models.
    \item[III.] Under the one-way constrained setting and general choice models, the adaptivity gap between \onesidedstatic~and \onesidedadaptive~is $e/(e-1)$.
\end{enumerate}

\subsubsection{Two-way Constrained Setting: Analysis under General Choice Models.}
Consider the following linear program, which is a relaxation for \Conesidedadaptive~under the two-way constrained model (when customers are the initiating side):
\begin{subequations}\label{eq:LP_relaxation_gap_os_oa_constrained}
\begin{align}
\max &\quad \sum_{j\in \S}\sum_{C\subseteq \C}f^K_j(C)\cdot \lambda_{j,C} \label{eq:obj_LP_relaxation_gap_os_oa_constrained}\\
s.t. &\quad \sum_{C\subseteq \C}\lambda_{j,C} = 1, \hspace{12em} \text{for all} \ j\in \S \label{eq:distlambda_LP_relaxation_gap_os_oa_constrained}\\
&\quad \sum_{C: C\ni i} \lambda_{j,C} = \sum_{S:|S|\leq K_i, S\ni j}\tau_{i,S}\cdot\phi_i(j,S), \hspace{3em} \text{for all} \ i\in \C, \ j\in \S \label{eq:flowconservation_LP_relaxation_gap_os_oa_constrained}\\
&\quad \sum_{S\subseteq \S:\; |S|\leq K_i} \tau_{i,S} =1, \hspace{10em}\text{for all} \ i\in \C, \label{eq:disttau_LP_relaxation_gap_os_oa_constrained}\\
&\quad \lambda_{j,C}, \ \tau_{i,S}\geq 0, \hspace{12em} \text{for all} \ j\in \S, \ C\subseteq \C, \ i\in \C, \ S\subseteq \S,\; |S|\leq K_i. \notag
\end{align}
\end{subequations}
In particular, we have the following lemma where we use $\OPT_{\eqref{eq:LP_relaxation_gap_os_oa_constrained}}$ to refer to  the objective value of the LP relaxation~\eqref{eq:LP_relaxation_gap_os_oa_constrained}. %
\begin{lemma}\label{lemma:ub_constrained_oa}
$\OPT_{\eqref{eq:LP_relaxation_gap_os_oa_constrained}}\geq \OPT_{\ConesidedA}$.    
\end{lemma}

The proof of Lemma \ref{lemma:ub_constrained_oa} follows the same structure as the proof of Lemma \ref{lemma:relaxation_onesided_customers}, where we show that the LP in \eqref{eq:relaxation_onesided_customers} is a relaxation of \Conesidedadaptive\ in the unconstrained setting.
We do not restate the full argument here, but highlight the main differences: (i) the expected objective now involves the constrained demand functions $f^K_j$ instead of the unconstrained demand functions $f_j$; (ii) the distribution $\tau_{i,S}$ is restricted to subsets $S\subseteq\S$ with $|S|\leq K_i$ for each $i\in\C$.

We remark that we cannot apply the correlation gap argument that we used in Lemma~\ref{lemma:LB_onesided_adaptive} with Relaxation~\eqref{eq:relaxation_onesided_customers} since the functions $f^K_j(\cdot)$'s are not necessarily submodular. We will resolve this challenge by connecting the concept of Contention Resolution (CR) schemes~\citep{chekuri2011submodular}, and the correlation gap.

\noindent{\bf Contention Resolution Schemes.} To ease the explanation, let us consider a fixed supplier, so for now we drop the index $j\in\S$. For a given budget $K\in\ZZ_+$, define $\mathcal{P}:=\{\mathbf{y}\in[0,1]^{|\C|}:\; \sum_{i\in\C}y_i\leq K\}$.
\begin{definition}[\citet{chekuri2011submodular}]\label{def:crs}
An $\alpha$-balanced CR scheme $\nu$ for $\mathcal{P}$ is a procedure that for every $\mathbf{y}\in\mathcal{P}$ and subset $C\subseteq\C$, returns a random set $\nu_{\mathbf{y}}(C)\subseteq C\cap\{i\in\C:\; y_i>0\}$ such that:
\begin{enumerate}
    \item[(i)] $|\nu_{\mathbf{y}}(C)|\leq K$ with probability 1, for all $C\subseteq\C$ and $\mathbf{y}\in\mathcal{P}$.
    \item[(ii)] For all $i\in\C$ with $y_i>0$, we have $\PP(i\in\nu_{\mathbf{y}}(R_{\mathbf{y}})\; | \; i\in R_{\mathbf{y}})\geq \alpha$, where $R_{\mathbf{y}}\subseteq\C$ is the random set such that each $i\in\C$ is included independently with probability $y_i$.
\end{enumerate}
\end{definition}
A CR scheme is said to be monotone if $\PP(i\in\nu_{\mathbf{y}}(C))\geq \PP(i\in\nu_{\mathbf{y}}(C'))$ for all $i\in C\subseteq C'$. With this definition, we focus on the following general result by \citet{chekuri2011submodular} (see their Theorem 1.3): 
\begin{theorem}[\citet{chekuri2011submodular}]\label{thm:crs_bound}
Consider a fractional point $\mathbf{y}\in\mathcal{P}$ and let $\nu$ be a monotone $\alpha$-balanced CR scheme for $\mathcal{P}$. Denote by $\tilde{C} = \nu_{\mathbf{y}}(C(\mathbf{y}))$, where $C(\mathbf{y})$ is the random set in which each $i\in\C$ is drawn independently with probability $y_i$. Then, for any monotone {submodular} set function $f:2^\C\to\RR_+$ we have
\[
\EE[f(\tilde{C})] \geq \alpha\cdot\sum_{C\subseteq\C}f(C)\cdot\prod_{i\in C}y_i\prod_{i\in \C\setminus C}(1-y_i).
\]
\end{theorem}
In particular, for cardinality constraints, \citet{chekuri2011submodular} shows that there exists a monotone $(1-1/e)$-balanced CR scheme (see their Corollary 4.8).

We now come back to our setting and prove the first part of Theorem~\ref{thm:cardinality_gap_os_oa}, i.e., that the adaptivity gap between \onesidedadaptive~and \onesidedstatic~is upper bounded by $(e/(e-1))^2$ and lower bounded by $e/(e-1)$.

\begin{proof}[Proof of Theorem~\ref{thm:cardinality_gap_os_oa}, Part I.]
The lower bound of $e/(e-1)$ is due to Lemma~\ref{lemma:UB_onesided_adaptive} which we proved for the unconstrained setting. Let us now prove the upper bound. 
For this, we consider the one-sided static policy formalized in Algorithm~\ref{alg:policy_constrained_gap_os_oa}.
\begin{algorithm}[htpb]
 \caption{One-sided static policy: Two-way constrained model}\label{alg:policy_constrained_gap_os_oa}
 \begin{algorithmic}[1]
 \color{black}
 \Require  Choice models $\phi$.
 \Ensure Assortments $S_i$ for each $i\in\C$ and $C_j$ for each $j\in\S$.
\State Solve \eqref{eq:LP_relaxation_gap_os_oa_constrained} and denote its optimal solution as $\lambda^\star,\tau^\star$.
\State For each customer $i\in\C$, display a random assortment $S_i$ sampled from distribution $\{\tau^\star_{i,S}\}_{S\subseteq\S}$.
\State Observe customers' choices. Let $R_j$ be the random set of customers that chose $j\in\S$.
\State For each $j\in\S$, solve the maximization problem that computes $f^K_j(R_j)$ and let $C_j$ be a corresponding optimal solution.
\State Display assortment $C_j$ to each supplier $j\in\S$.
 \end{algorithmic}
 \end{algorithm}
We remark that Steps 1 and 4 are not necessarily implementable in polynomial-time, however, for our purposes we only need to show the existence of a feasible one-sided static policy.
\begin{claim}\label{claim:constrained_distribution}
Let $\lambda^\star,\tau^\star$ be an optimal solution of \eqref{eq:LP_relaxation_gap_os_oa_constrained}. Then, we can construct a distribution $\lambda^\dag$ (not necessarily feasible in \eqref{eq:LP_relaxation_gap_os_oa_constrained}) %
that satisfies: 
\begin{itemize}
    \item[(i)] $\lambda_{j,C}^\dag =0$, for all $j\in\C, C\subseteq\C$ with $|C|>K_j$;
    \item[(ii)] $\sum_{C:|C|\leq K_j}\lambda^\dag_{j,C}= 1$, for all $j\in\S$;
    \item[(iii)] $\sum_{C:|C|\leq K_j,\;C\ni i}\lambda^\dag_{j,C} \leq \sum_{S:|S|\leq K_i,\; S\ni j}\tau^\star_{i,S}\cdot\phi_i(j,S)$, for all $i\in\C,j\in\S$;
    \item[(iv)] $\sum_{C\subseteq\C}f_j(C)\cdot \lambda^\dag = \sum_{C\subseteq\C}f_j^K(C)\cdot \lambda^\star$, for all $j\in\S$.
\end{itemize}
\end{claim}
We leave the proof of this claim to the end. Now, consider optimal solutions $\lambda^\star,\tau^\star$ of \eqref{eq:LP_relaxation_gap_os_oa_constrained} and $\lambda^\dag$ satisfying the properties in Claim~\ref{claim:constrained_distribution}. Denote the corresponding marginals as 
\[
x^\star_{ij} = \sum_{S: |S|\leq K_i, S\ni j}\tau_{i,S}^\star\cdot\phi_i(j,S)
\quad \text{and} \quad x^{\dag}_{ij}= \sum_{C:|C|\leq K_j,\;C\ni i} \lambda^\dag_{j,C}.
\]
Because of (iii) in Claim~\ref{claim:constrained_distribution} we know that $x^\dag_{i,j}\leq x^\star_{i,j}$. More importantly, for each $j\in\S$, we have
\begin{align*}
\sum_{i\in\C}x^\dag_{ij} = \sum_{i\in\C}\sum_{C:|C|\leq K_j,\;C\ni i} \lambda^\dag_{j,C} = \sum_{C\subseteq\C: |C|\leq K_j}\lambda^\dag_{j,C}\cdot |C|\leq K_j,
\end{align*}
where we use (i) and (ii) in Claim~\ref{claim:constrained_distribution}. This means that, for each $j\in\S$, we have $\mathbf{x}^\dag_{\cdot j}\in\mathcal{P}_j :=\{\mathbf{y}\in[0,1]^{|\C|}:\; \sum_{i\in\C}y_i\leq K_j\}$. %
Therefore, by using Theorem~\ref{thm:crs_bound} on $\mathbf{x}^\dag$, we get
\begin{equation}\label{eq:aux_constrained_gap}
\sum_{j\in\S}\EE[f_j(\nu_{\mathbf{x}^\dag}(R_j))]\geq \left(1-\frac{1}{e}\right)\cdot\sum_{j\in\S}\sum_{C\subseteq\C}f_j(C)\cdot\prod_{i\in C}x^\dag_{ij}\prod_{i\in\C\setminus C}(1-x^\dag_{ij}),
\end{equation}
where $R_j$ is the random set obtained in Step 4 of Algorithm~\ref{alg:policy_constrained_gap_os_oa}.
We conclude the proof of Theorem~\ref{thm:cardinality_gap_os_oa} as follows
\begin{align*}
\OPT_{\ConesidedS} &\geq \sum_{j\in\S}\sum_{C\subseteq\C}f_j^K(C)\cdot\prod_{i\in C}x^{\star}_{ij}\prod_{i\in\C\setminus C}(1-x^{\star}_{ij})\\
&\geq \sum_{j\in\S}\sum_{C\subseteq\C}f_j^K(C)\cdot\prod_{i\in C}x^{\dag}_{ij}\prod_{i\in\C\setminus C}(1-x^{\dag}_{ij})\\
&\geq \sum_{j\in\S}\sum_{C\subseteq\C}\EE[f_j(\nu_{\mathbf{x}^\dag}(C))\;|\; R_j=C]\cdot\prod_{i\in C}x^{\dag}_{ij}\prod_{i\in\C\setminus C}(1-x^{\dag}_{ij})\\
&=\sum_{j\in\S}\EE[f_j(\nu_{\mathbf{x}^\dag}(R_j))] \\
&\geq \left(1-\frac{1}{e}\right)\cdot\sum_{j\in\S}\sum_{C\subseteq\C}f_j(C)\cdot\prod_{i\in C}x^\dag_{ij}\prod_{i\in\C\setminus C}(1-x^\dag_{ij})\\
&\geq 
\left(1-\frac{1}{e}\right)^2\cdot\sum_{j\in\S}\sum_{C\subseteq\C}f_j(C)\cdot\lambda^\dag_{j,C}\\
&=
\left(1-\frac{1}{e}\right)^2\cdot\sum_{j\in\S}\sum_{C\subseteq\C}f^K_j(C)\cdot\lambda^\star_{j,C}\\
&= \left(1-\frac{1}{e}\right)^2\cdot\OPT_{\eqref{eq:LP_relaxation_gap_os_oa_constrained}}\geq \left(1-\frac{1}{e}\right)^2\cdot\OPT_{\ConesidedA},
\end{align*}
where the right-hand side of the first inequality is the expected value of  Algorithm~\ref{alg:policy_constrained_gap_os_oa} which is a feasible one-sided static policy. The next inequality is due to the monotonicity of $f^K_j$, $x^\dag_{ij}\leq x^\star_{ij}$ for all $i\in\C,j\in\S$ and the fact that multi-linear extensions preserve monotonicity. The following inequality holds because, for each $C\subseteq\C$, the outcome of the CR scheme is such that $\nu_{\mathbf{x}^\dag}(C)\subseteq C$ and $|\nu_{\mathbf{x}^\dag}(C)|\leq K_j$, which means this set is feasible in the maximization problem defined by $f^K_j(C)$. Therefore, $f^K_j(C)\geq f^K_j(\nu_{\mathbf{x}^\dag}(C)) = f_j(\nu_{\mathbf{x}^\dag}(C))$. So the inequality holds conditioned on $R_j=C$, where $R_j$ is the random set of customers that chose $j\in\S$. The fourth inequality is due to Theorem~\ref{thm:crs_bound}. In the fifth inequality, we use the correlation gap for monotone submodular functions~\citep{agrawal2010correlation}, where $\lambda^\dag$ is a feasible solution of 
\begin{equation}\label{eq:aux_LP_constrained_gap}
\max_{\mathbf{\lambda}\geq0}\left\{\sum_{j\in\S}\sum_{C\subseteq\C}f_j(C)\cdot\lambda_{j,C}:\; \sum_{C\subseteq\C}\lambda_{j,C}=1, \ \forall j\in\S, \sum_{C:C\ni i}\lambda_{j,C}=x^{\dag}_{ij}, \ \forall i\in\C,j\in\S \right\}.
\end{equation}
Note that this problem is not the same as \eqref{eq:LP_relaxation_gap_os_oa_constrained} since we replaced $f^K_j$ by $f_j$. Also, it is important to observe that we are using the marginal values $x_{ij}^\dag$, not $x^\star_{ij}$. The next equality follows because of (iv) in Claim~\ref{claim:constrained_distribution}. 
Finally, the last inequality follows from Lemma~\ref{lemma:ub_constrained_oa}.
\end{proof}

\begin{proof}[Proof of Claim~\ref{claim:constrained_distribution}.]
First, notice that for every set $C\subseteq\C$, there exists a subset $C^\star\subseteq C$, which is dependent on $C$, that is an optimal solution to problem $f^K_j(C)$ such that $|C^\star|\leq K_j$. Moreover, we have that $f^K_j(C) = f_j(C^\star)$. Observe that multiple sets $C$ may have the same associated optimal set $C^\star$.

Now, consider an optimal solution $\lambda^\star$ of \eqref{eq:LP_relaxation_gap_os_oa_constrained}. We construct $\lambda^\dag$ as follows: Initiate $\lambda^\dag$ equal to $\lambda^\star$. Then, for each $j\in\S$, $C\subseteq\C$ with $|C|>K_j$, where $C^\star(C)$ is the optimal set for $C$, we do $\lambda^\dag_{j,C^\star(C)}\leftarrow \lambda^\dag_{j,C^\star(C)} + \lambda^\star_{j,C}$ and $\lambda^\dag_{j,C}=0$. We repeat this procedure for every $\lambda^\star_{j,C}>0$ with $|C|>K_j$. Point (i) follows by construction, and points (ii) and (iv) by linearity. Finally, point (iii) also follows by linearity, but notice that we may set some $\lambda^\star$ to 0, which would decrease the left hand side in \eqref{eq:flowconservation_LP_relaxation_gap_os_oa_constrained}.
\end{proof}

\paragraph{Extension to other constraints.} Note that if we have a monotone CR scheme that is $\alpha$-balanced for a given family of constraints, then we can upper bound the adaptivity gap by $\frac{1}{\alpha}\cdot\frac{e}{e-1}$. In particular, \citet{chekuri2011submodular} shows that there exists a monotone $(1-1/e)$-balanced CR scheme for the general class of matroid constraints. This is relevant for matching platforms that require assortments to satisfy partition constraints which aims to diversify the groups of agents that are being displayed. Moreover, \citet{kashaev2023simple} also show that for a cardinality constraint with budget $K$, there exists a monotone CR scheme which is $1-\binom{n}{K}\left(1-\frac{K}{n}\right)^{n-K+1}\left(\frac{K}{n}\right)^K$-balanced, which can be lower bounded by $1-\mathcal{O}(1/\sqrt{K})$. 
In our setting, this means that in a market regime where assortments are large (e.g., $K_i=K_j=\sqrt{n}$), then the upper bound of the adaptivity gap is $\frac{e}{e-1}\cdot\frac{1}{1-\mathcal{O}(1/\sqrt{K})}$ which tends to $e/(e-1)$. 

\subsubsection{Two-way Constrained Setting: Analysis for MNL Choice Models.}\label{app:twowayconstrained_mnl}
In this section, we assume that customers' and suppliers' preferences are dictated by MNL choice models with weights $v_{ij}$ and $w_{ji}$ for all $i\in\C,j\in\S$, respectively. First, we show the following:
\begin{lemma}\label{lemma:MNL_constrained_submodularity}
Under the MNL choice model, the constrained demand function defined in \eqref{eq:cdf} is submodular.
\end{lemma}
\begin{proof}[Proof.]
Consider any supplier $j\in\S$ (a similar argument applies for customers). Then, for every $C\subseteq\C$, we have
\begin{align*}
    f^K_j(C)= \max\left\{\frac{w_j(C')}{1+w_j(C')}:\; |C'|\leq K_j, \ C'\subseteq C\right\}, \qquad \forall\; C\subseteq\C,
\end{align*}
where $w_j(C')=\sum_{i\in C'}w_{ji}$. Since the function $x/(1+x)$ is increasing, we note that the optimal assortments for the maximization problem defined in $f_j^K$ are also optimal for $h_j(C):=\max\{w_j(C'):\; |C'|\leq K_j,\; C'\subseteq C\}$, and vice versa. Therefore, we obtain 
\[
f^K(C) = \frac{h_j(C)}{1+h_j(C)}.
\]
Since $x/(1+x)$ is concave and $h_j(\cdot)$ is submodular~\citep{nemhauser1978analysis}, the result follows.
\end{proof}
We now prove that, under the two-way constrained model, the gap between \onesidedadaptive~and \onesidedstatic~ is $e/(e-1)$ when customers and suppliers follow MNL choice models.
\begin{proof}[Proof of Theorem~\ref{thm:cardinality_gap_os_oa}, Part II.]
Lemma~\ref{lemma:UB_onesided_adaptive} gives us a lower bound for the adaptivity gap. The upper bound of $e/(e-1)$ can be shown as in Lemma~\ref{lemma:LB_onesided_adaptive} by using Relaxation~\eqref{eq:LP_relaxation_gap_os_oa_constrained} and the submodularity of $f^K_j(\cdot)$ proved in Lemma~\ref{lemma:MNL_constrained_submodularity} (recall that this function is always monotone, independently of the choice model).
\end{proof}

\subsubsection{One-way Constrained Setting.} In this section, we assume that the constraints of the agents on the responding side (in this case, suppliers) are such that $K_j=+\infty$ for all $j\in\S$.
\begin{proof}[Proof of Theorem~\ref{thm:cardinality_gap_os_oa}, Part III.]
In this part, we show that gap is exactly $e/(e-1)$ under the one-way constrained model. As before, Lemma~\ref{lemma:UB_onesided_adaptive} gives us a lower bound. To prove the $e/(e-1)$ upper bound, we consider Relaxation~\eqref{eq:LP_relaxation_gap_os_oa_constrained}. Since $K_j=+\infty$, then $f^K_j \equiv f_j$ which we assumed is monotone and submodular. Therefore, the proof is analogous to Lemma~\ref{lemma:LB_onesided_adaptive}.

\end{proof}

\subsection{Proof of Theorem~\ref{thm:approx_oa_constrained_general}}\label{app:proof_cardinality_approximation_oa}
We prove this result for \Conesidedadaptive~and an analogous argument applies to \Sonesidedadaptive. 
Our main algorithm for Theorem~\ref{thm:approx_oa_constrained_general} is formalized in Algorithm~\ref{alg:staticpolicy_constrained_oa}. This method is analogous to Algorithm~\ref{alg:policy_constrained_gap_os_oa}, however, we need to approximately solve~\eqref{eq:LP_relaxation_gap_os_oa_constrained} (Step 1). Indeed, note that \eqref{eq:LP_relaxation_gap_os_oa_constrained} has exponentially many variables, but polynomial number of constraints. As we show next, we can obtain an approximate solution in polynomial time by using the ellipsoid method with an approximate separation oracle to solve its dual.
\begin{claim}\label{claim:approximate_solution_ellipsoid}
There exists a polynomial-time algorithm that returns a feasible solution $(\lambda,\tau)$ of \eqref{eq:LP_relaxation_gap_os_oa_constrained} such that its objective value is at least $(1-1/e)\cdot\OPT_{\eqref{eq:LP_relaxation_gap_os_oa_constrained}}$.
\end{claim}
We leave the proof of this claim to the end of this section.
We remark  that Step 4 in our algorithm is also implementable in polynomial time due to the oracle in Assumption~\ref{assumption:constrained_oa_oneway}. Later, we  discuss how this result extends when using an approximate oracle.

\begin{algorithm}[htpb]
\color{black}
 \caption{One-sided static policy for \Conesidedadaptive: Two-way constrained model}\label{alg:staticpolicy_constrained_oa}
 \begin{algorithmic}[1]
 \Require  Choice models $\phi$.
 \Ensure Assortments $S_i$ for each $i\in\C$ and $C_j$ for each $j\in\S$.
\State Obtain an approximate solution $(\lambda,\tau)$ of \eqref{eq:LP_relaxation_gap_os_oa_constrained} using Claim~\ref{claim:approximate_solution_ellipsoid}.
\State For each customer $i\in\C$, display a random assortment $S_i$ sampled from distribution $\{\tau_{i,S}\}_{S\subseteq\S}$.
\State Observe customers' choices. Let $R_j$ be the random set of customers that chose $j\in\S$.
\State For each $j\in\S$, let $C_j=\argmax\left\{\sum_{i\in C}\phi_j(i,C):\; C\subseteq R_j,\; |C|\leq K_j\right\}$.
\State Display assortment $C_j$ to each supplier $j\in\S$.
 \end{algorithmic}
 \end{algorithm}
\begin{proof}[Proof of Theorem~\ref{thm:approx_oa_constrained_general}.]
Let $(\lambda,\tau)$ be the approximate solution found in Step 1.
As in the proof of Theorem~\ref{thm:cardinality_gap_os_oa}, we can lower bound the expected objective value achieved by Algorithm~\ref{alg:staticpolicy_constrained_oa} as follows
\[
\sum_{j\in\S}\EE[f_j(C_j)] = \sum_{j\in\S}\sum_{C\subseteq\C}f_j^K(C)\cdot\prod_{i\in C}x_{ij}\prod_{i\in\C\setminus C}(1-x_{ij})\geq \left(1-\frac{1}{e}\right)^2\cdot\sum_{j\in\S}\sum_{C\subseteq\C}f_j(C)\cdot\lambda_{j,C},
\]
where $x_{ij}=\sum_{S:|S|\leq K_i,S\ni j}\tau_{i,S}\cdot\phi_i(j,S)$. We finish the proof by using Claim~\ref{claim:approximate_solution_ellipsoid} to show that
\[
\sum_{j\in\S}\sum_{C\subseteq\C}f_j(C)\cdot\lambda_{j,C}\geq \left(1-\frac{1}{e}\right)\cdot\OPT_{\eqref{eq:LP_relaxation_gap_os_oa_constrained}}\geq \left(1-\frac{1}{e}\right)\cdot\OPT_{\ConesidedA}.
\]
\end{proof}

\subsubsection{Missing Proofs.} 
Our goal in this section is to show Claim~\ref{claim:approximate_solution_ellipsoid}. First, consider the dual of \eqref{eq:LP_relaxation_gap_os_oa_constrained}:
\begin{subequations}\label{eq:dual_relaxation_onesided_customers_constrained}
\begin{align}
 \min &\quad \sum_{j\in \S}\beta_j + \sum_{i\in \C} \alpha_i \label{eq:obj_dual_relaxation_onesided_customers_constrained}\\
 s.t. &\quad  \beta_j +\sum_{i\in C}\gamma_{i,j} \geq f^K_j(C), \hspace{1.4em} \qquad \text{for all} \ C\subseteq\C, \ j\in \S, \label{eq:demand_dual_relaxation_onesided_customers_constrained}\\
 &\quad \alpha_i - \sum_{j\in S}\gamma_{i,j}\cdot \phi_{i}(j,S) \geq 0, \qquad \text{for all} \ S\subseteq\S, \ |S|\leq K_i,\ i\in \C,\label{eq:choicemodel_dual_relaxation_onesided_customers_constrained}\\
 &\quad \alpha_i, \ \beta_j, \ \gamma_{i,j}\in\RR, \qquad \hspace{5em}\text{for all } i\in \C, \ j\in \S.\notag
 \end{align}
\end{subequations}
We first observe the following:
\begin{claim}\label{claim:equivalent_duals}
Assume customers and suppliers follow weak-substitutable choice models. Then, problem \eqref{eq:dual_relaxation_onesided_customers_constrained} is equivalent to:
\begin{subequations}\label{eq:dual_relaxation_onesided_customers_constrained2}
\begin{align}
 \min &\quad \sum_{j\in \S}\beta_j + \sum_{i\in \C} \alpha_i \label{eq:obj_dual_relaxation_onesided_customers_constrained2}\\
 s.t. &\quad  \beta_j +\sum_{i\in C}\gamma_{ij} \geq f_j(C), \hspace{1.4em} \qquad \text{for all} \ C\subseteq\C,\; |C|\leq K_j, \ j\in \S, \label{eq:demand_dual_relaxation_onesided_customers_constrained2}\\
 &\quad \alpha_i - \sum_{j\in S}\gamma_{ij}\cdot \phi_{i}(j,S) \geq 0, \qquad \text{for all} \ S\subseteq\S, \ |S|\leq K_i,\ i\in \C,\label{eq:choicemodel_dual_relaxation_onesided_customers_constrained2}\\
 &\quad \alpha_i, \ \beta_j, \ \gamma_{i,j}\geq 0, \qquad \hspace{5em}\text{for all } i\in \C, \ j\in \S.\notag
 \end{align}
\end{subequations}
\end{claim}
The main differences between \eqref{eq:dual_relaxation_onesided_customers_constrained} and \eqref{eq:dual_relaxation_onesided_customers_constrained2} are: (i) the right-hand side of \eqref{eq:demand_dual_relaxation_onesided_customers_constrained2} is $f_j(C)$ instead of $f_j^K(C)$, and the constraint is only for subsets $C$ such that $|C|\leq K_j$; (ii) variables are non-negative. We leave the proof of this claim to the end. Now, let us complete the proof of Claim \ref{claim:approximate_solution_ellipsoid}.

\begin{proof}[Proof of Claim~\ref{claim:approximate_solution_ellipsoid}.]
To prove this result, we will make use of Problem~\eqref{eq:dual_relaxation_onesided_customers_constrained2}. We have two main observations regarding the separation of the constraints in~\eqref{eq:dual_relaxation_onesided_customers_constrained2}: Given a point $(\alpha,\beta,\gamma)$,
\begin{itemize}
    \item First, we can separate Constraints~\eqref{eq:choicemodel_dual_relaxation_onesided_customers_constrained2} in polynomial time as follows: Fix $i\in\C$ and find a subset $A\in\argmax\{\sum_{j\in S}\gamma_{ij}\cdot \phi_i(j,S):\; A\subseteq\S,\; |A|\leq K_i\}$, which can be done efficiently thanks to Assumption~\ref{assumption:constrained_oa_oneway}. Here, it is important to remark that the previous argument is possible because $\gamma_{ij}\geq0$ for all $i\in\C,j\in\S$. Then, if $\alpha_i< \sum_{j\in A}\gamma_{ij}\cdot\phi_i(j,A)$, we include the corresponding constraint; otherwise, there are no constraints to be added;
    \item Second, we can approximately separate Constraints~\eqref{eq:demand_dual_relaxation_onesided_customers_constrained2} in polynomial time as follows: Fix $j\in\S$ and, using the distorted greedy algorithm in \citep{harshaw2019submodular}, find a set $T\subseteq\C$ with $|T|\leq K_j$ such that 
\[
f_j(T)-\sum_{i\in T}\gamma_{ij}\geq \left(1-\frac{1}{e}\right)\cdot f_j(C) -\sum_{i\in C}\gamma_{ij},
\]
for all  $C\subseteq\C,\; |C|\leq K_j$. 
Then, if $\beta_i < f_j(T)-\sum_{i\in T}\gamma_{ij}$, we separate the corresponding constraint; otherwise, do nothing. As before, this argument is possible because $\gamma_{ij}\geq0$ for all $i\in\C,j\in\S$, which is necessary for the algorithm in \citep{harshaw2019submodular}.
\end{itemize}
Therefore, by using the ellipsoid method with both separation methods above as oracles (see e.g., Chapter 2 in~\citep{bubeck2015convex} and \citep{jansen2003approximate}), we can obtain a point $(\alpha,\beta,\gamma)$ such that 
\begin{align}
&\beta_j +\sum_{i\in C}\gamma_{ij}\geq (1-1/e)\cdot f_j(C) \quad\quad \text{for all} \ j\in\S,\; C\subseteq\C,\; |C|\leq K_j, \label{eq:properties_ellipsoid1}\\
&\alpha_i - \sum_{j\in S}\gamma_{ij}\cdot\phi_i(j,S)\geq 0 \hspace{4em} \text{for all} \ i\in\C,\; S\subseteq\S,\; |S|\leq K_i.     \label{eq:properties_ellipsoid2}
\end{align}
Given the above, let {\sf (Modified Primal)} be Problem~\eqref{eq:LP_relaxation_gap_os_oa_constrained} with the following modified constraints:
\begin{align}
\sum_{C\subseteq\C}\lambda_{j,C} &= \frac{e}{e-1}, \quad \text{for all} \ j\in\S,\label{eq:modified_primal1}\\
\sum_{S\subseteq\S}\tau_{i,S} &= \frac{e}{e-1}, \quad \text{for all} \ i\in\C.\label{eq:modified_primal2}
\end{align}
We denote the dual of {\sf (Modified Primal)} as {\sf (Modified Dual)} which corresponds to
\begin{align*}
 \min &\quad \frac{e}{e-1}\cdot\sum_{j\in \S}\beta_j + \frac{e}{e-1}\cdot\sum_{i\in \C} \alpha_i \notag\\
 s.t. &\quad  \beta_j +\sum_{i\in C}\gamma_{ij} \geq f_j(C), \hspace{1.4em} \qquad \text{for all} \ C\subseteq\C,\; |C|\leq K_j, \ j\in \S, \notag\\
 &\quad \alpha_i - \sum_{j\in S}\gamma_{ij}\cdot \phi_{i}(j,S) \geq 0, \qquad \text{for all} \ S\subseteq\S, \ |S|\leq K_i,\ i\in \C,\notag\\
 &\quad \alpha_i, \ \beta_j, \ \gamma_{i,j}\geq 0, \qquad \hspace{5em}\text{for all } i\in \C, \ j\in \S.\notag
\end{align*}
Let $(\alpha^0,\beta^0,\gamma^0)$ be the solution for {\sf (Modified Dual)} found by the ellipsoid algorithm with the separation procedure explained before and $\hat{\OPT}\text{\sf (Modified Dual)}$ be its objective value. Now, consider $(\alpha^1,\beta^1,\gamma^1) = \frac{e}{e-1}\cdot(\alpha^0,\beta^0,\gamma^0)$. We note that:
\begin{itemize}
    \item For all $j\in\S$, $C\subseteq\C, \ |C|\leq K_j$,
    \[
    \beta^1_j+\sum_{i\in C}\gamma_{ij}^1 = \frac{e}{e-1}\left(\beta^0_j+\sum_{i\in C}\gamma_{ij}^0\right) \geq f_j(C),
    \]
    where the inequality follows because $(\alpha^0,\beta^0,\gamma^0)$ satisfies \eqref{eq:properties_ellipsoid1}.
    \item For all $i\in\C$, $S\subseteq\S, \ |S|\leq K_i$,
    \[
    \alpha^1_i-\sum_{j\in S}\gamma_{ij}^1\cdot\phi_i(j,S) = \frac{e}{e-1}\left(\alpha^0_i-\sum_{j\in S}\gamma_{ij}^0\cdot\phi_i(j,S)\right) \geq 0,
    \]
    where the inequality follows because $(\alpha^0,\beta^0,\gamma^0)$ satisfies \eqref{eq:properties_ellipsoid2}.
\end{itemize}
This means that $(\alpha^1,\beta^1,\gamma^1)$ is feasible in \eqref{eq:dual_relaxation_onesided_customers_constrained}. Therefore, we have that
\begin{equation}\label{eq:aux_proof_ellipsoid}
\hat{\OPT}\text{\sf (Modified Dual)} = \frac{e}{e-1}\left(\sum_{i\in\C}\alpha_i^0+\sum_{j\in\S}\beta_j^0\right) =\sum_{i\in\C}\alpha_i^1+\sum_{j\in\S}\beta_j^1\geq \OPT_{\eqref{eq:dual_relaxation_onesided_customers_constrained}} = \OPT_{\eqref{eq:LP_relaxation_gap_os_oa_constrained}}, 
\end{equation}
where the last equality is due to strong duality.

Let us come back for a moment to the ellipsoid procedure. Observe that, at each step, this method finds some constraints (possibly for each $i\in\C$ and $j\in\S$) which are added to the linear program. Let $\mathcal{F}$ be the final set of constraints added to the program. Note that $\mathcal{F}$ is composed by a polynomial number of constraints. Let {\sf (Restricted Dual)} be {\sf (Modified Dual)} restricted only to those constraints in $\mathcal{F}$ and {\sf (Restricted Primal)} be its primal. Our goal now is to compare the optimal value of {\sf (Restricted Dual)} and the objective value of $(\alpha^0,\beta^0,\gamma^0)$ denoted as $\hat{\OPT}\text{\sf (Modified Dual)}$.
Note that $\hat{\OPT}\text{\sf (Modified Dual)}$ 
is the optimal value of a linear program similar to {\sf (Restricted Dual)}
where we only restrict variables to the constraints in $\mathcal{F}$ and the constraints of the form $\beta_j+\sum_{i\in C}\gamma_{ij}\geq f_j(C)$ are changed to $\beta_j+\sum_{i\in C}\gamma_{ij}\geq (1-1/e)\cdot f_j(C)$. This implies that the optimal solution of {\sf (Restricted Dual)} is feasible in this linear program, which means that 
\begin{equation}\label{eq:aux_proof_ellipsoid2}
\OPT\text{\sf (Restricted Primal)}=\OPT\text{\sf (Restricted Dual)}\geq \hat{\OPT}\text{\sf (Modified Dual)},
\end{equation}
where the equality is due to strong duality. Thanks to \eqref{eq:aux_proof_ellipsoid} and \eqref{eq:aux_proof_ellipsoid2}, we have that
\[
\OPT\text{\sf (Restricted Primal)}\geq \OPT_{\eqref{eq:LP_relaxation_gap_os_oa_constrained}}.
\]
Now, we solve {\sf (Restricted Primal)} to optimality, which we can do since it has a polynomial number of constraints and variables (the latter is because the constraints of {\sf (Restricted Dual)} are in $\mathcal{F}$). Let $(\lambda^0,\tau^0)$ be the optimal solution found, which satisfies \eqref{eq:modified_primal1} and \eqref{eq:modified_primal2}. Define $(\lambda^1,\tau^1)=\frac{e-1}{e}\cdot(\lambda^0,\tau^0)$ which is clearly feasible in \eqref{eq:LP_relaxation_gap_os_oa_constrained}. Therefore,
\begin{align*}
\sum_{j\in\S}\sum_{C\subseteq\C}f_j(C)\cdot\lambda^1_{j,C}&\geq \left(1-\frac{1}{e}\right)\cdot \sum_{j\in\S}\sum_{C\subseteq\C}f_j(C)\cdot\lambda^0_{j,C} \\&= \left(1-\frac{1}{e}\right)\cdot\OPT\text{\sf (Restricted Primal)}\\
&\geq \left(1-\frac{1}{e}\right)\cdot\OPT_{\eqref{eq:LP_relaxation_gap_os_oa_constrained}},
\end{align*}
which concludes the proof of the claim.
\end{proof}

\begin{proof}[Proof of Claim~\ref{claim:equivalent_duals}.]
Consider an optimal solution  $(\alpha, \beta, \gamma)$ for \eqref{eq:dual_relaxation_onesided_customers_constrained2}. It clearly satisfies the constraint \eqref{eq:choicemodel_dual_relaxation_onesided_customers_constrained}. It also satisfies the constraint \eqref{eq:demand_dual_relaxation_onesided_customers_constrained} because for $C \in \C$ and $j\in \S$, let $C^* \subseteq C$ such that $|C^*| \leq K_j$ and $f^K_j(C)=f_j(C^*)$. Then we get, $\beta_j + \sum_{j \in C} \gamma_{ij} \geq \beta_j + \sum_{j \in C^*} \gamma_{ij} \geq f_j(C^*) =f_j^K(C)$, where the first inequality uses the non-negativity of $\gamma_{ij}$ and the second one uses the feasibility of the solution for \eqref{eq:demand_dual_relaxation_onesided_customers_constrained2}. Therefore, we conclude that  $(\alpha, \beta, \gamma)$ is also feasible in \eqref{eq:dual_relaxation_onesided_customers_constrained}, which implies that $\OPT_{\eqref{eq:dual_relaxation_onesided_customers_constrained}}\leq \OPT_{\eqref{eq:dual_relaxation_onesided_customers_constrained2}}$. %

Our goal now is to show the opposite direction, i.e., $\OPT_{\eqref{eq:dual_relaxation_onesided_customers_constrained2}}\leq \OPT_{\eqref{eq:dual_relaxation_onesided_customers_constrained}}$. Consider an optimal solution $(\alpha,\beta,\gamma)$ of \eqref{eq:dual_relaxation_onesided_customers_constrained}. First, by taking $C=\emptyset$ in \eqref{eq:demand_dual_relaxation_onesided_customers_constrained} and $S=\emptyset$ in ~\eqref{eq:choicemodel_dual_relaxation_onesided_customers_constrained}, we obtain that $\alpha_i\geq0$ and $\beta_j\geq0$ for all $i\in\C$, $j\in\S$. Now, we construct the following solution $(\alpha,\beta,\gamma')$ with $\gamma'=\max\{\gamma,0\}$, where the max operator is component-wise, and show that is feasible for \eqref{eq:dual_relaxation_onesided_customers_constrained2}.

Since $\alpha$ and $\gamma$ satisfy \eqref{eq:choicemodel_dual_relaxation_onesided_customers_constrained}, then for all $i\in\C$ we have
\[
\alpha_i\geq \sum_{j\in S^\star_i}\gamma_{ij}\cdot\phi_i(j,S^\star_i),
\]
where $S^\star_i\in\argmax\left\{\sum_{j\in S}\gamma_{ij}\cdot\phi_i(j,S):\; S\subseteq\S,\; |S|\leq K_i\right\}$. We now prove that for all $j\in S^\star_i$ we have $\gamma_{ij}\geq 0$. Assume by contradiction that $N_i=\{j\in S^\star_i:\; \gamma_{ij}<0\}$ is non-empty. Then, note
\[
\sum_{j\in S^\star_i\setminus N_i}\gamma_{ij}\cdot\phi_i(j,S^\star_i\setminus N_i)\geq \sum_{j\in S^\star_i\setminus N_i}\gamma_{ij}\cdot\phi_i(j,S^\star_i)> \sum_{j\in S^\star_i}\gamma_{ij}\cdot\phi_i(j,S^\star_i),
\]
where the first inequality is due to weak-substitutability. This contradicts the optimality of $S^\star_i$. Therefore
\begin{align*}
\alpha_i&\geq\max\left\{\sum_{j\in S}\gamma_{ij}\cdot\phi_i(j,S):\; S\subseteq\{j\in\S:\gamma_{ij}\geq0\},\; |S|\leq K_i\right\} \\
&= \max\left\{\sum_{j\in S}\gamma'_{ij}\cdot\phi_i(j,S):\; S\subseteq\{j\in\S:\gamma_{ij}\geq0\},\; |S|\leq K_i\right\} \\
& \geq \max\left\{\sum_{j\in S}\gamma'_{ij}\cdot\phi_i(j,S):\; S\subseteq\S,\; |S|\leq K_i\right\},
\end{align*}
where the last inequality is due to weak-substitutability. Furthermore, for $C \subseteq \C$ and $j \in \S$ such that $|C|\leq K_j$, we have
$\beta_j + \sum_{i \in C} \gamma'_{ij} \geq \beta_j + \sum_{i \in C} \gamma_{ij} \geq f_j^K(C) \geq f_j(C) $,
where the first inequality holds because $\gamma'_{ij} \geq \gamma_{ij}$, the second one follows from \eqref{eq:demand_dual_relaxation_onesided_customers_constrained}, and the last one from the maximality of $f_j^K(C)$.

Finally, note that the objective values of $(\alpha,\beta,\gamma')$ and $(\alpha,\beta,\gamma)$ are the same, therefore 
\(
\OPT_{\eqref{eq:dual_relaxation_onesided_customers_constrained}} = \sum_{i\in\C}\alpha_i+\sum_{j\in\S}\beta_j \geq \OPT_{\eqref{eq:dual_relaxation_onesided_customers_constrained2}}
\). 
\end{proof}

\subsection{Proof of Corollary~\ref{coro:improved_guarantees_constrained}}\label{app:proof_corollary_improved_approx_constrained}
In the following, we focus on \Conesidedadaptive. The algorithm for \Sonesidedadaptive~is analogous.
\subsubsection{Algorithm for MNL choice models under the two-way constrained setting}\label{app:proof_cardinality_apx_oa_twoway}
For this setting, our method is analogous to Algorithm~\ref{alg:random_order_greedy} with two modifications: (i) in Step 3 the marginal values are defined as $f^K_j(C_j^{t-1}\cup\{i_t\})-f^K_j(C_j^{t-1})$; (ii) in Step 4, a cardinality constraint $|S|\leq K_i$ must be imposed when optimizing over assortments $S\subseteq\S$. For (i), we emphasize that these marginal values can be computed in polynomial time, since for any $C\subseteq\C$, an optimal solution in $f^K_j(C)$ can be constructed with the top-$K_j$ values among $w_{ji}$ with $i\in C$ (this follows from the proof of Lemma~\ref{lemma:MNL_constrained_submodularity} in Appendix~\ref{app:twowayconstrained_mnl}). For (ii), we note that, for any $i\in\C$, the maximization problem in Step 4 can be solved in polynomial-time thanks to the algorithm proposed by \citep{rusmevich_etal10} for MNL choice models. 

Given the algorithm above, we can prove the $1/2$ approximation factor in the same way we proved Theorem~\ref{thm:approxfactor_greedy}. For this, we use the dual formulation~\eqref{eq:dual_relaxation_onesided_customers_constrained}
 along with the submodularity of $f^K_j$ for the MNL choice model (Lemma~\ref{lemma:MNL_constrained_submodularity}). 

\subsubsection{One-way constrained setting.}\label{app:oneway_constrained_apx_oa} For this model, the dual formulation is similar to~\eqref{eq:dual_relaxation_onesided_customers_constrained} with the difference that in the right-hand side of the first constraint we write $f_j(C)$ instead of $f^K_j(C)$.

\noindent{\bf Algorithm.} For the one-way constrained setting, our method is analogous to Algorithm~\ref{alg:random_order_greedy} with only one modification: In Step 4, a cardinality constraint $|S|\leq K_i$ must be imposed when optimizing over assortments $S\subseteq\S$.

Given the algorithm above, we can prove the $1/2$ approximation factor in the same way we proved Theorem~\ref{thm:approxfactor_greedy}, thanks to the properties of the unconstrained demand function (Assumption~\ref{assumption:constrained_oa_oneway}).

\subsection{Proof of Theorem~\ref{thm:fullystatic_mnl_approx_cardinality}}\label{app:proof_cardinality_approximation_fs}

To prove this result, we follow the same approach as in the proof of Theorem~\ref{thm:fullystatic_mnl_approx}, dividing the analysis into different regimes: low-value customers and suppliers (Appendix~\ref{app:cardinality_fs_low_low}) and high-value customers (Appendix~\ref{app:cardinality_fs_high_suppliers}). For each regime, we present the corresponding relaxations and describe the associated randomized rounding procedures. With these components in place, the proof of Theorem~\ref{thm:fullystatic_mnl_approx_cardinality} is analogous to that of Theorem~\ref{thm:fullystatic_mnl_approx}. So it is sufficient to show the guarantees for each regime as we will describe next.

\subsubsection{Low value Customers and Low value Suppliers.}\label{app:cardinality_fs_low_low}
We start by presenting the main approximation for this regime.

\begin{lemma}\label{lemma:upper_bound_lowlow_cardinality}
    Suppose $w_{ji}\leq \alpha$ and $v_{ij}\leq \alpha$ for all $(i,j)\in \C \times \S $ for some $\alpha>0$. There exists a randomized polynomial time algorithm that gives a  $\frac{1}{(2+\alpha)^2}$-approximation to \fullystatic~under the two-way constrained model.
\end{lemma}

Under the two-way constrained setting, an upper bound to the \fullystatic~problem is given by the following linear relaxation
\begin{align} \label{eq:relax_static_simple_K}
z_{\sf LP}= \max &\quad \sum_{i=1}^n  \sum_{j=1}^m  v_{ij}w_{ji} y_{ij} &\\
&\quad y_{ij}+\sum_{j\ell=1}^m v_{i\ell}y_{i\ell} \leq 1, \qquad & \forall j \in [m], \notag \\
&\quad y_{ij}+\sum_{k=1}^n w_{jk}y_{kj} \leq 1, \qquad & \forall i \in [n], \notag \\
&\quad\sum_{j=1}^m y_{ij}\le K_i, \qquad & \forall i \in [n], \notag \\
&\quad\sum_{i=1}^n y_{ij}\le K_j, \qquad & \forall j \in [m], \notag\\
&\quad y_{ij}\ge0, \qquad & \forall i\in[n], \forall j \in [m]. \notag
\end{align}
Analogously to Lemma~\ref{lemma:upper_bound_static_lowlow} (proof in Appendix~\ref{apx:lemma:upper_bound_static_lowlow}), we can show that this relaxation is an upper bound. The only difference is that it incorporates cardinality constraints.\\
\noindent{\bf Algorithm.} Given a feasible solution $\mathbf{y}$ to Problem~\eqref{eq:relax_static_simple_K}, we use the polynomial-time dependent rounding scheme introduced in \citep{gandhi_etal06} to get a solution $\mathbf{x}\in\{0,1\}^{n\times m}$. In particular, this solution verifies for every $i\in[n],j\in[m]$:
\begin{itemize}
    \item $\EE[x_{ij}]=y_{ij}$;
    \item $\sum_{i\in[n]} x_{ij}\le K_j$ and $\sum_{j\in[m]} x_{ij}\le K_i$;
    \item {\bf Negative correlation.} $\PP(x_{i\ell}=b, x_{ij}=b)\leq\PP(x_{i\ell}=b)\cdot\PP(x_{ij}=b)$ for all $b\in\{0,1\}$, $\ell\in[m],\ell\neq j$. Similarly, $\PP(x_{kj}=b, x_{ij}=b)\leq\PP(x_{kj}=b)\cdot\PP(x_{ij}=b)$ for all $b\in\{0,1\}$, $k\in[n],k\neq i$.
\end{itemize}

\begin{proof}[Proof of Lemma~\ref{lemma:upper_bound_lowlow_cardinality}.]
Let $\mathbf{x}$ be the randomized solution obtained from the dependent rounding scheme applied to $\mathbf{y}$, a solution of Problem~\eqref{eq:relax_static_simple_K}.

Fix $i \in [n]$ and $j \in [m]$. Then, we obtain:
\begin{align*}
    \mathbb{E} \biggr[     \frac{w_{ji}}{1+   \sum_{k=1}^n  w_{jk} x_{kj}    }&\cdot\frac{v_{ij}}{1  +  \sum_{\ell=1}^m  v_{i\ell}  x_{i\ell}} \cdot  x_{ij}\biggr] \\
    &=     \mathbb{E} \biggr[     \frac{w_{ji}}{1+   \sum_{k=1}^n  w_{jk} x_{kj}    }\cdot\frac{v_{ij}}{1  +  \sum_{\ell=1}^m  v_{i\ell}  x_{i\ell}}\cdot  x_{ij} \biggr| x_{ij} =1    \biggr]  \mathbbm{P} ( x_{ij} =1) \\
    &=     \mathbb{E} \biggr[     \frac{w_{ji}}{1+ w_{ji} +  \sum_{k=1, k\neq i}^n  w_{jk} x_{kj}    }\cdot\frac{v_{ij}}{1  +  v_{ij} + \sum_{\ell=1, \ell \neq j}^m  v_{i\ell}  x_{i\ell}   }\biggr| x_{ij} =1     \biggr] \cdot y_{ij} \\
        & \geq        \frac{w_{ji}}{1+ w_{ji} +  \sum_{k=1, k\neq i}^n  w_{jk}   \mathbb{E} [x_{kj }\mid x_{ij}=1]    }  \cdot  \frac{v_{ij}}{1  +  v_{ij} + \sum_{\ell=1, \ell \neq j}^m  v_{i\ell}    \mathbb{E} [x_{i\ell}\mid x_{ij}=1]  }     \cdot y_{ij} \\
        & \geq        \frac{w_{ji}}{1+ w_{ji} +  \sum_{k=1, k\neq i}^n  w_{jk}   \mathbb{E} [x_{kj}]    } \cdot   \frac{v_{ij}}{1  +  v_{ij} + \sum_{\ell=1, \ell \neq j}^m  v_{i\ell}    \mathbb{E} [x_{i\ell}]  }     \cdot y_{ij} \\
                & = \frac{w_{ji}}{1+ w_{ji} +  \sum_{k=1, k\neq i}^n  w_{jk}  y_{kj}   }    \frac{v_{ij}}{1  +  v_{ij} + \sum_{\ell=1, \ell \neq j}^m  v_{i\ell}   {y}_{i\ell}   }     \cdot y_{ij} \\
                & \geq \frac{ w_{ji}}{1+\alpha+1} \cdot \frac{ v_{ij}}{1+\alpha+1}  y_{ij}\\
                & = \frac{1}{(2+\alpha)^2} \cdot w_{ji} v_{ij} y_{ij}.
\end{align*}
The first inequality is due to Jensen's inequality applied to the function (convex over the positive quadrant) $(x,y)\mapsto 1/(xy)$. The second is a consequence of the negative correlation of the $x_{ij}$. 
\end{proof}

\subsubsection{High value Customers.}\label{app:cardinality_fs_high_suppliers}
In this section, we focus on the following result: 
\begin{lemma} \label{lemma:high_value_suppliers_cardinality}
Suppose $w_{ji} \geq \alpha$ for all $(i,j) \in \C \times \S$. There exists a randomized polynomial-time algorithm that achieves a $\left(\frac{e-1}{e}\cdot\frac{\alpha}{1+\alpha}\right)$-approximation for \fullystatic~under the two-way constrained model.
\end{lemma}
To show this lemma, we consider the following upper bound of \fullystatic~under this regime.
\begin{equation} \label{eq:concave_K_cardinality}
z_{\sf C} = \max \left\{\;\; \sum_{i=1}^n \frac{\sum_{j=1}^m v_{ij} y_{ij}}{1+\sum_{j=1}^m v_{ij} y_{ij}}:
 \sum_{i=1}^n y_{ij} \leq 1, \ \forall j \in [m], \;\; \sum_{j=1}^m y_{ij}\le K_i,\; \forall i\in[n], \ 
\; \mathbf{y} \in\{0,1\}^{n\times m}\right\}.
\end{equation}
Analogous to the unconstrained case, one can easily show that \eqref{eq:concave_K_cardinality} upper bounds \fullystatic. The constraints in this problem correspond to the intersection of two partition matroids: (a) one over the set of suppliers, $\sum_{i=1}^n y_{ij} \leq 1,$ for all $j \in [m]$; (b) another over the set of customers, $\sum_{j=1}^m y_{ij}\le K_i$ for all $i\in[n]$. Furthermore, observe that the demand function of each customer $f_i(\mathbf{y}_i)=\sum_{j=1}^m v_{ij} y_{ij}/(1+\sum_{j=1}^m v_{ij} y_{ij})$ is monotone, submodular, and is defined only on one part of the partition matroids, specifically, constraints~(a). Given this, we can use the randomized algorithm proposed by~\citep{chekuri2010dependent} (see their Theorem II.4) which, for our problem, guarantees the following:
\begin{lemma}[\citet{chekuri2010dependent}]\label{lemma:rounding_matroid}
        There exists a randomized polynomial-time algorithm for Problem~\eqref{eq:concave_K_cardinality} that achieves a $1-1/e$ approximation factor.
\end{lemma}
Finally, as in Lemma~\ref{thm:high_value_suppliers}, any solution $\mathbf{y}$ to Problem~\eqref{eq:concave_K_cardinality} with objective value {\sf val} will have an objective value of at least $\frac{\alpha}{1+\alpha}\cdot {\sf val}$. Putting this and Lemma~\ref{lemma:rounding_matroid} together, the proof of Lemma~\ref{lemma:high_value_suppliers_cardinality} follows.%
\endproof

\section{Appendix to  Section~\ref{sec:experiments}}\label{app:experiments}

\subsection{Missing Benchmarks}\label{app:experiments_benchmarks}

\subsubsection{Bilinear Formulation for \fullystatic.}\label{app:fs_bilinear}
We introduce a reformulation of Problem~\eqref{eq:main_static} as a disjoint bilinear program, with linear constraints formed by the product of two packing polytopes. Formally,
\begin{align}
\label{eq:bilinear}
z_B= \max &\quad \sum_{i\in[n],j\in[m]} y_{ij}z_{ij} &\\
&\quad 0\le y_{ij}\le v_{ij}\left(1-\sum_{\ell=1}^m y_{i\ell}\right), \quad & \forall i \in [n],\forall j \in [m], \notag \\
&\quad 0\le z_{ij}\le w_{ji}\left(1-\sum_{k=1}^n z_{kj}\right), \quad & \forall i \in [n],\forall j \in [m]. \notag
\end{align}
\begin{theorem}
    Problems \eqref{eq:main_static} and \eqref{eq:bilinear} are equivalent.
\end{theorem}
\begin{proof}[Proof.]
    Consider a solution $\mathbf{x}\in\{0,1\}^{n\times m}$ of Problem~\eqref{eq:main_static}. We construct a feasible solution of Problem~\eqref{eq:bilinear}. For any $i\in[n], j\in[m]$, define
    \begin{align*}
        y_{ij}&=\frac{v_{ij}x_{ij}}{1+\sum_{\ell=1}^m v_{i\ell}x_{i\ell}},&
        z_{ij}&=\frac{w_{ji}x_{ij}}{1+\sum_{k=1}^n w_{jk}x_{kj}}.
    \end{align*}
    It is clear that $(\mathbf{y},\mathbf{z})$ is a solution to Problem~\eqref{eq:bilinear} with same objective value as $x$ for Problem~\eqref{eq:main_static} since $x_{ij}^2 = x_{ij}$. To prove the other direction, notice that Problem~\eqref{eq:bilinear} can be rewritten as
    \begin{align*}
        z_B=\max_{\mathbf{y}\in Y}\left(\max_{\mathbf{z}\in Z}\sum_{i\in[n],j\in[m]} y_{ij}z_{ij}\right)=\max_{\mathbf{z}\in Z}\left(\max_{\mathbf{y}\in Y}\sum_{i\in[n],j\in[m]} y_{ij}z_{ij}\right),
    \end{align*}
    where $Y$ and $Z$ are the corresponding packing constraints for variables $\mathbf{y}$ and $\mathbf{z}$, respectively. The equality follows because of separability of $Y$ and $Z$. Moreover, each inner
     maximization problem  is a linear program. This implies that an optimal solution $(\mathbf{y},\mathbf{z})$ of Problem~\eqref{eq:bilinear} is a vertex in $Y\times Z$. 
     
     Let us now study the optimal vertex. The vector $\mathbf{y}$ has $mn$ coordinates and the system of constraints over $\mathbf{y}$ in Problem~\eqref{eq:bilinear} features $2mn$ constraints that are grouped in pair $y_{ij}\ge 0$ and $y_{ij}\le v_{ij}(1-y_{i1}...-y_{im})$. Observe that $y_{i1}+...+y_{im}$ is always strictly smaller than one, and thus at most one constraint of the pair can be satisfied, except if $v_{ij}=0$, in which case the two constraints are linearly dependent. Thus, a vertex has at least one constraint active per pair, meaning that  for all $i\in[n]$ and for all $j\in[m]$
    \begin{align*}
        y_{ij}=0\text{\quad or \quad} y_{ij}=v_{ij}\left(1-\sum_{\ell=1}^m y_{i\ell}\right).
    \end{align*}
    A similar reasoning can be applied to $\mathbf{z}$. Moreover, if $z_{ij}=0$ or $y_{ij}=0$, then both can be set to $0$ without altering the objective value. Finally, we define $x_{ij}=0$ if $y_{ij}$ and $z_{ij}$ are $0$, otherwise $x_{ij}=1$. We can easily verify that $\mathbf{x}$ is a solution for Problem~\eqref{eq:main_static} with the same objective value as $(\mathbf{y},\mathbf{z})$ for Problem~\eqref{eq:bilinear}.
\end{proof}

\subsubsection{Upper Bound for \onesidedadaptive.}\label{app:ub_oa}
Assume that customers and suppliers follow MNL choice models. Consider the following concave program 
\begin{align}
\label{eq:ub_oa}
z_C= \max &\quad \sum_{j\in\S} \frac{z_j}{1+z_j} &\\
&\quad y_{ij} + \sum_{\ell\in\S} v_{i\ell} y_{i\ell}\le 1, \qquad & \forall i \in\C,\forall j\in\S, \notag \\
&\quad y_{ij}\ge 0, \qquad & \forall i \in\C,\forall j\in\S \notag,\\
&\quad z_{j}=\sum_{i\in\C} v_{ij}w_{ji} y_{ij}\qquad & \forall j \in\S.\notag
\end{align}
This was proposed as an upper bound for \Conesidedstatic~by \citet{torrico21} and we now show that it, in fact, upper bounds  \Conesidedadaptive. 
\begin{lemma}
    The optimal value $z_C$ of Problem~\eqref{eq:ub_oa} upper bounds $\OPT_{\ConesidedA}$.
\end{lemma}
\begin{proof}[Proof.]
For any $S\subseteq\S$ and $C\subseteq\C$, we denote $v_i(S)=\sum_{j\in S} v_{ij}$ and $w_j(C)=\sum_{i\in C}w_{ji}$. Take $\lambda_{j,C},\tau_{i,S}$ a solution to Problem~\eqref{eq:relaxation_onesided_customers}, which is a relaxation of \Conesidedadaptive. For every $i,j$, define
\begin{align*}
    y_{ij}=\sum_{S\ni j}\frac{\tau_{i,S}}{1+v_i(S)}\ge 0.
\end{align*}
We have that
\begin{align*}
    y_{ij}+\sum_{\ell\in\S} v_{i\ell}y_{i\ell}&=\sum_{S\ni j}\frac{\tau_{i,S}}{1+v_i(S)} + \sum_{\ell\in\S}\sum_{S\ni\ell} v_{i\ell}\frac{\tau_{i,S}}{1+v_i(S)}\\
    &=\sum_{S\ni j}\frac{\tau_{i,S}}{1+v_i(S)}+\sum_{S\subseteq\S}\frac{\tau_{i,S}v_i(S)}{1+v_i(S)}\\
    &\le\sum_{S\subseteq\S}\tau_{i,S}\le 1.
\end{align*}
Now, for any $j\in\S$, Jensen's inequality ensures that
\begin{align*}
\sum_{C\subseteq\C}\lambda_{j,C}\frac{w_j(C)}{1+w_j(C)}&\le\frac{\sum_{C\subseteq\C}\lambda_{j,C}\cdot w_j(C)}{1+\sum_{C\subseteq\C}\lambda_{j,C}\cdot w_j(C)}.
\end{align*}
Moreover, we obtain 
\begin{align*}
    \sum_{C\subseteq\C}\lambda_{j,C}\cdot w_j(C)=\sum_{C\subseteq\C}\lambda_{j,C}\sum_{i\in C} w_{ji}=\sum_{i\in\C}w_{ji}\sum_{C\ni i}\lambda_{j,C}&=\sum_{i\in\C}w_{ji}\sum_{S\ni j}\tau_{i,S}\frac{v_{ij}}{1+v_i(S)}\\
    &=\sum_{i\in\C}v_{ij}w_{ji}y_{ij}=z_j,
\end{align*}
where the third equality is due the second constraints in Problem~\eqref{eq:relaxation_onesided_customers}. Thus,
\begin{align*}
\sum_{j\in\S}\sum_{C\subseteq\C}\lambda_{j,C}\cdot \frac{w_j(C)}{1+w_j(C)}\le\sum_{j\in\S}\frac{z_j}{1+z_j}\le z_C, 
\end{align*}
which concludes that the objective value $z_C$ uppers bounds $\OPT_{\ConesidedA}$.
\end{proof}
By symmetry of the problem, it is possible to upper bound $\OPT_{\SonesidedA}$ and, consequently, $\OPT_{\onesidedA}$ by taking the maximum of $z_C$ and $z_S$, where $z_S$ is the objective value of the symmetric problem of Problem~\eqref{eq:ub_oa} (exchanging customers and suppliers).

\subsubsection{Upper Bound for \fullyadaptive.}\label{app:ub_fa}
Assume that customers and suppliers follow MNL choice models. Consider the following linear program
\begin{align}
\label{eq:ub_fa}
z= \max &\quad \sum_{i,j\in\C\times\S} x_{ij} &\\
&\quad x_{ij} \le v_{ij}\cdot\left(1-\sum_{\ell\in\S} x_{i\ell}\right), \qquad & \forall i \in\C,\forall j\in\S, \notag \\
&\quad x_{ij}\leq w_{ji}\cdot\left(1-\sum_{k\in\C} x_{kj}\right), \qquad & \forall i \in\C,\forall j\in\S, \notag \\
&\quad x_{ij}\ge 0, \qquad & \forall i \in\C,\forall j\in\S. \notag
\end{align}
\begin{lemma}
    The optimal value $z$ of Problem~\eqref{eq:ub_fa} upper bounds $\OPT_{\fullyA}$.
\end{lemma}
\begin{proof}[Proof.]
    Let $\pi$ be an optimal policy for \fullyadaptive. We define $x_{ij}=\PP(i\leftrightarrow j)$ over the sample paths of policy $\pi$, where $i\leftrightarrow j$ means that $i$ and $j$ choose each other, i.e., they match. We denote $y_{ij}=\PP(i\rightarrow j)$ and $z_{ij}=\PP(i\leftarrow j)$, i.e., the probability that $i$ chooses $j$ and $j$ chooses $i$, respectively. By inclusion, we have that $x_{ij}\le\min\{y_{ij},z_{ij}\}$. Fix $i\in\C$. Consider the probability space $\Tilde{\Omega}$ of sample paths $\Tilde{\omega}$ truncated just before $i$ makes her choice within the assortment $S_i$ displayed by policy $\pi$. This probability space admits a natural coupling with the non-truncated one $\Omega$. For any $\tilde{\omega}\in\tilde{\Omega}$ and  $j\in\S$, we obtain
    \begin{align*}
        \PP_{\tilde{\omega}}(i\rightarrow j)=\frac{v_{ij}\cdot\one_{\{j\in S_i\}}}{1+v_i(S_i)}
    \end{align*}
    where the probability is taken over sample paths $\omega\in\Omega$ such that the truncation of $\omega$ is $\tilde{\omega}$. We have
    \begin{align*}
        \PP_{\tilde{\omega}}(i\rightarrow j)+v_{ij}\cdot\sum_{\ell\in\S}\PP_{\tilde{\omega}}(i\rightarrow\ell)&=v_{ij}\cdot\frac{\one_{\{j\in S_i\}}+\sum_{\ell\in\S} v_{i\ell}\cdot\one_{\{\ell\in S_i\}} }{1+v_i(S_i)}\le v_{ij}.
    \end{align*}
    By averaging over all $\tilde{\omega}$ in the above inequality, we get
    \begin{align*}
        y_{ij}+v_{ij}\sum_{\ell\in\S} y_{i\ell}\le v_{ij}.
    \end{align*}
    which proves the set of first constraints in \eqref{eq:ub_fa}.
    The other constraints are obtained in a similar fashion. Furthermore, $\pi$ and $\mathbf{x}$ have same objective value for their respective optimization problems.
\end{proof}

\subsection{Additional Tables}\label{app:additional_tables}
We present the running time performance in Table~\ref{table:adaptivetimes} for computing the optimal solutions of \onesidedadaptive\ and \fullyadaptive, as well as their respective relaxation-based upper bounds.

\begin{table}[htpb]
{
    \caption{\bf Average time performance}\label{table:adaptivetimes}
    \centerline{
    \scalebox{0.9}{
      \begin{tabular}{cc|c|c|c}
        \toprule
               $m=n$                  & $\OPT_\onesidedA$ & $\UB_\onesidedA$ & $\OPT_\fullyA$ & $\UB_\fullyA$ \\ 
         \midrule
        2 & 0 & 0.039 & 0.002 & 0.013\\ 
        3 & 0.015 & 0.055 & 0.037 & 0.016  \\ 
        4 & 0.063 & 0.086 & 1.381 & 0.016  \\ 
        5 & 1.101 & 0.129 & 60.727 & 0.016  \\ 
        6 & 23.34 & 0.162 & 3301 & 0.019 \\ 
        7 & 571.0 & 0.251 & - & 0.021 \\ 
        8 & - & 0.335 & - & 0.023  \\ 
        10 & - & 0.455 & - & 0.028 \\ 
        15 & - & 0.934 & - & 0.062 \\ 
        20 & - & 1.553 & - & 0.132 \\ 
        \bottomrule
    \end{tabular}}}
     \vspace{1em}
                \centerline{\begin{minipage}[t]{13cm}
                \footnotesize{\emph{Note.} 
		Average time performance (in seconds) for each adaptive policy and each upper bound. The symbol `-' indicates that the algorithm goes past the limit of 2 hours. We do not report the time performance of $\OPT_{\onesidedS}$ since it is not the main focus of our work.  
		}
                \end{minipage}}
                }
\end{table}


\color{black}
\section{Additional Experiments: Sensitivity Analysis}\label{appendix:sensitivity}
In this appendix, our goal is to empirically analyze the impact of several market features on the algorithms' performance. Following the terminology of Section~\ref{sec:experiments}, we compare the performance of $\alg_\fullyS$, $\alg_\onesidedS$, and $\alg_\onesidedA$. Note that in our paper, $\alg_\onesidedA$ and $\alg_\fullyA$ correspond to the same algorithm, so we do not include the latter. 
In our analysis, we aim to test their performance with respect to two metrics:
\begin{itemize}
    \item the \textit{heterogeneity} of the agents' choice models;
    \item the value of the outside option.
\end{itemize}

{\bf Experimental Setup.} We consider the unconstrained setting with both sides having the same size, i.e., $m=n$. Customers and suppliers make choices according to MNL choice models. Specifically, the value that a customer $i\in\C$ has over a supplier $j\in\S$ is sampled from the uniform distribution on $[0,1]$, i.e., $v_{ij}\sim~U[0,1]$. On the other hand, the value that a supplier $j\in\S$ has over a customer $i\in\C$ is sampled from the exponential distribution of parameter 1, i.e., $w_{ji}\sim \text{Exp}(1)$. To introduce a more complex behavior, for each customer, we choose uniformly at random a factor $a_i\in\{0.5,1,2\}$
that scales all their value parameters. Each of these scaling factors attempt to represent different levels of attraction: pickier (0.5), moderate (1), and less picky (2). 
The same is done for suppliers. Furthermore, each $v_{ij}$ and $w_{ji}$ is randomly set to $0$ independently with probability $\sqrt{2}/2\approx0.707$. This ensures that each pair $(i,j)\in\C\times\S$ consider each other acceptable with probability $0.5$. 
As opposed to Section~\ref{sec:experiments} where the MNL parameters were generated without imposing any structure, in this setup, we generate the parameters in such a way that can capture different levels of heterogeneity. More concretely, to determine the preference values of a given customer (same applies to each supplier), we introduce a parameter $q\in\NN$ that represents the number of different MNL choice models that this customer could follow. Before setting any preference values, $q$ models are sampled independently according to the described distributions above. Then, each customer is assigned to one of the $q$ models independently and uniformly at random. The same applies independently for each supplier and we use the same parameter $q$. On the one hand, when $q=1$, all customers share the exact same choice model; and the same holds for the supplier side. In other words, users on each side are homogeneous. On the other hand, when $q=10$, the resulting instance is heterogeneous enough and hardly indistinguishable from the case when $q=\infty$.\footnote{We did not consider higher values of $q$ since the results did not significantly change.} 
We consider four market sizes: $m=n=20$, $m=n=50$, $m=n=200$, $m=n=500$. All the experiments where performed on a machine with an Intel(R) i5, 4 cores clocked at $2.40$GHz, and $8$GB of RAM. To solve linear and bilinear programs we use Gurobi $11.0.3$ and to solve convex programs we use CVXPY $1.5.3$.

\begin{figure}[htbp]
    \centering
    \begin{subfigure}{0.4\textwidth}
        \centering
        \includegraphics[width=\linewidth]{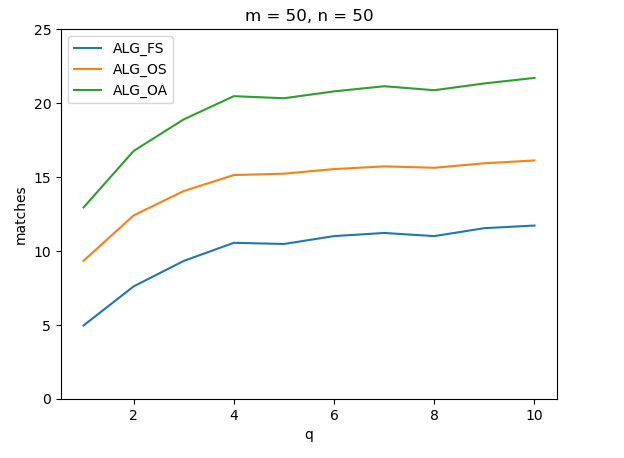}
    \end{subfigure}
    \begin{subfigure}{0.4\textwidth}
        \centering
        \includegraphics[width=\linewidth]{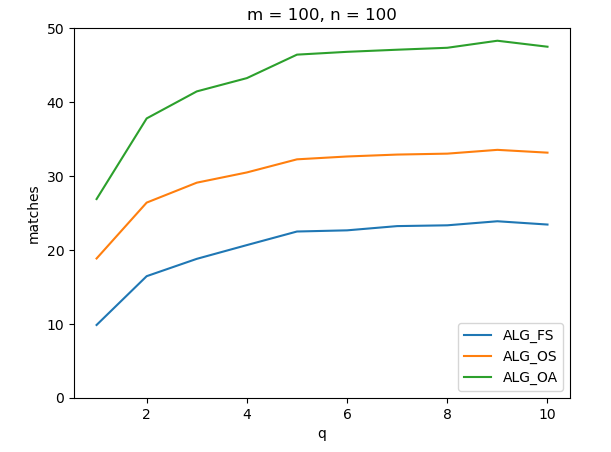}
    \end{subfigure}

    \vspace{0.5em}

    \begin{subfigure}{0.4\textwidth}
        \centering
        \includegraphics[width=\linewidth]{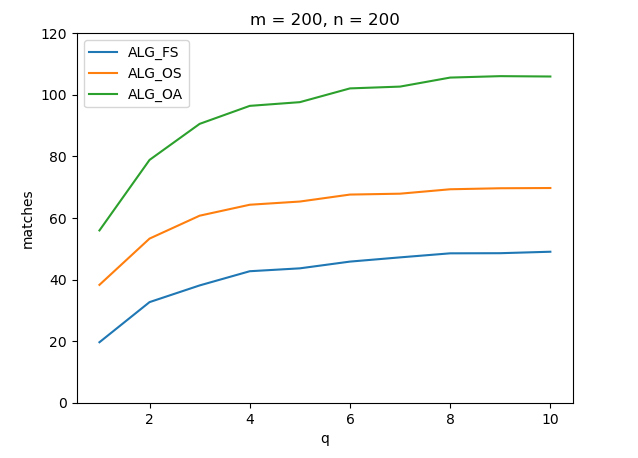}
    \end{subfigure}
    \begin{subfigure}{0.4\textwidth}
        \centering
        \includegraphics[width=\linewidth]{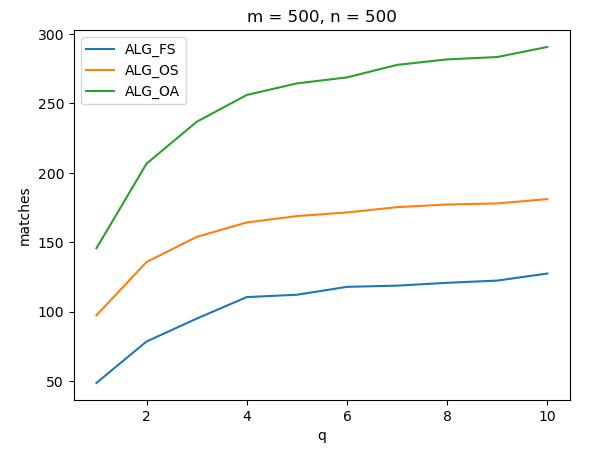}
    \end{subfigure}
    \caption{\textcolor{black}{Impact of the heterogeneity parameter $q$ for different market sizes.}}
    \label{fig:q}
\end{figure}

\vspace{1em}
\noindent{\bf Impact of the Heterogeneity Level.} For this experiment, we fix the outside option value to $v_0=1$. We present our results in Figure~\ref{fig:q} where parameter $q$ ranges from $1$ to $10$. For each value of $q$, we sample between 100 and 200 instances. In particular, to account for the higher variance of the sampling of instances in the low heterogeneity regime, we use a varying number of samples equal to $100+\left\lfloor \frac{100}{q^2}\right\rfloor$. For each value of $q$, we report the average number of matches over all these samples. From Figure~\ref{fig:q}, we observe that all three algorithms' performance improves when the heterogeneity increases. One likely interpretation is that the increased diversity facilitates the pairing of agents by reducing choice congestion. However, the three algorithms preserve their relative performance with respect to each other, where we note that higher levels of adaptivity yields better results.

\begin{figure}[htbp]
    \centering
    \begin{subfigure}{0.4\textwidth}
        \centering
        \includegraphics[width=\linewidth]{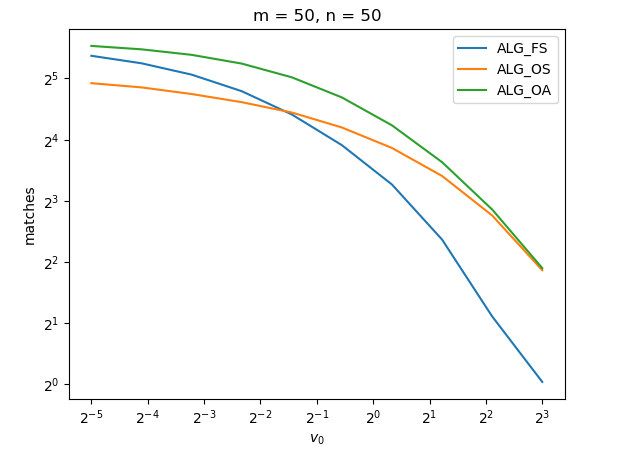}
    \end{subfigure}
    \begin{subfigure}{0.4\textwidth}
        \centering
        \includegraphics[width=\linewidth]{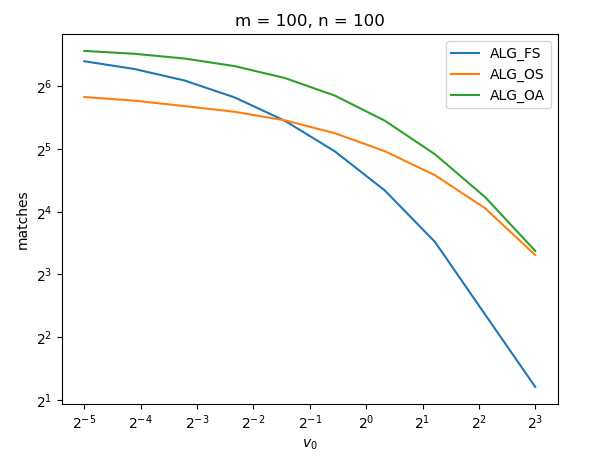}
    \end{subfigure}

    \vspace{0.5em}

    \begin{subfigure}{0.4\textwidth}
        \centering
        \includegraphics[width=\linewidth]{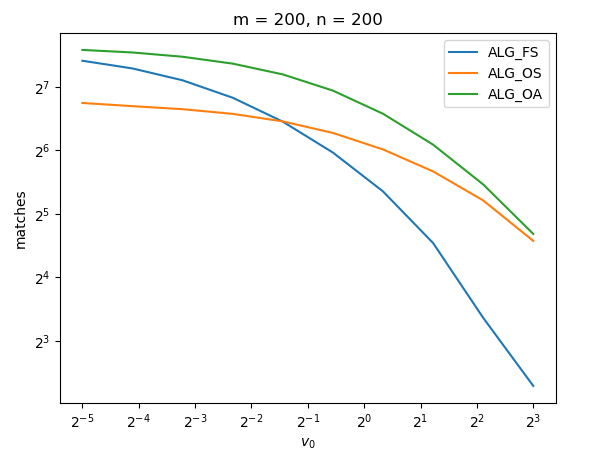}
    \end{subfigure}
    \begin{subfigure}{0.4\textwidth}
        \centering
        \includegraphics[width=\linewidth]{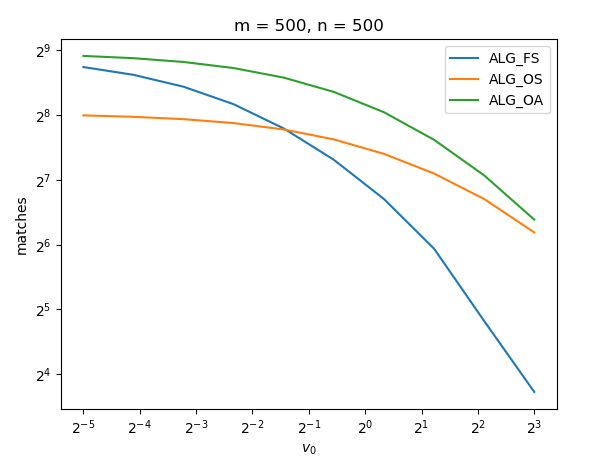}
    \end{subfigure}
    \caption{\textcolor{black}{Impact of the outside option parameter $v_0$ for different market sizes.}}
    \label{fig:v_0}
\end{figure}

\vspace{1em}
\noindent{\bf Impact of the Outside Option Value.} In these experiments, we fix $q=10$ but vary the outside option value $v_0$ from $1/32$ to $8$. To generate each performance curve, we start by generating $50$ instances with $q=10$ then we test all three algorithms for each value $v_0$. 
In Figure~\ref{fig:v_0}, 
we report the results with logarithmic axis to ease the visualization of our results. We first remark that each curve is decreasing, which matches our intuition that making the agents pickier should reduce the number of matches. Furthermore, in the large $v_0$ regime, all curves tend to $0$ as it becomes impossible to get any matches. In the small $v_0$ regime, we remark that $\alg_\fullyS$ and $\alg_\onesidedA$ converge to the same objective value which corresponds to a maximal matching where each edge $(i,j)$ is present if and only if $v_{ij}\cdot w_{ji}>0$. 
On the other hand, we remark that $\alg_\onesidedS$ converges to smaller values than the objectives of $\alg_\fullyS$ and $\alg_\onesidedA$ when $v_0$ is very small. Intuitively, $\OPT_{\onesidedS}$, $\OPT_\fullyS$ and $\OPT_\onesidedA$  should all converge to  the same value in this regime, however  we would like to note that the curves we are plotting correspond to approximate solutions and not the optimal ones. Hence the order of adaptivity is not always necessarily preserved by the algorithms.

\color{black}